\documentclass[preprint,3p,authoryear,times]{elsarticle}

\usepackage[fleqn]{amsmath}
\usepackage[]{amsfonts}
\usepackage[]{amsthm}
\usepackage[]{amssymb}
\usepackage{empheq}
\usepackage{lscape} 

\usepackage[ruled,longend,linesnumbered]{algorithm2e}
\usepackage{graphicx} % Required for inserting images

\usepackage[table]{xcolor}
\usepackage[normalem]{ulem}

\newcommand{\rev}[1]{#1}
\newcommand{\delrev}[1]{}

\usepackage[]{algorithm2e}
\usepackage[english]{babel}
\usepackage[T1]{fontenc}

\usepackage{hyperref} % To include references inside the pdf file
\usepackage{graphicx}
\usepackage{float} % 加载float包
\usepackage{caption,subcaption}
\usepackage{tikz}
\usepackage{pstricks}
\usepackage{xstring}
\usepackage{multirow}
\usepackage{tabularx}
\usepackage{hhline}
\usepackage{url}

\definecolor{darkblue}{rgb}{0,0,0.5}

\hypersetup{colorlinks=true,linkcolor=black,citecolor=darkblue,urlcolor=darkblue} % Set options for the hyperref package

\newcommand{\newref}[2]{\hyperref[#2]{#1~\ref*{#2}}} % hyperref: \newref{labeltext (e.g. table, figure, ...)}{labelname}

\usepackage{longtable}
\usepackage{pdflscape}

\usepackage[toc,page]{appendix}

\usepackage{bbding}
\usepackage{tabularray}
\usepackage{booktabs} % 导入三线表需要的宏包

\begin{document}

\begin{frontmatter}

\title{\fontsize{16pt}{18pt}\selectfont
On the quantification of the spatiotemporal impact of a single discretionary lane change on upstream traffic safety and efficiency%
\tnoteref{titleNote}}

\tnotetext[titleNote]{The paper was published in \textit{Accident Analysis \& Prevention} in 2026.
\url{https://doi.org/10.1016/j.aap.2026.108774}}

% Authors:
\author[1,2]{Kangning Hou}
\author[1,2]{Jia Zou}
\author[1,2]{Fangfang Zheng$^*$}
\author[1,2]{Xiaobo Liu}
\author[3]{Zhengbing He}

%\cortext[cor1]{Corresponding author.}

% Affiliations:
\address[1]{School of Transportation and Logistics, Southwest Jiaotong University, Sichuan, China}

\address[2]{Intelligent Comprehensive Transportation Key Laboratory of Sichuan Province, Southwest Jiaotong University, Sichuan, China}

\address[3]{Faculty of Science and Engineering, University of Nottingham Ningbo China, Zhejiang, China}

% \begin{frontmatter}

% % Title:
% % \title{\fontsize{20pt}{28pt}\selectfont On the quantification of the spatiotemporal impact of a single discretionary lane change on upstream traffic safety and efficiency\footnote[$^\dagger$]{The paper is published by Accident Analysis & Prevention in 2026  \url{https://doi.org/10.1016/j.aap.2026.108774}}}

% \title{\fontsize{16pt}{18pt}\selectfont
% On the quantification of the spatiotemporal impact of a single discretionary lane change on upstream traffic safety and efficiency%
% \textsuperscript{\dagger}%
% \begingroup
% \renewcommand{\thefootnote}{\dagger}
% \footnotetext{The paper was published in \textit{Accident Analysis \& Prevention} in 2026. \url{https://doi.org/10.1016/j.aap.2026.108774}}
% \endgroup
% }

% % Authors:
% \author[1,2]{Kangning Hou}
% \author[1,2]{Jia Zou}
% \author[1,2]{Fangfang Zheng$^*$}
% %\ead{fzheng@swjtu.edu.cn}
% \author[1,2]{Xiaobo Liu}
% \author[3]{Zhengbing He}

% % Affiliations:
% \address[1]{School of Transportation and Logistics, Southwest Jiaotong University, Sichuan, China}

% \address[2]{Intelligent Comprehensive Transportation Key Laboratory of Sichuan Province, Southwest Jiaotong University, Sichuan, China}

% \address[3]{Faculty of Science and Engineering, University of Nottingham Ningbo China, Zhejiang, China}

\begin{sloppypar} 
\begin{abstract}

Lane-changing is a critical driving maneuver, and understanding its effects on traffic safety and efficiency is essential for effective traffic management and optimization. \rev{However, existing studies provide limited means to account for natural traffic dynamics when identifying disturbances associated with lane changes.} Moreover, methods for measuring the spatial extent and duration of the impact of a single discretionary lane change, as well as comprehensive metrics for quantifying its overall spatiotemporal impact, remain underdeveloped. To address these gaps, this study proposes a comprehensive analytical framework to evaluate the spatiotemporal impact of a single discretionary lane change on upstream traffic. The framework identifies affected following vehicles and their impact durations from vehicle trajectories before and after the lane change. A Travel Distance Bias indicator is introduced to measure the motion deviation of following vehicles relative to local reference vehicles, while vehicle-specific pre-event fluctuation envelopes and run-length persistence filters are used to \rev{reduce the misattribution of background traffic fluctuations to the lane-change event}. Based on the affected vehicles and affected intervals, two aggregate indicators, the Total Efficiency Impact Magnitude (TEIM) and the Total Safety Impact Magnitude (TSIM), are developed to quantify efficiency- and safety-related impact magnitudes. Matched no-lane-change controls and method-comparison experiments provide \rev{comparative evidence that the affected-vehicle identification procedure reduces, but does not eliminate, false-positive detections caused by background traffic fluctuations}. Empirical results based on Zen traffic trajectory data from Japan show that the identified impact in the target lane lasts on average 23.8 seconds, affects approximately 5.6 vehicles, and yields average TEIM and TSIM values of -10.8 meters and 4.6\%, respectively, corresponding to average efficiency loss and increased TTC-based unsafe exposure. In the original lane, the identified impact persists for about 25 seconds, affects 5.3 vehicles, and yields average TEIM and TSIM values of 4.7 meters and 2.1\%, respectively, corresponding to average efficiency gain and increased TTC-based unsafe exposure.

\vspace{0.1cm}
\end{abstract}


\end{sloppypar} 

\end{frontmatter}

%%%%%%%%%%%%%%%%%%%%%%%%%%%%%%%%%%%%%%%%%%%%%%%%%%%%%%%%%%%%%
%% CONTENTS
%%%%%%%%%%%%%%%%%%%%%%%%%%%%%%%%%%%%%%%%%%%%%%%%%%%%%%%%%%%%%
%\linenumbers  Running line numbers.

%%%%%%%%%%%%%%%%%%%%%%%%%%%%%%%%%%%%%%%%%%%%%%%%%%%%%%%%%%%%%

\newpage
\begin{sloppypar} 
% \section*{Highlights} % 3 to 5 highlights
% \begin{itemize}
% \item {We evaluate the impact of lane configurations on bus energy consumption.}

% \end{itemize}

%%%%%%%%%%%%%%%%%%%%%%%%%%%%%%%%%%%%%%%%%%%%%%%%%%%%%%%%%%%%%%%%%%%%%%%%%%%%%%%%%%%%%%%%%%%%%%%%%%%%%%%%%%%%%
\section{Introduction}
Lane-changing is a critical component of traffic dynamics, significantly influencing traffic flow efficiency and road safety. Despite substantial advancements in developing algorithms for lane-changing decision making \citep{ali2021clacd,li2022decision,zhang2023learning,shi2024deepad}, trajectory planning \citep{zong2022dynamic,liu2022dynamic,gao2023dual,hou2024cooperative}, and lane-changing prediction \citep{xing2020ensemble,ali2022predicting,huang2024driver}, the magnitude of the effects of a single discretionary lane change on traffic safety and efficiency has not yet been satisfactorily quantified \citep{zheng2014recent,an2023vehicle}. Precise observation and measurement underscore the need for accurate modeling of lane-change impacts to devise effective traffic control strategies and enhance road safety measures. Consequently, exploring quantitative approaches to assess the impact of a lane change is of significant importance. 

Existing studies on lane-change impact quantification can be divided into macroscopic and microscopic approaches. Macroscopic methods evaluate the effects of lane changes from an aggregate traffic dynamics perspective, primarily focusing on traffic efficiency, safety, and fuel consumption \citep{jin2010kinematic,jin2013multi,xiao2014stochastic,feng2015traffic,pan2016modeling,li2017studies,zhu2022flow}. These studies specifically examine the influence of all lane-changing vehicles on the overall traffic flow. For instance, \citet{jin2010kinematic} introduced a location-dependent lane-change intensity variable into the fundamental diagram. Building on this, \citet{jin2013multi} developed a multicommodity kinematic wave (KW) model that treated lane-changing vehicles and non-lane-changing vehicles as separate commodities to explore their impact on traffic capacity. \citet{feng2015traffic} quantified the lane-change impact by calculating the average delay using statistical analysis of simulation data. \citet{li2017studies} applied car-following and lane-change models to assess the impact of lane changes on traffic flow rate, average vehicle speed, traffic safety, and fuel consumption in scenarios involving lane drops and moving bottlenecks. \citet{zhu2022flow} proposed a cooperative merging model designed to minimize the total delay caused by merging maneuvers for all vehicles passing through the ramp merging area. However, these macroscopic methods do not provide insights at the individual vehicle level, which limits their effectiveness in detailed microscopic lane-changing decision studies.

Microscopic approaches, on the other hand, focus on the specific effects of a lane change at the level of individual vehicles, analyzing the complex behaviors and interactions of the lane-changing vehicle and its surrounding vehicles during the lane-changing process. These approaches can be further classified into three categories: instantaneous quantification, fixed spatial-temporal quantification, and dynamic impact quantification. In the category of instantaneous quantification, the Minimizing Overall Braking Induced by Lane Change (MOBIL) model, introduced by \citet{kesting2007general}, utilized the anticipated acceleration of the "following vehicle in the target lane (referred to as TFV)" as an indicator of lane-change impact to guide lane-changing decisions. \citet{xie2019data} considered the states of the lane-changing vehicle’s immediately adjacent vehicles (such as TFV) to be pivotal features influencing lane-changing decisions, proposing a deep belief network-based lane-changing decision model. \citet{hou2023cooperative} developed utility functions that incorporate lane-change impacts in terms of efficiency and comfort for the TFV, aiding vehicles in deciding between maintaining longitudinal movement and executing a lane change. However, these methods primarily focus on the immediate following vehicle and provide only an instantaneous assessment of the lane-change impact. 

In the category of fixed spatial-temporal quantification, \citet{coifman2006impact} introduced the concept of delay induced by a lane change and proposed a method to estimate it between two detectors, considering the distance between two detectors as the spatial extent of the lane-change impact. \citet{yang2019examining} evaluated the lane-change impact on the immediate following vehicle using metrics such as speed change rate, braking timestamp, and time-to-collision (TTC) during the lane-changing duration. \citet{an2019lane} developed an optimization approach to minimize the total deceleration effort of the vehicle directly influenced by the lane-changing maneuver in the target lane within a specified duration, thereby generating an optimal lane-change trajectory. \citet{li2020short} predicted the lane-change impact on crash risks and flow change, employing a fixed spatiotemporal impact region, assuming a fixed number of 5 TFVs affected by a single lane change for 9 seconds, with an equal impact duration for each vehicle. This assumption overlooks the fact that the impact duration may vary among different TFVs. \citet{chen2023gaps,chen2023modeling} measured the lane-change impact within the pre-insertion process (referred to as anticipation) by examining the difference in reaction time of the immediate TFV at the start and end moments of this anticipation process. A critical limitation of these methods is that selecting inappropriate spatial or temporal thresholds can lead to either overestimation or underestimation of the impact magnitude, failing to capture the dynamic progression of the impact.

In the category of dynamic quantification, two notable studies \rev{stand} out \citep{zheng2013effects,he2023impact}. \citet{zheng2013effects} highlighted that driver characteristics, such as reaction time and minimum spacing in Newell’s car-following model, vary over time. By leveraging the relationship between maximum passing rate and response time proposed by \citet{chiabaut2009fundamental}, they calibrated time-dependent reaction times and assessed the lane-change impact on the immediate follower in the target lane. However, this study does not specify a temporal or spatial range for the lane-change impact. \citet{he2023impact} estimated the number of vehicles affected by a lane change and the duration of its effect on each vehicle by examining the evolving pattern of minimum spacing. Nonetheless, this method cannot differentiate between positive and negative lane-change impacts (e.g., travel time reduction vs. increase, or risk mitigation vs. escalation), as it focuses solely on the affected duration or the number of affected vehicles. 

In summary, \rev{existing studies provide limited methodological support for characterizing the temporal and spatial impacts associated with a single lane change} (\hyperref[table1]{Table~\ref{table1}}), leaving two major research gaps. First, \rev{few studies explicitly account for background traffic fluctuations when identifying lane-change-associated disturbances, which may result in such fluctuations being misattributed to the lane-change event}. Second, prior research has seldom jointly characterized the spatial extent and duration of the identified disturbances across the original and target lanes, quantified their cumulative spatiotemporal magnitudes, or distinguished between beneficial and adverse changes. \rev{Addressing these gaps requires a framework that reduces such misattribution while providing a comprehensive and interpretable characterization of lane-change-associated spatiotemporal impacts.} 

\begin{table}[h!]
\centering
\captionsetup{font=small} % 设置标题字体大小为10pt
\caption{Comparative study of dynamic quantification for lane-change impact}\label{table1}
\small
\begin{tabular}{llll}
  \toprule[1.0pt] % 第一行较粗的线
  \textbf{Methodological capabilities} & \citet{zheng2013effects} & \citet{he2023impact} & \textbf{This study}\\
  \midrule[0.5pt] % 中间较细的线
  Quantify impact duration (temporal)   & \ding{55}   & \ding{51}   & \ding{51}   \\
  Quantify number of affected vehicles (spatial)   & \ding{55}    & \ding{51}   & \ding{51}   \\
  Quantify cumulative spatiotemporal impact   & \ding{55}   & \ding{55}  & \ding{51}   \\
  \rev{Reduce misattribution of background traffic fluctuations}  & \ding{55}   & \ding{55}   & \ding{51}  \\
  Distinguish between positive and negative impacts  & \ding{55}   & \ding{55}   & \ding{51}  \\
  \bottomrule[1.0pt] % 最后一行较粗的线
\end{tabular}
\end{table}

To address the above gaps, this paper proposes a trajectory-based methodology to quantify the spatiotemporal impact of a single discretionary lane change on traffic safety and efficiency at the vehicle level. The analysis relies on naturalistic trajectory data before and after the lane change. The main contributions are summarized as follows:

\begin{enumerate}[(1)]
  \item \textbf{Analytical Approach}: We develop a fluctuation-aware analytical approach to identify affected upstream following vehicles and their impact durations, and to estimate the spatial and temporal impact range of a single discretionary lane change on upstream vehicles in both the original and target lanes. To reduce the misattribution of background traffic fluctuations to lane-change impacts, the approach combines local dynamic reference trajectories, vehicle-specific pre-event TDB fluctuation envelopes, and vehicle-specific run-length persistence filters.  
  \item \textbf{Efficiency and Safety Metrics}: Two new metrics are introduced—the Total Efficiency Impact Magnitude (TEIM) and the Total Safety Impact Magnitude (TSIM)—to quantify the continuous spatiotemporal effects of lane changes on upstream traffic. Both metrics incorporate baseline-excess corrections \rev{to reduce the contribution of background traffic fluctuations to the estimated lane-change impact magnitudes.}
  \item \textbf{Validation Design}: We introduce a validation design based on matched no-lane-change controls to assess the reliability of the proposed quantification framework in the absence of ground-truth affected-vehicle labels. Specifically, the design evaluates the trade-off between the affected-vehicle identification rate in real lane-change events and the false positive rate in matched no-lane-change pseudo-events.  
\end{enumerate}

This study provides new insights into the spatiotemporal impact of a single discretionary lane change from both safety and efficiency \rev{perspectives}, highlighting clear differences between the target and original lanes.

The remainder of this paper is organized as follows: \hyperref[Section2]{Section~\ref{Section2}} outlines the methodology used to quantify the spatiotemporal impact range of lane changes at the individual vehicle level. \hyperref[Section3]{Section~\ref{Section3}} introduces the proposed traffic efficiency and safety metrics. \hyperref[Section4]{Section~\ref{Section4}} describes the trajectory data and the methods for obtaining lane-change data and estimating the lane-change start time. \hyperref[Section5]{Section~\ref{Section5}} provides an in-depth analysis of the lane-change impact through data examples. \hyperref[Section6]{Section~\ref{Section6}} reports matched no-lane-change control validation, key-parameter sensitivity analyses, method comparisons, and robustness assessment. Finally, \hyperref[Section7]{Section~\ref{Section7}} concludes with main findings and directions for future research.

%%%%%%%%%%%%%%%%%%%%%%%%%%%%%%%%%%%%%%%%%%%%%%%%%%%%%%%%%%%%%%%%%%%%%%%%%%%%%%%%%%%%%%%%%%%%%%%%%%%%%
\section{Methodology for quantifying the spatiotemporal impact range}\label{Section2}
\begin{figure}[!htbp]
  \centering
  \includegraphics[width=\textwidth]{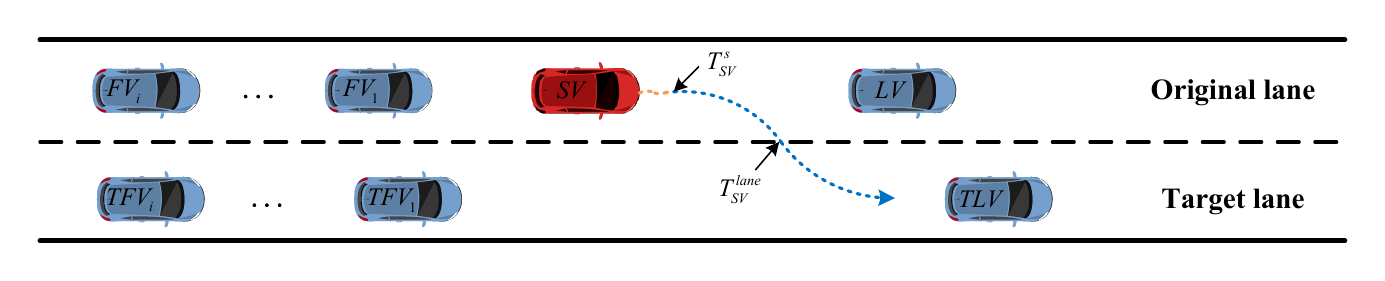}
  \captionsetup{font={small}}
  \caption{Illustration of a typical lane-changing scenario}
  \label{Fig1}
\end{figure}

This section presents a trajectory-based analytical framework designed to quantify the spatiotemporal impact range of a single lane-change by analyzing vehicle trajectories before and after the maneuver, specifically identifying the number of affected vehicles and their impact durations. \hyperref[Fig1]{Fig.~\ref{Fig1}} illustrates a typical lane-changing scenario, where the \rev{subject} vehicle (\textit{SV}) intends to change lanes. \textit{TLV} and ${TFV}_i$ denote the leading and the $i^{th}$ following vehicles of \textit{SV} in the target lane, while \textit{LV} and ${FV}_i$ represent the leading and following vehicles of \textit{SV} in the original lane. $T_{SV}^{lane}$ denotes the moment when \textit{SV} crosses the lane marking, and $T_{SV}^s$ denotes the moment when \textit{SV} initiates the lane-changing process. 

The proposed framework comprises four critical components: 1) Demarcation time identification: Using a kinematic wave theory-based method, the demarcation time $T_{\ell,i}^s$ is determined to separate potentially affected and unaffected data for each following vehicle $F_{\ell,i}$ in both the target lane and the original lane. 2) Efficiency metric design: Travel Distance Bias (TDB) is used to capture relative travel distance deviations of vehicles. 3) Judgment criteria: Two criteria are applied to determine whether $F_{\ell,i}$ is affected by the \textit{SV}’s lane change and to estimate the lane-specific impact duration. 4) Spatiotemporal impact range: the total number of affected vehicles and overall impact durations are obtained. For clarity, the main variables and lane-specific notation used in the methodology are summarized in \hyperref[table2]{Table~\ref{table2}}.

\begin{table}[!h]
  \centering
  \captionsetup{font=small} % 设置标题字体大小为10pt
  \caption{Notation List}\label{table2}
  \renewcommand{\arraystretch}{1.1} % 调整行高为1.5倍
  \resizebox{\textwidth}{!}{ % 调整表格大小以适应页面宽度
  %\small % 设置表格文本字体大小为 10pt
  \begin{tabular}{p{2cm}p{14.5cm}}
    \toprule[1.0pt] % 第一行较粗的线
    \textbf{Notation} & \textbf{Description} \\
    \midrule[0.5pt] % 中间较细的线
    $T_{SV}^s$ & Time when \textit{SV} initiates the lane-change execution process \\
    $\ell \in \{T,O\}$ & Lane side, with $T$ denoting the target lane and $O$ denoting the original lane \\
    $F_{\ell,i}$ & The $i^{th}$ following vehicle in lane $\ell$; $F_{T,i}=TFV_i$ and $F_{O,i}=FV_i$ \\
    $L_\ell$ & Local leading reference vehicle in lane $\ell$; $L_T=TLV$ and $L_O=LV$ \\
    $T_{\ell,i}^s$ & Demarcation time separating potentially unaffected and affected data of $F_{\ell,i}$ \\
    ${TDB}_{\ell,i}(k)$ & Travel distance bias at the $k^{th}$ interval, representing the delay-distance deviation of $F_{\ell,i}$ relative to $L_\ell$ \\
    ${TDB}_{\ell,i}^\ast$ & Threshold of TDB for $F_{\ell,i}$ \\
    $\eta$ & Multiplier controlling the width of the pre-event TDB fluctuation envelope \\
    $\Delta t$ & Time interval \\
    $\Theta_{\ell,i}(k)$ & Indicator of potentially affected status of $F_{\ell,i}$ at the $k^{th}$ time interval \\
    $\mathbf{\Theta}_{\ell,i}^f, \mathbf{\Theta}_{\ell,i}^r$ & Sets of $\Theta_{\ell,i}(k)$ in potentially unaffected and affected data segments, respectively \\
    $n_{\ell,i}^f, n_{\ell,i}^r$ & Number of $\Delta t$ intervals in potentially unaffected and affected data segments, respectively \\
    $\Omega_{\ell,i}^{f\ast}$ & Maximum number of consecutive occurrences of $\Theta_{\ell,i}(k)=1$ in the potentially unaffected segment for $F_{\ell,i}$ \\
    $\mathbf{c}_{\ell,i,q}^{f}, \mathbf{c}_{\ell,i,q}^{r}$ & Sets of $k$ values corresponding to $\Omega_{\ell,i,q}^f, \Omega_{\ell,i,q}^r$, respectively \\
    $\mathbf{K}_{\ell,i}^{A}$ & Set of estimated affected time intervals for $F_{\ell,i}$ \\
    $\Upsilon_{\ell,i}$ & Estimated affected status of $F_{\ell,i}$ \\
    $T_{\ell,i}^A$ & Impact duration of the \textit{SV}'s lane change on $F_{\ell,i}$ \\
    $m_{\mathrm{stop}}$ & Number of consecutive unaffected following vehicles required to stop the upstream impact search \\
    ${CTDB}_{\ell,i}(k)$ & Corrected Travel Distance Bias by integrating ${TDB}_{\ell,i}(k)$ with correction factor $\delta_{\ell,i}(k)$ \\
    $\delta_{\ell,i}(k)$ & Correction factor for CTDB \\
    $w_{\ell,i}^A$ & Efficiency impact magnitude of the \textit{SV}'s lane change on $F_{\ell,i}$, i.e., the sum of ${CTDB}_{\ell,i}(k)$ during impact \\
    $N_\ell^A$ & Total number of affected following vehicles in lane $\ell$ \\
    $W_{\ell,effi}^A$ & Total efficiency impact magnitude (TEIM) of the \textit{SV}'s lane change on all affected $F_{\ell,i}$ \\
    $T_\ell^A$ & Overall duration of the \textit{SV}'s lane change on lane $\ell$ \\
    $r_{\ell,i}^A$ & Net safety impact magnitude of the \textit{SV}'s lane change on $F_{\ell,i}$ \\
    $R_{\ell,unsafe}^A$ & Total safety impact magnitude (TSIM) of the \textit{SV}'s lane change on all affected $F_{\ell,i}$ \\
    \bottomrule[1.0pt] % 最后一行较粗的线
  \end{tabular}
  }
\end{table}

\subsection{Determining \texorpdfstring{$T_{\ell,i}^s$}{Tli\_s} using kinematic wave theory}\label{Section2.1} 
To illustrate the methodology, consider a lane-specific following vehicle $F_{\ell,i}$. Since the exact moment when $F_{\ell,i}$ is first affected by the \textit{SV}'s lane change cannot be precisely determined, a demarcation time $T_{\ell,i}^s$ is introduced to distinguish potentially unaffected from affected data for $F_{\ell,i}$. Specifically, data prior to $T_{\ell,i}^s$ are considered potentially unaffected data, while data after $T_{\ell,i}^s$ are treated as potentially affected data. The accuracy of the impact assessment depends on aligning $T_{\ell,i}^s$ as closely as possible with the actual onset of disturbance for $F_{\ell,i}$.

\begin{figure}[h]
  \centering
  \includegraphics[width=0.65\textwidth]{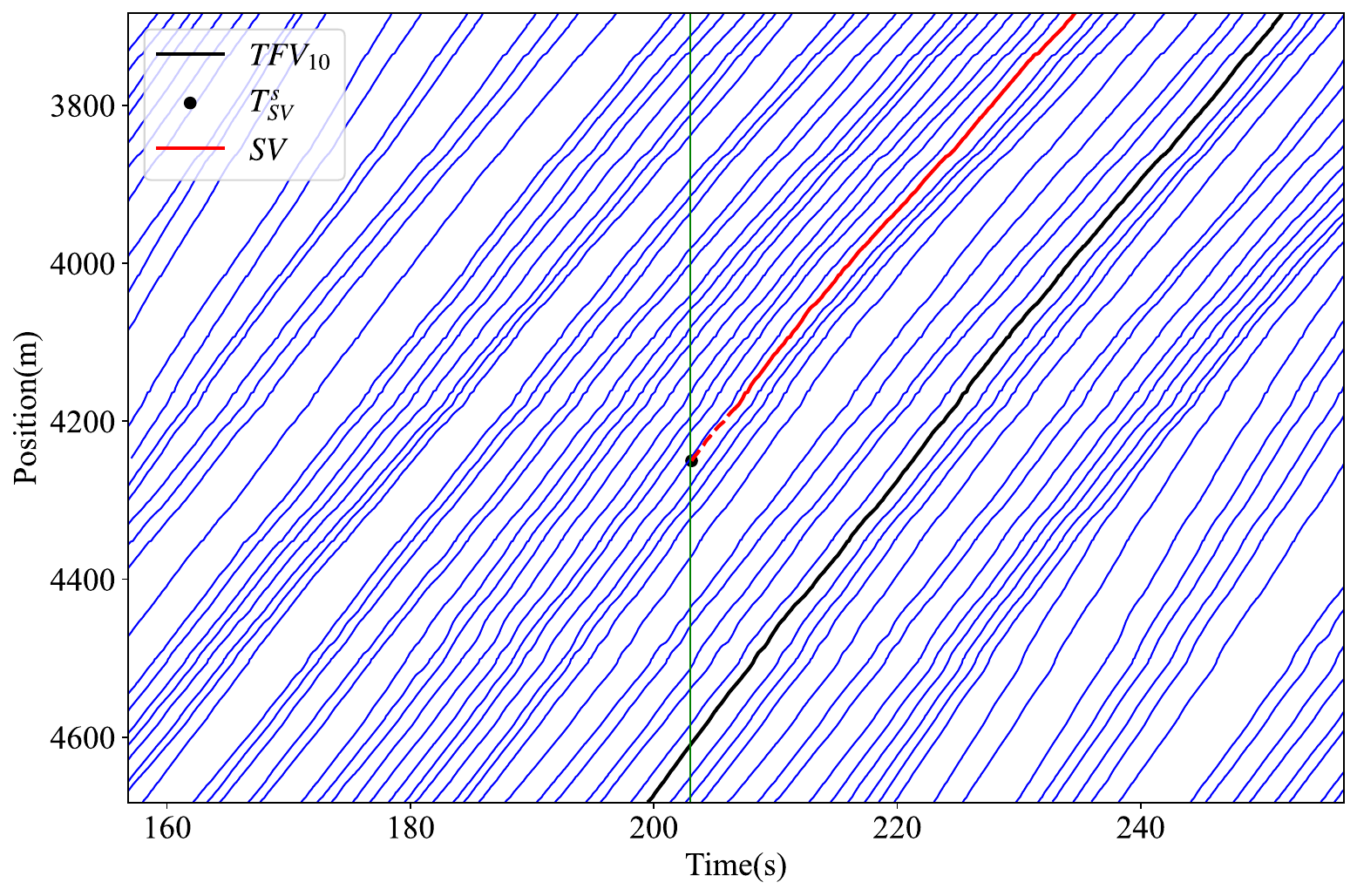}
  \captionsetup{font={small}}
  \caption{Trajectories from an extracted lane-changing instance}
  \label{Fig2}
\end{figure}

\hyperref[Fig2]{Fig.~\ref{Fig2}} illustrates an extracted lane-changing instance. The red and black lines represent the trajectory curves of the \textit{SV} and ${TFV}_{10}$, respectively, with the black dot marking the moment when the \textit{SV} initiates the lane-changing process ($T_{SV}^s$). If $T_{SV}^s$ were directly used as $T_{\ell,i}^s$, vehicles farther from the \textit{SV} (e.g. ${TFV}_{10}$) would have little or no potentially unaffected data available, reducing the accuracy of impact estimation. 

To address this issue, a refined method based on kinematic wave theory is proposed. The physical disturbance differs by lane: in the target lane, \textit{SV} insertion compresses the available gap and may induce deceleration responses, whereas in the original lane, \textit{SV} departure releases space and may induce acceleration and gap-closing responses. In both cases, the response propagates upstream through the lane-specific car-following chain. The propagation time between two adjacent vehicles ($t_w$) can be calculated as \citep{holland1998generalised,qin2021lighthill}:
\begin{equation}\label{Eq1}
t_w = \frac{h}{v - c_w} = \left[\frac{dV}{dh}\right]^{-1} = \left[V'(h)\right]^{-1}
\end{equation}
where $h$ and $v$ represent the equilibrium spacing and speed, respectively, and $c_w$ is the wave speed relative to the road, which is a positive value. $V(h)$ is the equilibrium speed-spacing function, i.e., $v=V(h)$. Following \citet{qin2021lighthill}, we can obtain:
\begin{equation}\label{Eq2}
[V'(h)]^{-1} = h^v
\end{equation}
where $h^v$ represents the derivative of $h$ with respect to $v$. Thus, the demarcation time $T_{\ell,i}^s$ is calculated as:
\begin{equation}\label{Eq3}
T_{\ell,i}^s = T_{SV}^s + \sum_{j=1}^{i} t_{\ell,w,j}, \quad i \geq 1
\end{equation}
where $t_{\ell,w,j}$ represents the time required for the kinematic wave to propagate between $F_{\ell,j}$ and its preceding vehicle $F_{\ell,j-1}$. According to \hyperref[Eq1]{Eqs.~(\ref{Eq1})} and \hyperref[Eq2]{(\ref{Eq2})}, $t_{\ell,w,j}=h_{\ell,j}^v$, which can be calibrated using real trajectory data based on the lower order Newell model \citep{newell2002simplified}:
\begin{equation}\label{Eq4}
x_{\ell,j}(t + \tau_{\ell,j}) = x_{\ell,j-1}(t) - d_{\ell,j}
\end{equation}
where $\tau_{\ell,j}$ and $d_{\ell,j}$ represent the reaction time and minimum space of $F_{\ell,j}$, respectively. $x_{\ell,j}$ and $x_{\ell,j-1}$ denote the positions of $F_{\ell,j}$ and $F_{\ell,j-1}$, respectively. In the Newell model, the relationship between $h_{\ell,j}$ and $v$ at equilibrium is given by:
\begin{equation}\label{Eq5}
h_{\ell,j} = d_{\ell,j} + v\tau_{\ell,j}
\end{equation}

From \hyperref[Eq5]{Eq.~(\ref{Eq5})}, $h_{\ell,j}^v$ is obtained as $\tau_{\ell,j}$; thus, $t_{\ell,w,j}=\tau_{\ell,j}$, and $T_{\ell,i}^s=T_{SV}^s+\sum_{j=1}^{i} \tau_{\ell,j}$. The parameters $\tau_{\ell,j}$ and $d_{\ell,j}$ are calibrated to minimize the discrepancy between the observed trajectories and those predicted by the Newell model, using the preceding vehicle's position data as input. The objective function for the calibration is the sum of squared errors (SSE):
\begin{equation}\label{Eq6}
(\tau_{\ell,j}^*, d_{\ell,j}^*) = \underset{(\tau_{\ell,j}, d_{\ell,j})}{\text{argmin}} \sum_{t \in T_{\text{all}}^{\ell,j}} (x_{\ell,j}(t) - x_{\ell,j}'(t))^2
\end{equation}
where $x_{\ell,j}'(t)$ represents the predicted position of $F_{\ell,j}$ at time $t$, and $T_{\text{all}}^{\ell,j}$ represents the total time that $F_{\ell,j}$ is within the selected trajectory set. Parameters were calibrated independently for each following vehicle in each lane, with bounds set as $\tau_{\ell,j}$ $\in [0.1, 5]$ and $d_{\ell,j}$ $\in [0.1, 10]$.

\subsection{Designing the Travel Distance Bias (TDB) metric}\label{Section2.2} 
\begin{figure}[h]
  \centering
  \includegraphics[width=\textwidth]{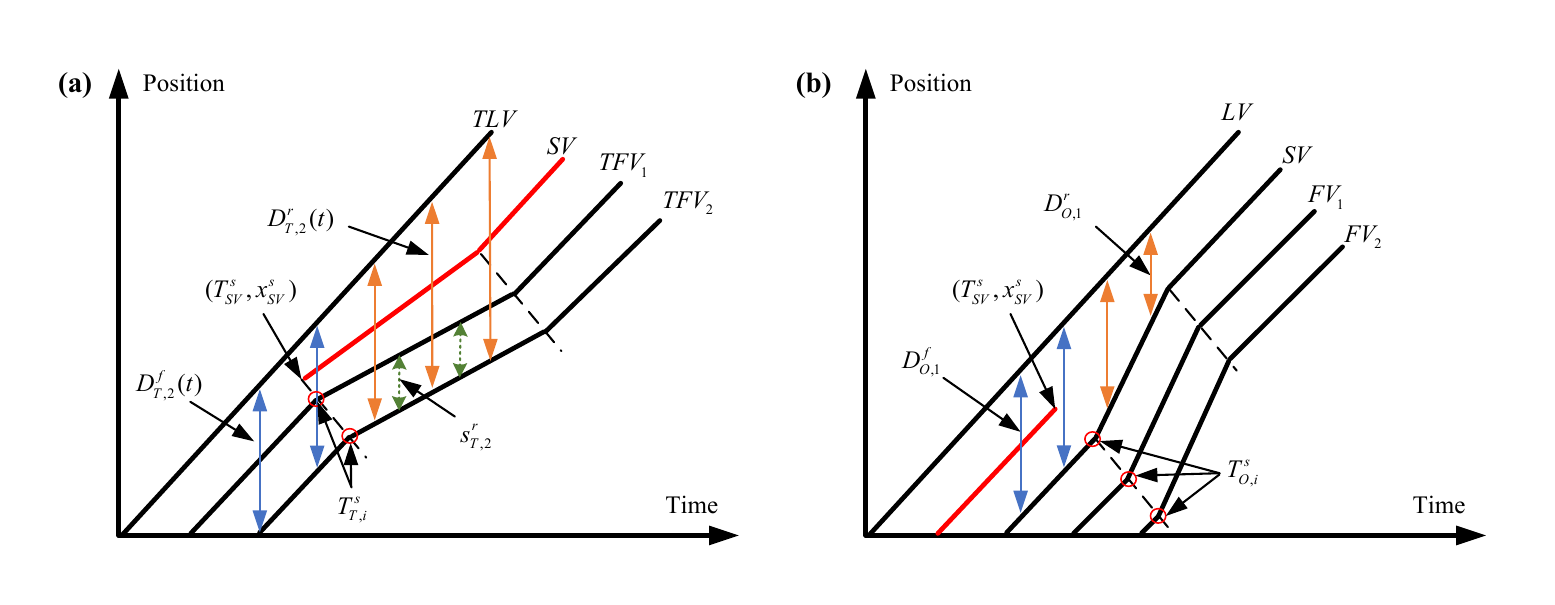}
  \captionsetup{font={small}}
  \caption{Schematic diagram of vehicle trajectories before and after a lane change: (a) \textit{TFV}s and (b) \textit{FV}s}
  \label{Fig3}
\end{figure}

\hyperref[Fig3]{Fig.~\ref{Fig3}} illustrates vehicle trajectories in the target and original lanes before and after a lane change. In \hyperref[Fig3]{Fig.~\ref{Fig3}a}, $D_{T,2}^f$, $D_{T,2}^r$, and $s_{T,2}^r$ denote, respectively, the distance between ${TFV}_2$ and \textit{TLV} before $T_{T,2}^s$, the distance between ${TFV}_2$ and $TLV$ after $T_{T,2}^s$, and the distance between ${TFV}_2$ and ${TFV}_1$ after $T_{T,2}^s$. The blue, orange, and green curves correspond to $D_{T,2}^f$, $D_{T,2}^r$, and $s_{T,2}^f$. It is observed that $D_{T,2}^r$ changes continuously after $T_{T,2}^s$, whereas $s_{T,2}^r$ remains relatively stable. This stability aligns with the Newell model, which suggests that the trajectory of a following vehicle tends to mimic that of its leader. Thus, in a car-following scenario, the spacing between a vehicle and its immediate predecessor fluctuates only slightly around the equilibrium spacing, producing minimal distance variation over short intervals. Accordingly, we focus on the relative distance variation between $F_{\ell,i}$ and the lane-specific local leading reference $L_\ell$, rather than the short-term spacing between $F_{\ell,i}$ and its immediate predecessor, to construct an efficiency metric that reflects the lane-change impact. This metric, termed travel distance bias (TDB), is defined as:
\begin{equation}\label{Eq7}
{TDB}_{\ell,i}(k) = \int_{T_{\ell,i}^s+(k-1)\Delta t}^{T_{\ell,i}^s+k\Delta t} \Delta V_{\ell,i}(t) \, dt
\end{equation}
\begin{equation}\label{Eq8}
\Delta V_{\ell,i}(t) = v_{F_{\ell,i}}(t) - v_{L_\ell}(t)
\end{equation}
where $v_{F_{\ell,i}}(t)$ and $v_{L_\ell}(t)$ represent the speeds of $F_{\ell,i}$ and $L_\ell$ at time $t$, respectively. The TDB measures the delay distance of $F_{\ell,i}$ compared with its lane-specific local leading reference $L_\ell$ over each time interval $\Delta t$. A negative value indicates that $F_{\ell,i}$ travels a shorter distance than $L_\ell$, while a positive value suggests a longer distance. Thus, Eq.~(\ref{Eq8}) uses \textit{TLV} as the local dynamic reference for $TFV_i$ in the target lane and \textit{LV} as the local dynamic reference for $FV_i$ in the original lane.

At the macroscopic level, traditional traffic delay is defined as the difference between actual and ideal travel times, with the ideal case assuming constant travel speed (i.e., free flow speed or non-affected speed) \citep{xu2020bi,kodupuganti2023facilities}. If TDB were defined in this manner, $v_{L_\ell}$ in \hyperref[Eq8]{Eq.~(\ref{Eq8})} would be considered constant, serving as an ideal reference. 

At the microscopic level, however, the speed of $L_\ell$ may vary due to traffic conditions ahead, traffic oscillations, or local congestion propagation. Therefore, $L_\ell$ is used as a moving local benchmark that contains part of the background traffic fluctuations also experienced by the upstream following vehicle $F_{\ell,i}$. By comparing the traveled distance of $F_{\ell,i}$ with that of $L_\ell$, TDB partly cancels the motion changes shared with local traffic evolution and retains the excess deviation relative to this local benchmark. This design reduces the misattribution of speed variations shared with the local leading vehicle and estimates the excess disturbance associated with the \textit{SV}'s lane change. Defining TDB in the traditional way \rev{could increase the misattribution} of these fluctuations to the \textit{SV}'s influence.

\hyperref[Fig3]{Fig.~\ref{Fig3}b} shows a comparable case for \textit{FV}s in the original lane, where \textit{FV}s might accelerate post $T_{O,i}^s$ to close the gap with \textit{LV}, suggesting that the \textit{SV}'s lane change could potentially enhance efficiency for \textit{FV}s. Thus, the same signed TDB structure is used for comparable accounting in both lanes, while the sign is interpreted according to the lane-specific physical mechanism.

\subsection{Establishing judgment criteria for determining affected time of \texorpdfstring{$F_{\ell,i}$}{Fli}}\label{Section2.3}

This subsection outlines the procedure for determining the affected time of $F_{\ell,i}$. In the potentially unaffected pre-demarcation segment, $F_{\ell,i}$ generally follows the lane-specific local reference $L_\ell$ with only normal car-following fluctuations. Nevertheless, minor fluctuations may still occur, causing $\Delta V_{\ell,i}(k)$ to deviate from zero. Consequently, ${TDB}_{\ell,i}(k)$ may also deviate from zero in these segments, partly due to inherent traffic noise in trajectory data. To determine whether $F_{\ell,i}$ is affected, it is necessary to analyze ${TDB}_{\ell,i}(k)$ in the potentially unaffected data segments and establish a threshold, denoted as ${TDB}_{\ell,i}^\ast$.

Given that ${TDB}_{\ell,i}(k)$ can be positive or negative when unaffected by a lane change, directly aggregating these values would bias the threshold estimation. To address this, positive and negative ${TDB}_{\ell,i}(k)$ values are separated in the potentially unaffected data segment as follows: 
\begin{equation}\label{Eq9}
\mathbf{TDB}_{\ell,i}^- = \left\{ TDB_{\ell,i}(k) \,|\, TDB_{\ell,i}(k) < 0, \, k \in [1, n_{\ell,i}^f] \right\}
\end{equation}
\begin{equation}\label{Eq10}
\mathbf{TDB}_{\ell,i}^+ = \left\{ TDB_{\ell,i}(k) \,|\, TDB_{\ell,i}(k) \geq 0, \, k \in [1, n_{\ell,i}^f] \right\}
\end{equation}
\begin{equation}\label{Eq11}
n_{\ell,i}^f=\left \lfloor \frac{T_{\ell,i}^s-T_{\ell,i}^{lb}}{\Delta t} \right \rfloor
\end{equation}
where $T_{\ell,i}^{lb}$ represents the lower time boundary of the extracted trajectory data for $F_{\ell,i}$. $n_{\ell,i}^f$ represents the total number of intervals within [$T_{\ell,i}^{lb}$, $T_{\ell,i}^s$]. The operator $\left\lfloor\cdot\right\rfloor $ represents the floor function.

The means ($\mu_{\ell,i}^+$, $\mu_{\ell,i}^-$) and standard deviations ($\sigma_{\ell,i}^+$, $\sigma_{\ell,i}^-$) of these sets are calculated as:
\begin{equation}\label{Eq12}
\mu_{\ell,i}^+ = \frac{1}{m_{\ell,i}} \sum_{{TDB}_{\ell,i}(k) \in \mathbf{TDB}_{\ell,i}^+} {TDB}_{\ell,i}(k)
\end{equation}\label{Eq13}
\begin{equation}
\sigma_{\ell,i}^+ = \sqrt{\frac{1}{m_{\ell,i}}\sum_{{TDB}_{\ell,i}(k) \in \mathbf{TDB}_{\ell,i}^+} {({TDB}_{\ell,i}(k) - \mu_{\ell,i}^+)}^2}
\end{equation}
\begin{equation}\label{Eq14}
\mu_{\ell,i}^- = \frac{1}{(n_{\ell,i}^f-m_{\ell,i})} \sum_{{TDB}_{\ell,i}(k) \in \mathbf{TDB}_{\ell,i}^-} {TDB}_{\ell,i}(k)
\end{equation}
\begin{equation}\label{Eq15}
\sigma_{\ell,i}^- = \sqrt{\frac{1}{n_{\ell,i}^f-m_{\ell,i}}\sum_{{TDB}_{\ell,i}(k) \in \mathbf{TDB}_{\ell,i}^-} {({TDB}_{\ell,i}(k) - \mu_{\ell,i}^-)}^2}
\end{equation}
\begin{equation}\label{Eq16}
m_{\ell,i} = \left| \mathbf{TDB}_{\ell,i}^+ \right|
\end{equation}

We construct a lane- and vehicle-specific pre-event TDB fluctuation envelope from the potentially unaffected segment. This envelope is denoted by the TDB threshold (${TDB}_{\ell,i}^{\ast}$) and is defined as:
\begin{equation}\label{Eq17}
{TDB}_{\ell,i}^* = 
\begin{cases}
[\mu_{\ell,i}^+ - \eta\sigma_{\ell,i}^+,\mu_{\ell,i}^++\eta\sigma_{\ell,i}^+ ], & \text{if } TDB_{\ell,i} (k) \ge 0 \\
[\mu_{\ell,i}^--\eta\sigma_{\ell,i}^-,\mu_{\ell,i}^-+\eta\sigma_{\ell,i}^- ], & \text{if } TDB_{\ell,i} (k) < 0
\end{cases}
\end{equation}
where $\eta$ controls the envelope width. This envelope represents the normal amplitude range of TDB under comparable local traffic and trajectory-noise conditions. Accordingly, a time interval is labeled as potentially affected only when its TDB exceeds this pre-event envelope.

The potentially affected status of $F_{\ell,i}$ at the $k^{th}$ time interval is defined as:
\begin{equation}\label{Eq18}
\Theta_{\ell,i}(k) = 
\begin{cases}
1, & \text{if } TDB_{\ell,i} (k) \notin TDB_{\ell,i}^* \\
0, & \text{otherwise}
\end{cases}
\end{equation}
where $\Theta_{\ell,i}$ is a binary variable, with 1 and 0 representing potentially affected and unaffected states, respectively. However, because TDB is computed discretely, $\Theta_{\ell,i}(k)$ may oscillate rapidly between 0 and 1, which does not always reflect the actual driving behavior. Furthermore, disturbances unrelated to the \textit{SV}'s lane change may occasionally induce large fluctuations in $F_{\ell,i}$'s speed, resulting in unexpected occurrences of $\Theta_{\ell,i}\left(k\right)=1$ within unaffected segments. Thus, relying solely on the TDB threshold to determine whether $F_{\ell,i}$ is affected is insufficient. 

To overcome this, we further construct a lane- and vehicle-specific run-length persistence filter, a discriminant index $\Omega_{\ell,i}^{f\ast}$, defined as the maximum number of consecutive $\Theta_{\ell,i}(k)=1$ occurrences in the potentially unaffected segment:
\begin{equation}\label{Eq19}
\Omega_{\ell,i}^{f*} = \max (\Omega_{\ell,i,1}^f, \Omega_{\ell,i,2}^f, \ldots, \Omega_{\ell,i,q}^f, \ldots), \quad q = 1, 2, 3, \ldots
\end{equation}
where $\Omega_{\ell,i,q}^f$ represents the length of the $q^{th}$ consecutive sequence of $\Theta_{\ell,i}(k)=1$ of $F_{\ell,i}$ in the potentially unaffected data segment. The potentially unaffected segment is used as a vehicle-specific empirical reference for the normal persistence of abnormal TDB states. This filter is not a formal parametric statistical test or an optimal theoretical threshold; rather, it is a nonparametric, vehicle-specific persistence rule that uses the pre-demarcation segment to reduce short-lived abnormal sequences comparable to pre-event fluctuations. The rule is that $F_{\ell,i}$ is considered affected only when the number of consecutive occurrences of $\Theta_{\ell,i}\left(k\right)=1$ exceeds $\Omega_{\ell,i}^{f\ast}$. Thus, a post-demarcation abnormal sequence is retained only when it persists longer than the longest abnormal sequence already observed before $T_{\ell,i}^{s}$; otherwise, it is treated as indistinguishable from background traffic fluctuation or trajectory noise. Accordingly, in both unaffected and affected data segments, any sequence of $\Theta_{\ell,i}\left(k\right)=1$ shorter than $\Omega_{\ell,i}^{f\ast}$ is reclassified as $\Theta_{\ell,i}\left(k\right)=0$. This adjustment guarantees that all values in the potentially unaffected data segment remain $\Theta_{\ell,i}\left(k\right)=0$. The pseudocode for calculating $\Omega_{\ell,i}^{f\ast}$ and adjusting $\Theta_{\ell,i}^f\left(k\right)$ is presented in \hyperref[Algorithm1]{Algorithm~\ref{Algorithm1}}.
\vspace{0.3cm}

\begin{algorithm}[H]
\caption{Calculation of $\Omega_{\ell,i}^{f\ast}$ and adjustment of $\Theta_{\ell,i}^f(k)$ by applying $\Omega_{\ell,i}^{f\ast}$}
\label{Algorithm1}
    \KwIn{$\mathbf{\Theta}_{\ell,i}^f$}
    \KwOut{$\Omega_{\ell,i}^{f\ast}$, $\mathbf{\Theta}_{\ell,i}^f$}
    \Begin{
        $\mathbf{K}_{\ell,i}^f = \{k \mid\mathbf{\Theta}_{\ell,i}^f(k) == 1, k \in [1, n_{\ell,i}^f]\}$\;
        $z_1 = 1$, $q = 1$\;
        \While{$z_1 \leq \text{length}(\mathbf{K}_{\ell,i}^f)$}{
            $z_2 = 1$\;
            \While{$z_2 \leq \text{length}(\mathbf{K}_{\ell,i}^f) - z_1$ \textbf{and} $\mathbf{K}_{\ell,i}^f(z_1) + z_2 == \mathbf{K}_{\ell,i}^f(z_1 + z_2)$}{
                $z_2 = z_2 + 1$\;
            }
            \If{$z_2 \geq 1$}{
                $\mathbf{c} _{\ell,i,q}^{f} = \mathbf{K}_{\ell,i}^f(z_1 : 1 : z_1 + z_2 - 1)$\;
                $\Omega_{\ell,i,q}^f = \textit{length}(\mathbf{c} _{\ell,i,q}^{f})$\;
                $q = q + 1$\;
            }
            $z_1 = z_1 + z_2$\;
        }
        $\Omega_{\ell,i}^{f\ast} = \max(\Omega_{\ell,i,1}^f, \Omega_{\ell,i,2}^f, \ldots, \Omega_{\ell,i,q}^f, \ldots)$\;
        $\mathbf{C} _{\ell,i}^{f} = \{\mathbf{c} _{\ell,i,1}^{f}, \cdots, \mathbf{c} _{\ell,i,q}^{f}, \cdots\}$\;
        \For {each $\mathbf{c} _{\ell,i,q}^{f} \subseteq \mathbf{C} _{\ell,i}^{f}$}{
            \If{$\Omega_{\ell,i,q}^f \leq \Omega_{\ell,i}^{f\ast}$}{
                \For {each $k \in \mathbf{c} _{\ell,i,q}^{f}$}{
                    $\mathbf{\Theta}_{\ell,i}^f(k) = 0$\;
                }
            }
        }
        $\mathbf{\Theta}_{\ell,i}^f = \{\mathbf{\Theta}_{\ell,i}^f(k) \mid k \in [1, n_{\ell,i}^f]\}$\;
    }
\end{algorithm}
\vspace{0.3cm}

We employ a method analogous to \hyperref[Algorithm1]{Algorithm~\ref{Algorithm1}} to process the data following time $T_{\ell,i}^s$, referred to as the potentially affected data segment:
\begin{equation}\label{Eq20}
\Omega_{\ell,i}^r = \left\{ \Omega_{\ell,i,1}^r, \Omega_{\ell,i,2}^r, \ldots, \Omega_{\ell,i,q}^r, \ldots \right\}, \quad q = 1,2,3,\ldots
\end{equation}
where \(\Omega_{\ell,i,q}^r\) represents the number of consecutive occurrences of \(\Theta_{\ell,i}(k) = 1\) in the $q^{th}$ sequence of \(F_{\ell,i}\) within the potentially affected data segment.

Based on the threshold (\(\Omega_{\ell,i}^{f\ast}\)), the estimated affected time intervals \(\mathbf{t}_{\ell,\mathbf{i},\mathbf{q}}^\mathbf{A}\) of \(F_{\ell,i}\) can be determined as:
\begin{equation}\label{Eq21}
\mathbf{t}_{\ell,\mathbf{i},\mathbf{q}}^\mathbf{A} = 
\begin{cases}
    \mathbf{c} _{\ell,i,q}^{r}, & \text{if } \Omega_{\ell,i,q}^r > \Omega_{\ell,i}^{f\ast} \\
    \emptyset, & \text{otherwise}
\end{cases}
\end{equation}
where $\mathbf{c} _{\ell,i,q}^{r}$ represents the set of $k$ values corresponding to $\Omega_{\ell,i,q}^r$. The detailed calculation of $\mathbf{c} _{\ell,i,q}^{r}$ can be found in \hyperref[Algorithm1]{Algorithm~\ref{Algorithm1}}. The union of all affected time intervals is given by: 
\begin{equation}\label{Eq22}
\mathbf{K}_{\ell,i}^A = \bigcup \mathbf{t}_{\ell,\mathbf{i},\mathbf{q}}^\mathbf{A}
\end{equation}

The criterion for determining whether $F_{\ell,i}$ is affected by the lane change is defined as:
\begin{equation}\label{Eq23}
\Upsilon_{\ell,i} = 
\begin{cases}
    1, & \text{if } \mathbf{K}_{\ell,i}^A \neq \emptyset \\
    0, & \text{otherwise}
\end{cases}
\end{equation}
where $\Upsilon_{\ell,i}$ is a binary variable: a value of 1 indicates that $F_{\ell,i}$ is affected, whereas a value of 0 indicates that it is unaffected.

The impact duration of the \textit{SV}'s lane change on $F_{\ell,i}$ is calculated as:
\begin{equation}\label{Eq24}
    T_{\ell,i}^A = t_{\ell,i}^e - t_{\ell,i}^s
\end{equation}
with 
\begin{equation}\label{Eq25}
    t_{\ell,i}^s = T_{\ell,i}^s + \left( \min{\left( \mathbf{K}_{\ell,i}^A \right)} - 1 \right) \cdot \Delta t
\end{equation}
\begin{equation}\label{Eq26}
    t_{\ell,i}^e = T_{\ell,i}^s + \max{\left(\mathbf{K}_{\ell,i}^A \right)} \cdot \Delta t
\end{equation}
where $t_{\ell,i}^s$ and $t_{\ell,i}^e$ represent the starting and ending times of the affected period of $F_{\ell,i}$, respectively. Specifically, $T_{\ell,i}^A$ represents the duration between the first and last affected moments of $F_{\ell,i}$. \hyperref[table3]{Table~\ref{table3}} provides a numerical example illustrating the process of determining $\Upsilon_{\ell,i}$ and $T_{\ell,i}^A$.

\begin{table}[h!]
  \centering
  \renewcommand{\arraystretch}{1.1} % Adjust row height
  \captionsetup{font=small} % 设置标题字体大小为10pt
  \caption{Numerical Example of Determining $\Upsilon_{\ell,i}$ and $T_{\ell,i}^A$}\label{table3}
  \begin{tabular}{p{2.5cm}|p{5.8cm}|p{5.8cm}}
    \hline
    \textbf{Variables} & \textbf{Potentially unaffected segment} & \textbf{Potentially affected segment} \\
    \hline
    \rowcolor{gray!25} % 第2行底纹颜色
    $\mathbf{\Theta}_{\ell,i}$ & \multicolumn{2}{c}{$\mathbf{\Theta}_{\ell,i}=\{0,1,0,1,1,0,0,1,1,1,0,1,1,0\}$} \\
    \hline
    $n_{\ell,i}^f, n_{\ell,i}^r$ & $n_{\ell,i}^f=6$ & $n_{\ell,i}^r=8$ \\
    \hline
    $\mathbf{\Theta}_{\ell,i}^f, \mathbf{\Theta}_{\ell,i}^r$ & $\mathbf{\Theta}_{\ell,i}^f=\{0,1,0,1,1,0\}$ & $\mathbf{\Theta}_{\ell,i}^r=\{0,1,1,1,0,1,1,0\}$ \\
    \hline
    $\Omega_{\ell,i,q}^f, \Omega_{\ell,i,q}^r$ & $\Omega_{\ell,i,1}^f=1, \Omega_{\ell,i,2}^f=2$ & $\Omega_{\ell,i,1}^r=3, \Omega_{\ell,i,2}^r=2$ \\
    \hline
    $\mathbf{c}_{\ell,i,q}^f, \mathbf{c}_{\ell,i,q}^r$ & $\mathbf{c}_{\ell,i,1}^f=\{2\}, \mathbf{c}_{\ell,i,2}^f=\{4,5\}$ & $\mathbf{c}_{\ell,i,1}^r=\{8,9,10\}, \mathbf{c}_{\ell,i,2}^r=\{12,13\}$ \\
    \hline
    \rowcolor{gray!25} % 第2行底纹颜色
    $\Omega_{\ell,i}^{f\ast}$ & \multicolumn{2}{c}{2} \\
    \hline
    $\mathbf{t}_{\ell,\mathbf{i},\mathbf{q}}^\mathbf{A}$ & - & $\mathbf{t}_{\ell,\mathbf{i},\mathbf{1}}^\mathbf{A}=\{8,9,10\}, \mathbf{t}_{\ell,\mathbf{i},\mathbf{2}}^\mathbf{A}=\emptyset$ \\
    \hline
    $\mathbf{K}_{\ell,i}^A$ & - & $\{8,9,10\}$ \\
    \hline
    \rowcolor{gray!25} % 第2行底纹颜色
    $\Upsilon_{\ell,i}$ & \multicolumn{2}{c}{1} \\
    \hline
    \rowcolor{gray!25} % 第2行底纹颜色
    $T_{\ell,i}^A$ & \multicolumn{2}{c}{$3\Delta t$} \\
    \hline
  \end{tabular}
\end{table}

\subsection{Determining the number of affected vehicles and the overall impact duration}\label{Section2.4} 
To mitigate random effects, we introduce a stopping parameter $m_{\mathrm{stop}}$, which denotes the number of consecutive upstream following vehicles in lane $\ell$ that must remain unaffected before the upstream impact search is stopped. In the main analysis, $m_{\mathrm{stop}}=2$, corresponding to the assumption that two consecutive unaffected following vehicles indicate that the lane-specific disturbance generated by the lane change has been absorbed or dissipated. Accordingly, the total number of affected following vehicles in lane $\ell$, denoted as $N_\ell^A$, is determined by:
\begin{equation}\label{Eq27}
N_\ell^A =
\begin{cases}
 i_\ell^{\mathrm{stop}}-1, & \text{if } i_\ell^{\mathrm{stop}} \text{ exists},\\
 N_\ell, & \text{otherwise}
\end{cases}
\end{equation}
with
\begin{equation}\label{Eq28}
i_\ell^{\mathrm{stop}} =
\min\left\{i \mid \sum_{j=0}^{m_{\mathrm{stop}}-1}\Upsilon_{\ell,i+j}=0,\;1\le i\le N_\ell-m_{\mathrm{stop}}+1 \right\}
\end{equation}
where $N_\ell$ represents the number of following vehicles in lane $\ell$ within the extracted trajectory data.

Let $t_{\ell,\text{last}}$ denote the difference between the affected end time of the last affected following vehicle in lane $\ell$ ($t_{\ell,E}$) and the affected start time of the first affected following vehicle in lane $\ell$ ($t_{\ell,S}$), given by:
\begin{equation}\label{Eq29}
t_{\ell,\text{last}} = t_{\ell,E} - t_{\ell,S}
\end{equation}
with
\begin{equation}\label{Eq30}
t_{\ell,E}=T_{\ell,N_\ell^A}^{s}+  \max \left ( \mathbf{K}_{\ell,N_\ell^A}^A \right ) \cdot \Delta t
\end{equation}
\begin{equation}\label{Eq31}
t_{\ell,S} = T_{\ell,initial}^s + \min \left ( \mathbf{K}_{\ell,initial}^A \right )\cdot \Delta t
\end{equation}
where $T_{\ell,N_\ell^A}^s$ and $T_{\ell,initial}^s$ represent the demarcation times distinguishing between potentially unaffected and affected data of the last affected $F_{\ell,N_\ell^A}$ and the first affected $F_{\ell,initial}$. $\mathbf{K}_{\ell,N_\ell^A}^A$ and $\mathbf{K}_{\ell,initial}^A$ represent the corresponding sets of affected intervals, respectively.

Since $t_{\ell,\text{last}}$ may differ from the maximum $T_{\ell,i}^A$, the overall duration of the \textit{SV}'s lane change on lane $\ell$ is calculated as:
\begin{equation}\label{Eq32}
T_\ell^A = \max{\left(t_{\ell,\text{last}},\max\limits_{1\le i\le N_\ell^A}{(T_{\ell,i}^A)}\right)}
\end{equation}

\section{Metrics of traffic efficiency and safety}\label{Section3}
Building on the spatiotemporal impact ranges identified in \hyperref[Section2]{Section~\ref{Section2}}, this section develops quantitative metrics to evaluate the overall impact magnitude of a single discretionary lane change on upstream traffic from both efficiency and safety perspectives. To reduce the contribution of traffic fluctuations unrelated to the lane-change event, the proposed metrics incorporate baseline-excess corrections, so that the estimated efficiency and safety magnitudes can better reflect the excess changes associated with the lane-changing maneuver.

\subsection{Traffic efficiency metric}\label{Section3.1} 
To evaluate the efficiency impact, we introduce the Corrected Travel Distance Bias (CTDB), which extends the previously developed TDB by incorporating a correction factor $\delta_{\ell,i}(k)$. This factor reduces the contribution of normal TDB fluctuations unrelated to the \textit{SV}'s lane change—an aspect overlooked in earlier methods based on evolving reaction time or minimum spacing calibration \citep{zheng2013effects,he2023impact}.

The CTDB for $F_{\ell,i}$ at time step $k$ is defined as:
\begin{equation}\label{Eq33}
{CTDB}_{\ell,i}(k) = {TDB}_{\ell,i}(k) - \delta_{\ell,i}(k)
\end{equation}
where the correction factor $\delta_{\ell,i}(k)$ is defined as:
\begin{equation}\label{Eq34}
\delta_{\ell,i}(k) = 
\begin{cases}
\mu_{\ell,i}^- - \eta\sigma_{\ell,i}^- & \text{if } TDB_{\ell,i}(k) < \mu_{\ell,i}^- - \eta\sigma_{\ell,i}^- \\
\mu_{\ell,i}^- + \eta\sigma_{\ell,i}^- & \text{if } \mu_{\ell,i}^- + \eta\sigma_{\ell,i}^- < TDB_{\ell,i}(k) < 0 \\
\mu_{\ell,i}^+ - \eta\sigma_{\ell,i}^+ & \text{if } 0 < TDB_{\ell,i}(k) < \mu_{\ell,i}^+ - \eta\sigma_{\ell,i}^+ \\
\mu_{\ell,i}^+ + \eta\sigma_{\ell,i}^+ & \text{if } TDB_{\ell,i}(k) > \mu_{\ell,i}^+ + \eta\sigma_{\ell,i}^+ \\
TDB_{\ell,i}(k) & \text{otherwise}
\end{cases}
\end{equation}
\indent This correction can be viewed as a baseline-excess CTDB correction. When ${TDB}_{\ell,i}(k)$ lies within the normal pre-event envelope, $\delta_{\ell,i}(k)={TDB}_{\ell,i}(k)$ and therefore ${CTDB}_{\ell,i}(k)=0$; when ${TDB}_{\ell,i}(k)$ exceeds the envelope, only the excess beyond the corresponding boundary is accumulated.

The efficiency impact magnitude (EIM) of the \textit{SV}'s lane change on $F_{\ell,i}$ is obtained by summing ${CTDB}_{\ell,i}(k)$ over the impact duration:
\begin{equation}\label{Eq35}
w_{\ell,i}^A = \sum_{k \in \mathbf{K}_{\ell,i}^A} { {CTDB}_{\ell,i}(k) }
\end{equation}

The total efficiency impact magnitude (TEIM) of the \textit{SV}'s lane change on all affected $F_{\ell,i}$ in lane $\ell$ is then calculated as:
\begin{equation}\label{Eq36}
W_{\ell,effi}^A = \sum_{i=1}^{N_\ell^A} w_{\ell,i}^A \cdot \Upsilon_{\ell,i}
\end{equation}

\subsection{Traffic safety metric}\label{Section3.2} 
To assess safety impacts, we extend beyond traditional time-to-collision (TTC) indicators to capture the continuous spatiotemporal influence of lane changes on upstream vehicles. While TTC measures the instantaneous risk between a leading vehicle and a following vehicle, it does not account for cumulative or longer-term effects across the traffic stream. Therefore, we integrate TTC with the overall unsafe rate metric proposed in \citet{hou2024enhancing}, which quantifies the proportion of time vehicles spend in unsafe states. 

The unsafe proportion for each $F_{\ell,i}$ during its affected period $\left[t_{\ell,i}^s,t_{\ell,i}^e\right]$ is defined as:
\begin{equation}\label{Eq37}
r_{\ell,i}^{LC} = 
\frac{
\displaystyle \sum_{k = t_{\ell,i}^s / \Delta t}^{t_{\ell,i}^e / \Delta t} \delta_{\ell,i,k}
}{
\left( t_{\ell,i}^e - t_{\ell,i}^s \right) / \Delta t
},
\quad
\delta_{\ell,i,k} =
\begin{cases}
1, & 0 < {\rm TTC}_{\ell,i}(k)  \le {\rm TTC}^\ast, \\
0, & \text{otherwise.}
\end{cases}
\end{equation}

\begin{equation}\label{Eq38}
{\rm TTC}_{\ell,i}(k) =
\begin{cases}
\dfrac{x_{\ell,i-1}(k) - x_{\ell,i}(k) - \mathcal{L}_{\ell,i-1}}{v_{\ell,i}(k) - v_{\ell,i-1}(k)}, & v_{\ell,i}(k) > v_{\ell,i-1}(k), \\[10pt]
\infty, & \text{otherwise.}
\end{cases}
\end{equation}
where $x_{\ell,i}(k)$ and $v_{\ell,i}(k)$ are the position and speed of $F_{\ell,i}$, and $\mathcal{L}_{\ell,i-1}$ is the length of the preceding vehicle $F_{\ell,i-1}$. A larger TTC indicates a safer condition. Although TTC thresholds of 1.5--3 s are commonly used in safety analyses \citep{li2017evaluationcacc,yao2023platoon}, small TTC values can be sparse in non-crash or normal-driving contexts \citep{rahman2018longitudinal,rahman2019safety}. Therefore, the TTC threshold was set to 4.5 s in this study. This choice is consistent with \citet{lu2021performance}, who reported that a TTC threshold could be higher than 4 s and identified 4.5 s in their calibration results.

Because $F_{\ell,i}$ may also experience unsafe conditions caused by natural traffic oscillations (unrelated to the \textit{SV}'s maneuver), a baseline unsafe rate is calculated for each $F_{\ell,i}$ within an equal-duration pre-event baseline window immediately before the affected period:
\begin{equation}\label{Eq39}
r_{\ell,i}^{ULC} =
\frac{
\displaystyle \sum_{k = (2t_{\ell,i}^s - t_{\ell,i}^e)/\Delta t}^{t_{\ell,i}^s / \Delta t} \delta_{\ell,i,k}
}{
\left( t_{\ell,i}^e - t_{\ell,i}^s \right) / \Delta t
}
\end{equation}

The net safety impact magnitude (SIM) associated with the \textit{SV}'s lane change is then derived as:
\begin{equation}\label{Eq40}
r_{\ell,i}^A=r_{\ell,i}^{LC}-r_{\ell,i}^{ULC}
\end{equation}
\indent Thus, SIM represents the net change in TTC-based unsafe exposure relative to the pre-event baseline, rather than the absolute unsafe proportion during the affected period.

A positive $r_{\ell,i}^A$ indicates an increase in risk, while a negative value denotes a risk reduction. Finally, the total safety impact magnitude (TSIM) across all affected $F_{\ell,i}$ in lane $\ell$ is calculated as:
\begin{equation}\label{Eq41}
R_{\ell,unsafe}^A=\sum_{i=1}^{N_\ell^A}{r_{\ell,i}^A\cdot\Upsilon_{\ell,i}}
\end{equation}

Overall, the proposed framework does not rely on a single mechanism \rev{to reduce the misattribution of background traffic fluctuations}. Instead, it mitigates their influence through several connected mechanisms. First, the TDB formulation in Eqs.~(\ref{Eq7})--(\ref{Eq8}) uses the lane-specific local leading reference $L_\ell$ (i.e., \textit{TLV} or \textit{LV}) as a local dynamic reference, which reduces the misattribution of common local speed fluctuations to the subject vehicle's lane change. Second, the pre-event TDB envelope in Eq.~(\ref{Eq17}) establishes a vehicle-specific range of normal TDB variation before the potential disturbance. Third, the run-length filter in Eq.~(\ref{Eq19}) removes short abnormal sequences that may arise from trajectory noise or natural traffic fluctuations. Fourth, the CTDB correction in Eqs.~(\ref{Eq33})--(\ref{Eq34}) accumulates only the excess TDB beyond the baseline fluctuation envelope when computing EIM and TEIM. Finally, the safety metric subtracts the pre-event unsafe-time proportion in Eqs.~(\ref{Eq39})--(\ref{Eq40}) so that SIM and TSIM represent net changes in TTC-based unsafe exposure. Together, these mechanisms are intended to reduce the misattribution of background traffic oscillations to lane-change-associated disturbances and to improve the interpretability of the identified impact range and aggregate impact metrics.

%%%%%%%%%%%%%%%%%%%%%%%%%%%%%%%%%%%%%%%%%%%%%%%%%%%%%%%%%%%%%%%%%%%%%%%%%%%%%%%%%%%%%%%%%%%%%%%%%%%%%
\section{Data Preparation}\label{Section4}
\subsection{The Zen Traffic Data}\label{Section4.1} 
For a comprehensive analysis of the temporal and spatial impact of a single lane change on upstream traffic, access to extensive data is essential. Commonly used traffic trajectory datasets such as HighD and NGSIM are constrained by their limited road lengths—420 meters for HighD, and 640 meters for US101 and 500 meters for I80 in NGSIM, respectively. These constraints limit their effectiveness for in-depth analysis. To overcome this limitation, this study utilizes the Zen Traffic Data (ZTD)\footnote{Zen Traffic Data (ZTD) is an open-source dataset available at \url{https://zen-traffic-data.net/english/}.} from Japan, which is known for its extensive coverage and is thus more suitable for examining the comprehensive effects of lane changes. 

\begin{figure}[h!]
  \centering
  \includegraphics[width=0.9\textwidth]{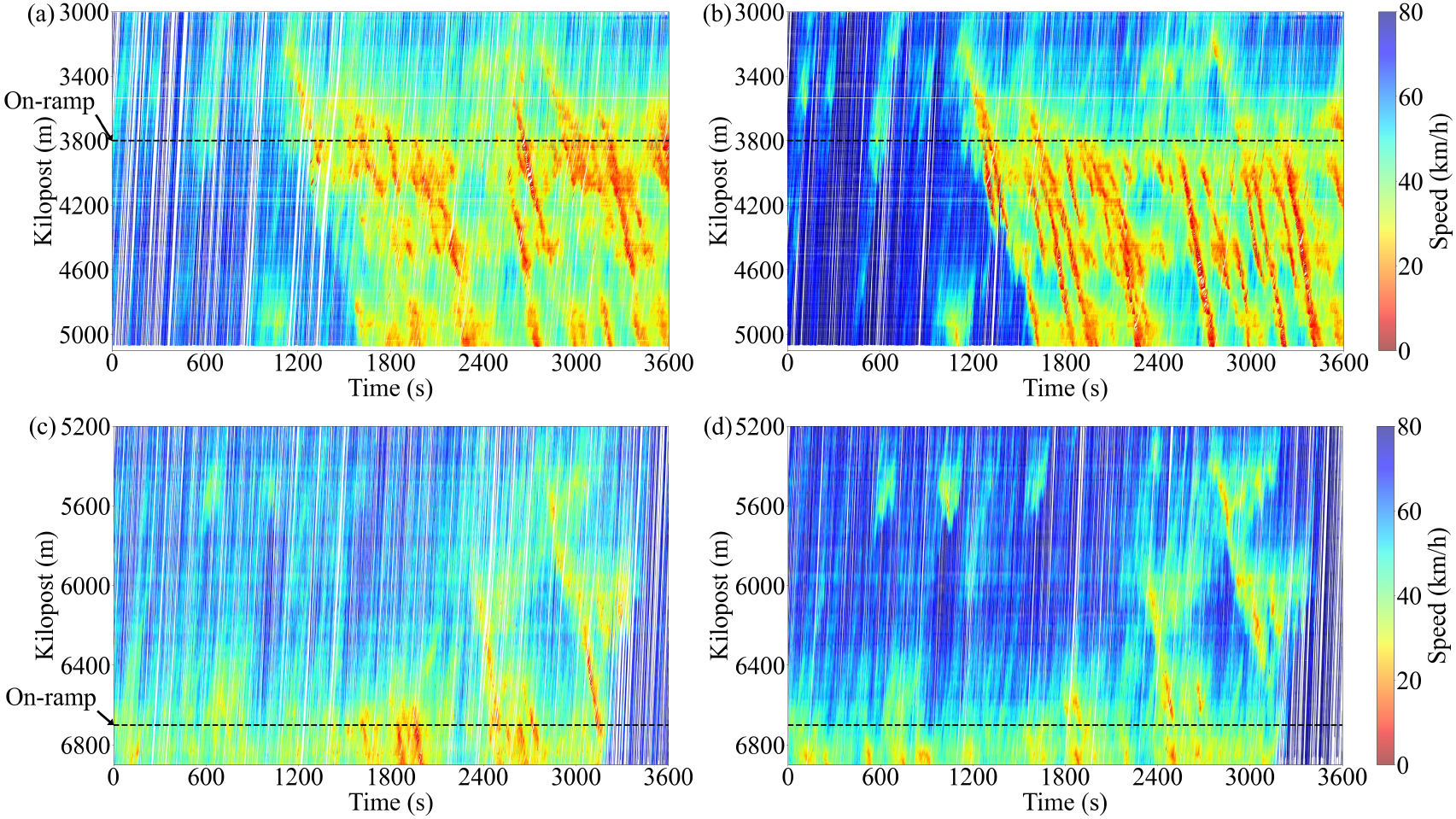}
  \captionsetup{font={small}}
  \caption{Vehicle trajectories of different routes: (a) Ikeda Route \#11 (F001) driving lane; (b) Ikeda Route \#11 (F001) passing lane; (c) Wangan Route \#4 (F001) driving lane; (d) Wangan Route \#4 (F001) passing lane (from \url{https://zen-traffic-data.net/english/outline/})}
  \label{Fig4}
\end{figure}

The ZTD dataset consists of vehicle travel data collected from three road segments in Japan, specifically Routes \#4, \#11, and \#13. Each segment encompasses two main lanes (a driving lane and a passing lane) and one ramp lane. As illustrated in \hyperref[Fig4]{Fig.~\ref{Fig4}}, the data comprises both free-flow traffic and congested traffic, enabling a comprehensive analysis of the lane-change impact across various traffic conditions. The driving lane typically serves slower-moving traffic, while the passing lane accommodates faster-moving vehicles. The data collection durations were as follows: 5 hours for each of the first two segments and 6 hours for the third segment, with recordings taken at 0.1-second intervals. The lengths of the segments for Routes \#4, \#11 and \#13 are 2 km, 1.6km and 2.7 km, respectively. Key data attributes include vehicle ID, datetime, vehicle type, velocity, laneID, kilopost (distance from the downstream endpoint of the expressway route), and geographic coordinates (latitude and longitude). It is important to note that traffic in Japan moves on the left, and the kilopost readings decrease as the vehicle travels along the route.

To enhance the reliability of the speed within the dataset, we apply a filtering method as described by \citet{montanino2015trajectory}. This method helps to smooth the speed profiles of the trajectories, thereby reducing the impact of random noise and improving the accuracy of our analysis.

\subsection{Lane-changing data extraction}\label{Section4.2} 
To quantify the impact of the \textit{SV}'s lane change on upstream traffic, lane-change instances involving the SV and its surrounding vehicles were extracted based on the following criteria:  

\textbf{Criterion 1}: Data from all vehicles in all lanes are collected within a 50-second window before and after $T_{SV}^{lane}$ and within a 500-meter longitudinal range before and after the position corresponding to $T_{SV}^{lane}$. The moment $T_{SV}^{lane}$ is identified by a change in the laneID attribute of the \textit{SV}. 

\textbf{Criterion 2}: Only single discretionary lane-change events are included. Mandatory lane changes (e.g., merges and diverges) and consecutive lane changes (e.g., two or more continuous lane crossings or overtaking followed by return) are excluded.  

\textbf{Criterion 3}: To reduce confounding from nearby lane changes, only instances where no upstream vehicle changes lanes into either the target or original lanes within the event isolation window $[T_{SV}^{s}-T_{\mathrm{win}},T_{SV}^{s}+T_{\mathrm{win}}]$ are retained. In the main analysis, $T_{\mathrm{win}}=50$ s.

After applying these criteria, a total of 228 valid lane-change instances involving the \textit{SV} and its surrounding vehicles were obtained. 

\subsection{Determining the lane-change start time}\label{Section4.3} 
Estimating $T_{SV}^s$ is crucial for quantifying the lane-change impact of the \textit{SV}. Since $T_{SV}^s$ is not directly available in the raw ZTD dataset, it is estimated through the following steps: 

\vspace{0.1cm}
\textbf{Step 1: Lateral position calculation.} Data of lane-keeping vehicles in the left lane (driving lane) are first extracted, characterized by constant laneID values. The average longitude and latitude of these vehicles define the lane centerline. The lateral offset of each vehicle from this centerline is then computed as its lateral position at each time step. 

\vspace{0.1cm}
\textbf{Step 2: Identification of $T_{SV}^s$}. The start of the lane-changing process $T_{SV}^s$ is determined by analyzing the oscillation frequency of lateral position or lateral speed \citep{dong2021application,shangguan2022proactive,venthuruthiyil2022interrupted}. In this study, we adopt the method based on the lateral position. Specifically, the instant when the \textit{SV}'s lateral position ceases oscillating and begins to decrease or increase monotonically is defined as $T_{SV}^s$.

\begin{figure}[h]
  \centering
  \begin{subfigure}{0.45\linewidth}
    \centering
    \includegraphics[width=\linewidth]{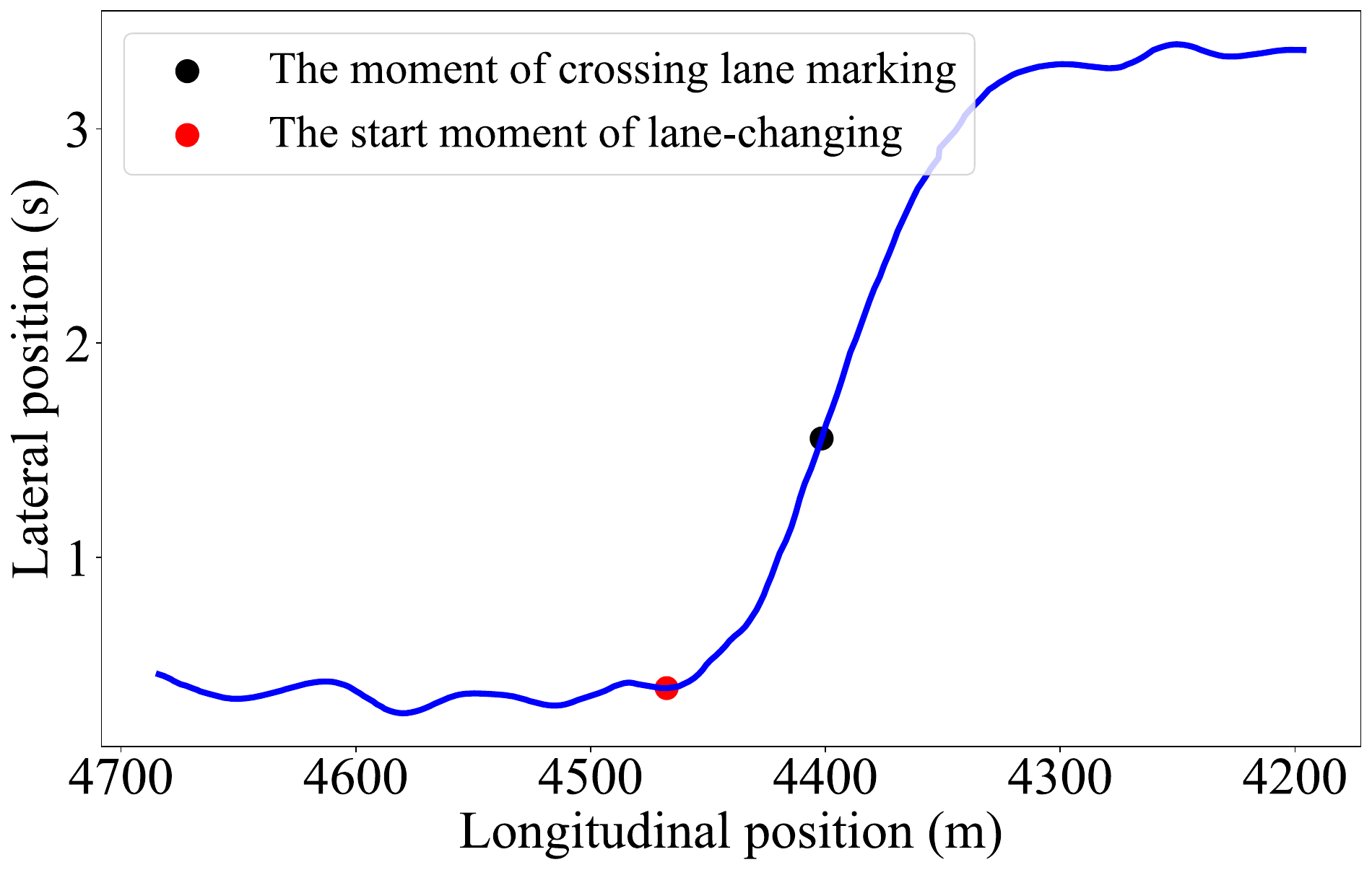}
    \caption{}
    \label{Fig5a}
  \end{subfigure}
  \hspace{0.01\linewidth} % 添加一些水平间距
  \begin{subfigure}{0.45\linewidth}
    \centering
    \includegraphics[width=\linewidth]{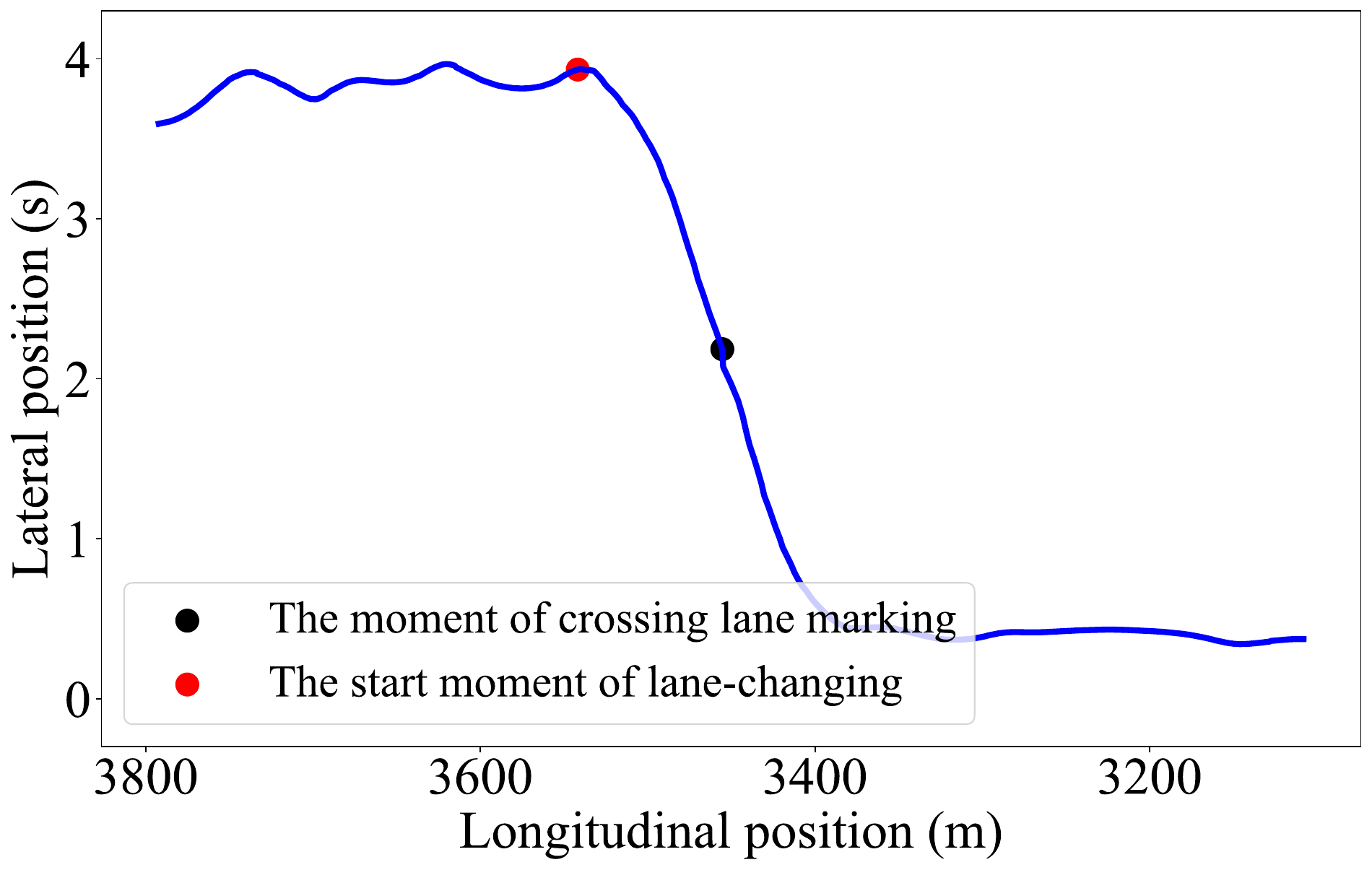}
    \caption{}
    \label{Fig5b}
  \end{subfigure}
  \caption{Examples of identifying the start time of the lane-changing process}
  \label{Fig5}
\end{figure}

\vspace{0.1cm}
\hyperref[Fig5]{Fig.~\ref{Fig5}} shows two examples of identifying the lane-change start time using this approach. 

%%%%%%%%%%%%%%%%%%%%%%%%%%%%%%%%%%%%%%%%%%%%%%%%%%%%%%%%%%%%%%%%%%%%%%%%%%%%%%%%%%%%%%%%%%%%%%%%%%%%%
\section{Case study of lane-change impact}\label{Section5}
To demonstrate the applicability of the proposed methodology, a case study was conducted using the extracted lane-change instances from \hyperref[Section4]{Section~\ref{Section4}}. The time interval $\Delta t$ was set to 0.5 seconds. A larger interval would smooth out the effects of fluctuations, while a smaller interval would overemphasize noise in TDB variations, both leading to inaccurate representations of traffic dynamics. The envelope multiplier $\eta$ in \hyperref[Eq17]{Eq.~(\ref{Eq17})} was set to 1, corresponding to a $\mu \pm \sigma$ pre-event TDB fluctuation envelope.

\subsection{Visualization of lane-change impact quantification process }\label{Section5.1} 
\begin{figure}[h]
  \centering
  \begin{subfigure}{0.47\linewidth}
    \centering
    \includegraphics[width=\linewidth]{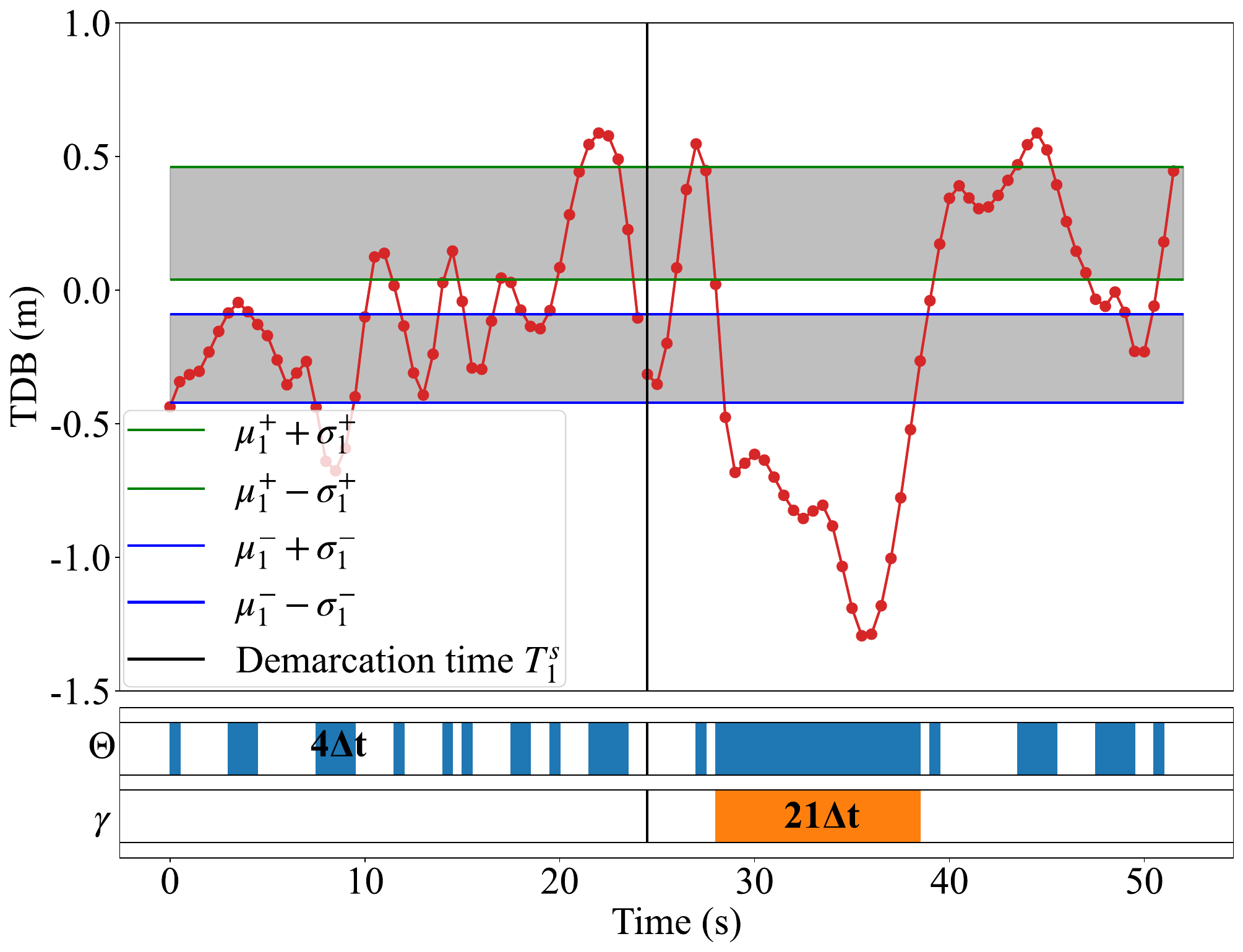}
    \caption{${TFV}_1$}
    \label{Fig6a}
  \end{subfigure}
  \hspace{0.01\linewidth} % 添加一些水平间距
  \begin{subfigure}{0.47\linewidth}
    \centering
    \includegraphics[width=\linewidth]{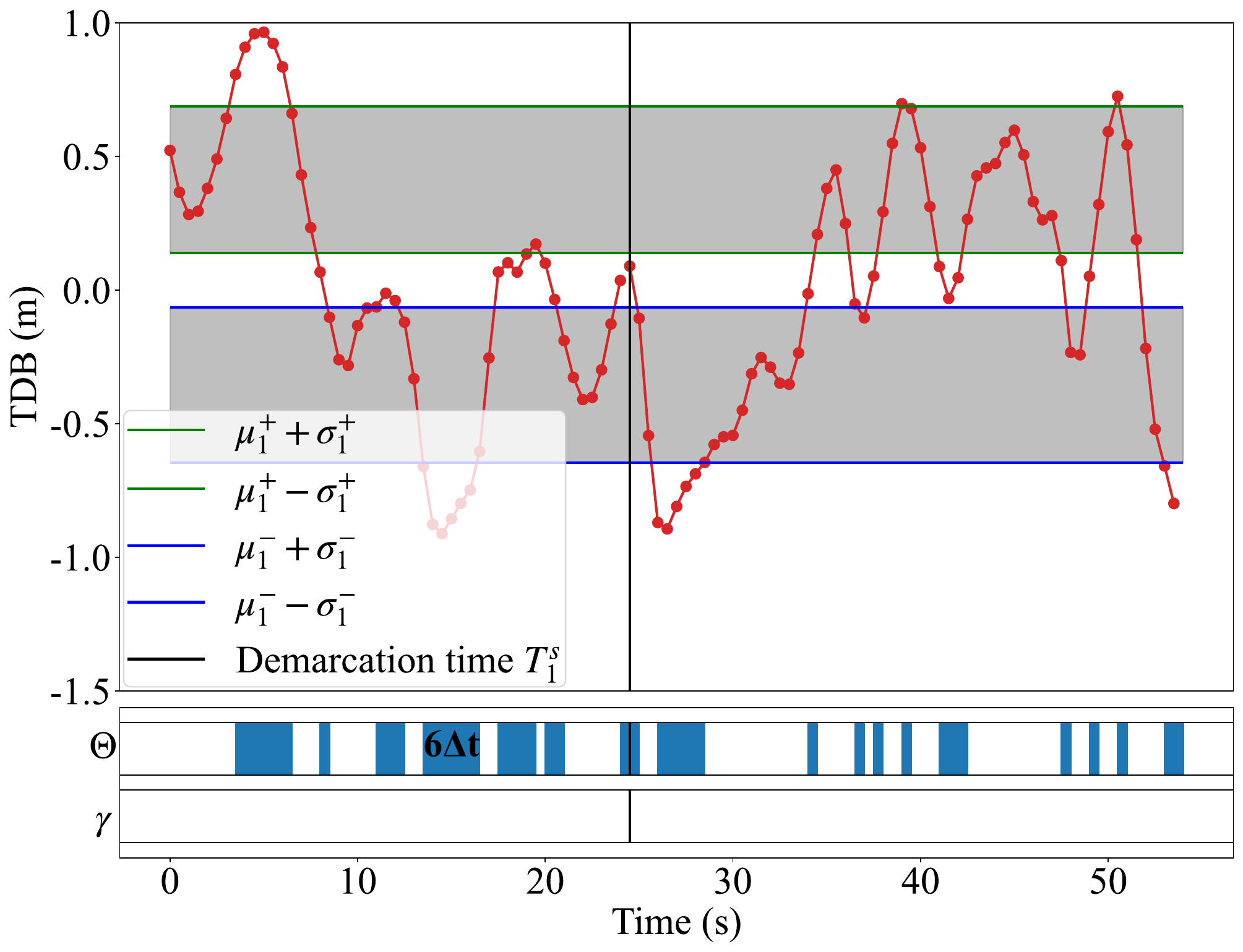}
    \caption{${FV}_1$}
    \label{Fig6b}
  \end{subfigure}
  \caption{Quantification process of lane-change impact in a lane-changing instance (vehicle 2783 in the L002\_F002\_ALL data set). White regions denote $\Theta_{\ell,i}=0$ or $\Upsilon_{\ell,i}=0$; colored regions denote $\Theta_{\ell,i}=1$ or $\Upsilon_{\ell,i}=1$}
  \label{Fig6}
\end{figure}

To effectively visualize and interpret the lane-change impact quantification, we selected one representative lane-change instance (\hyperref[Fig6]{Fig.~\ref{Fig6}}). \hyperref[Fig6a]{Fig.~\ref{Fig6a}} illustrates the process for ${TFV}_1$, and \hyperref[Fig6b]{Fig.~\ref{Fig6b}} illustrates the process for ${FV}_1$. Taking ${TFV}_1$ as an example, \hyperref[Fig6a]{Fig.~\ref{Fig6a}} shows the TDB, the potentially affected status $\Theta_{T,1}$, and estimated affected status $\Upsilon_{T,1}$ at each time interval. Colored regions in the lower part of the diagram represent $\Theta_{\ell,i}=1$ or $\Upsilon_{\ell,i}=1$, while white regions indicate $\Theta_{\ell,i}=0$ or $\Upsilon_{\ell,i}=0$. The potentially affected status ($\Theta_{T,1}$) of ${TFV}_1$ at each time interval, determined using the first judgment criterion defined in \hyperref[Eq17]{Eqs.~(\ref{Eq17})} and \hyperref[Eq18]{(\ref{Eq18})}, is shown in the second subfigure of \hyperref[Fig6a]{Fig.~\ref{Fig6a}}, where $\Omega_{T,1}^{f\ast}=4$. After the demarcation time $T_{T,1}^s$, the estimated affected status $\Upsilon_{T,1}$ is adjusted according to the second judgment criterion outlined in \hyperref[Eq21]{Eqs.~(\ref{Eq21})}-\hyperref[Eq23]{(\ref{Eq23})}, depicted in the third subfigure of \hyperref[Fig6a]{Fig.~\ref{Fig6a}}. The impact duration for ${TFV}_1$ was $21\Delta t$ ($10.5$ seconds), whereas that for ${FV}_1$ was $0$ seconds.

A comparison between the second and third subfigures in \hyperref[Fig6a]{Fig.~\ref{Fig6a}} reveals that the potentially affected status in the second subfigure fluctuates excessively within short periods (between 10 and 20 seconds), which is inconsistent with real driving behavior. In real scenarios, drivers do not respond instantaneously to lane changes; instead, their adjustments are gradual and occur over time. These frequent transitions fail to capture realistic driver adaption. By introducing the second discriminant index $\Omega_{\ell,i}^{f\ast}$ (\hyperref[Eq19]{Eq.~(\ref{Eq19})}), the procedure suppresses these short-term transitions and produces a more temporally stable estimated affected-status sequence, as shown in the third subfigure of \hyperref[Fig6a]{Fig.~\ref{Fig6a}}. 

By construction, the filter makes short status reversals less likely in the estimated sequence. Hence, the procedure provides a more temporally stable characterization of lane-change-associated disturbances while suppressing short abnormal runs under the adopted criterion.

\begin{figure}[h]
    \centering
    \begin{subfigure}{0.32\linewidth}
        \centering
        \includegraphics[width=\linewidth]{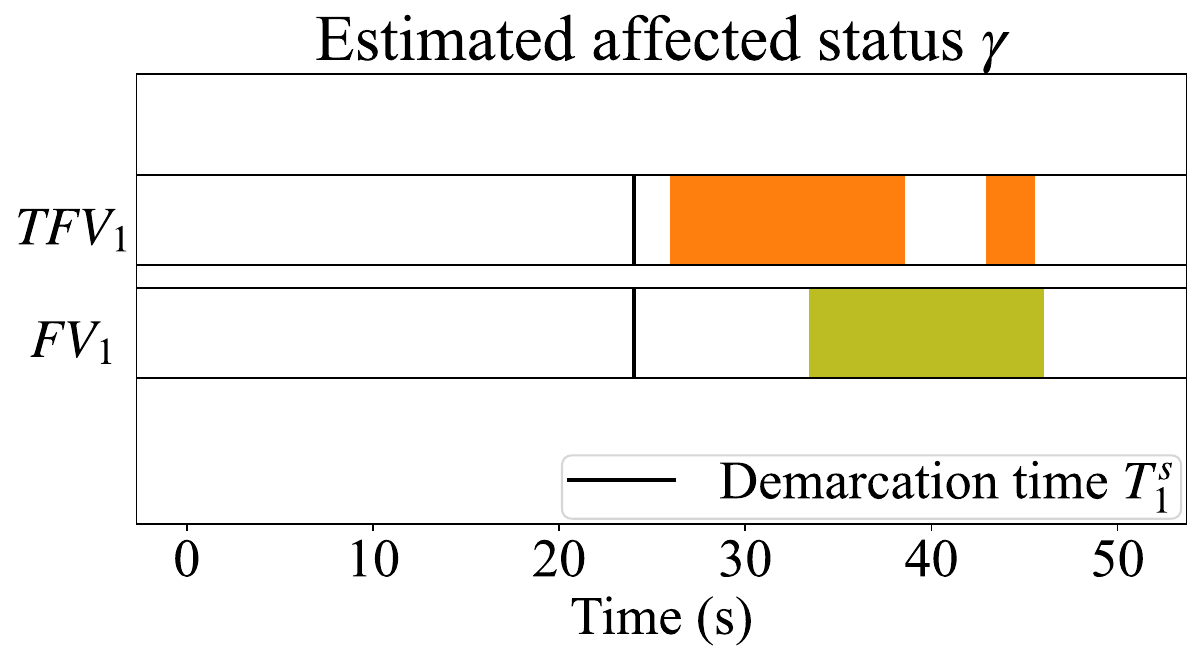}
        \caption{}
        \label{Fig7a}
    \end{subfigure}
    %\hspace{0.01\linewidth} % 添加一些水平间距
    \begin{subfigure}{0.32\linewidth}
        \centering
        \includegraphics[width=\linewidth]{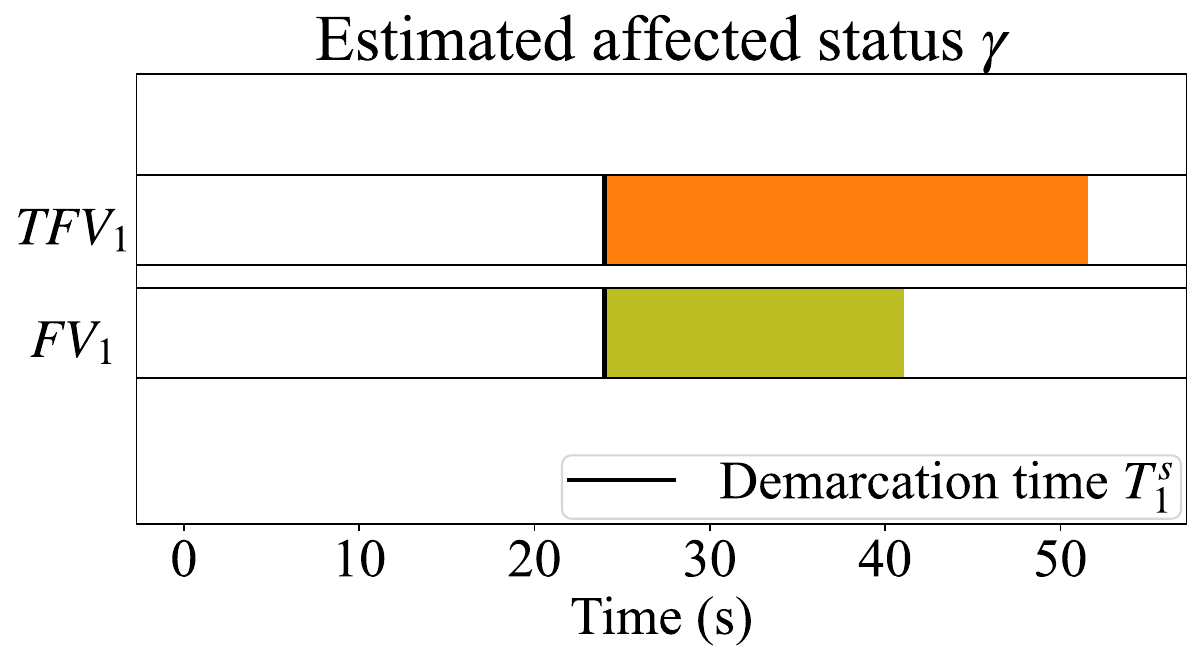}
        \caption{}
        \label{Fig7b}
    \end{subfigure}
    %\hspace{0.01\linewidth} % 添加一些水平间距
    \begin{subfigure}{0.32\linewidth}
        \centering
        \includegraphics[width=\linewidth]{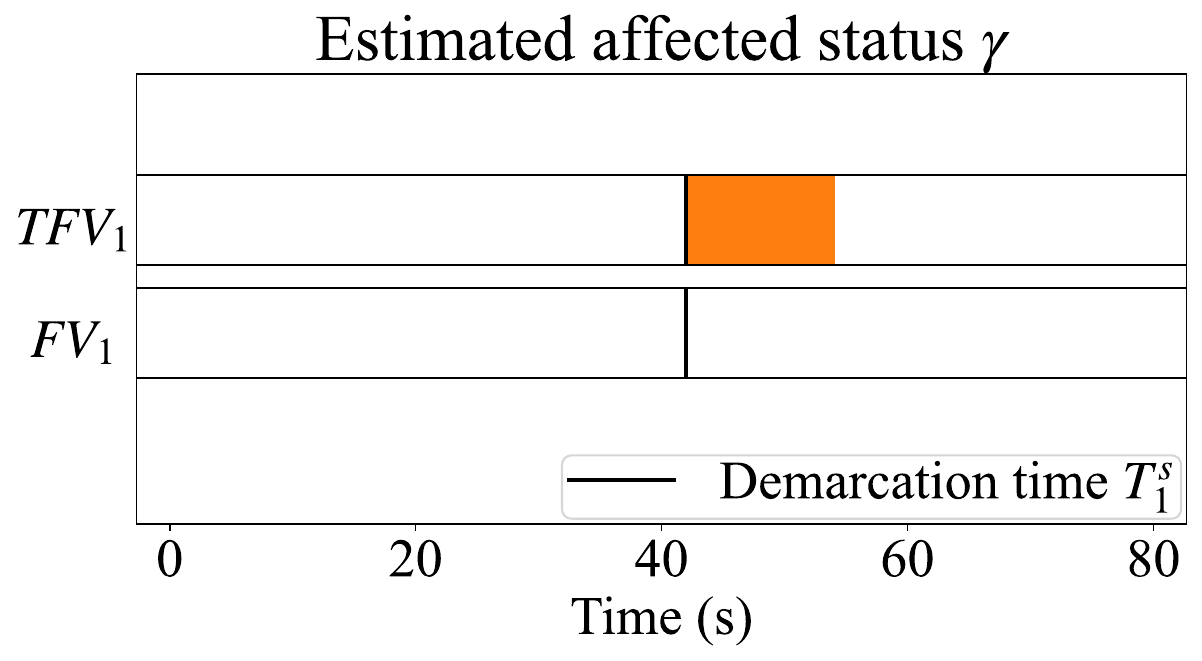}
        \caption{}
        \label{Fig7c}
    \end{subfigure}
    \caption{Estimated affected statuses of ${TFV}_1$ and ${FV}_1$ in three lane-changing instances (white: $\Upsilon_{\ell,i}=0$; colored: $\Upsilon_{\ell,i}=1$): (a) Vehicle 205 in the L001\_F002\_ALL data set; (b) Vehicle 699 in the L001\_F002\_ALL data set; and (c) Vehicle 2891 in the L002\_F002\_ALL data set}
    \label{Fig7}
\end{figure}

\hyperref[Fig7]{Fig.~\ref{Fig7}} illustrates the estimated affected statuses of ${TFV}_1$ and ${FV}_1$ in three distinct lane-changing events, highlighting the varied responses of following vehicles. The colored regions represent $\Upsilon_{\ell,i}=1$ (affected status), while white regions indicate $\Upsilon_{\ell,i}=0$ (unaffected status). When the results in \hyperref[Fig7]{Fig.~\ref{Fig7}} are analyzed alongside the third subfigures of \hyperref[Fig6a]{Figs.~\ref{Fig6a}} and \hyperref[Fig6b]{\ref{Fig6b}}, it becomes evident that the initial times when vehicles are affected differ among following vehicles, indicating heterogeneity in their responses. These initial affected times do not necessarily coincide with $T_{\ell,1}^s$ (the moment when the following vehicle first encounters the lane-specific disturbance), as shown in \hyperref[Fig7a]{Fig.~\ref{Fig7a}} and the third subfigure of \hyperref[Fig6a]{Fig.~\ref{Fig6a}}. This suggests that the following vehicle may not respond immediately to the lane-changing vehicle's maneuver. Such lag behavior aligns with the relaxation phenomenon in lane-changing processes, where drivers initially \rev{accept} shorter headways before gradually restoring normal spacings \citep{laval2008microscopic,zheng2013effects}. 

Furthermore, \hyperref[Fig7a]{Fig.~\ref{Fig7a}} reveals that ${TFV}_1$ transitions from being affected to unaffected and then back to affected. This pattern suggests that drivers may adjust their acceptable spacing through multiple steps rather than a single adjustment, highlighting the dynamic and complex nature of driver behavior during the lane-changing process. This complexity cannot be fully captured by simplistic models that assume instantaneous and single-step responses. Such discontinuous influence processes are not accounted for in existing methods \citep{zheng2013effects,he2023impact}, showcasing the superiority of our proposed approach in detecting a more nuanced impact process. In \hyperref[Fig7c]{Fig.~\ref{Fig7c}}, while ${TFV}_1$ is affected, ${FV}_1$ remains unaffected, indicating that no ${FV}_1$ response was detected under the adopted criteria.  

These observations underscore that the following vehicles exhibit varying degrees of sensitivity and responsiveness to lane-change events, influenced by factors such as driver perception, vehicle dynamics, or prevailing traffic conditions.

\subsection{Quantitative results}\label{Section5.2} 
\subsubsection{Impact on immediate followers}\label{Section5.2.1}

\begin{figure}[!htbp]
    \centering
    \begin{subfigure}{0.35\linewidth}
        \centering
        \includegraphics[width=\linewidth]{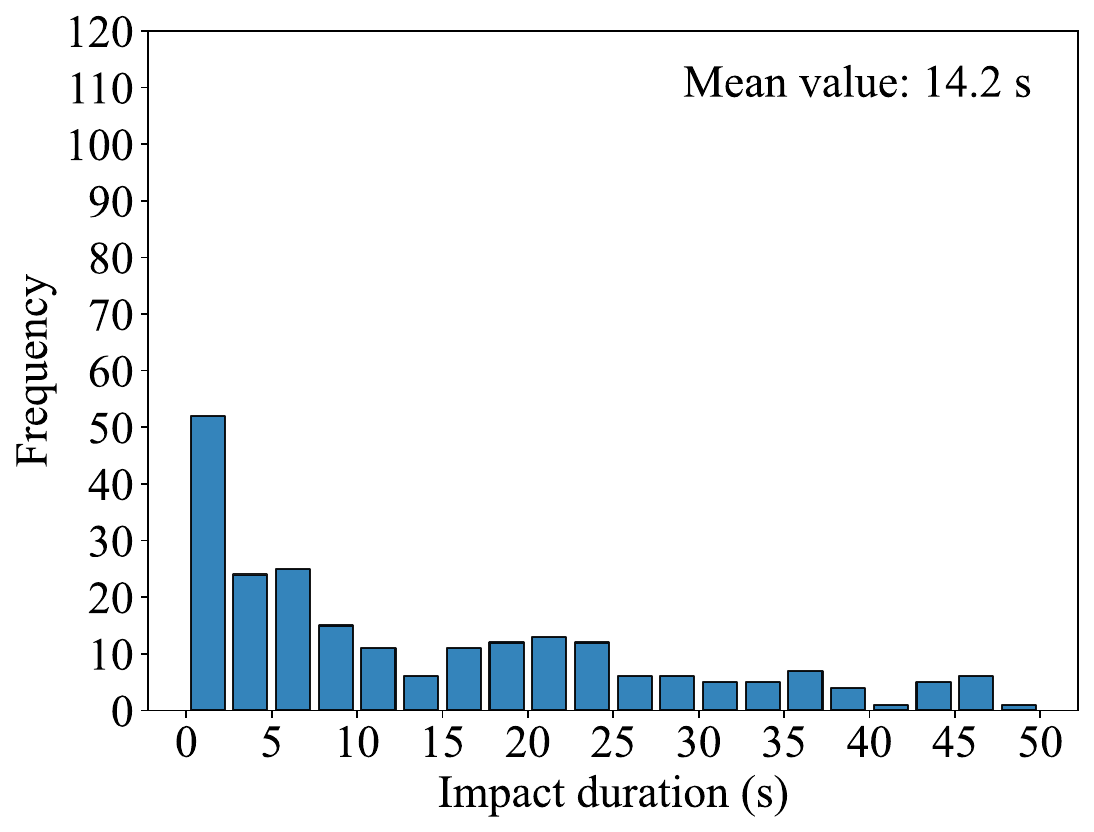}
        \caption{}
        \label{Fig8a}
    \end{subfigure}
    \hspace{0.01\linewidth} % 添加一些水平间距
    \begin{subfigure}{0.35\linewidth}
        \centering
        \includegraphics[width=\linewidth]{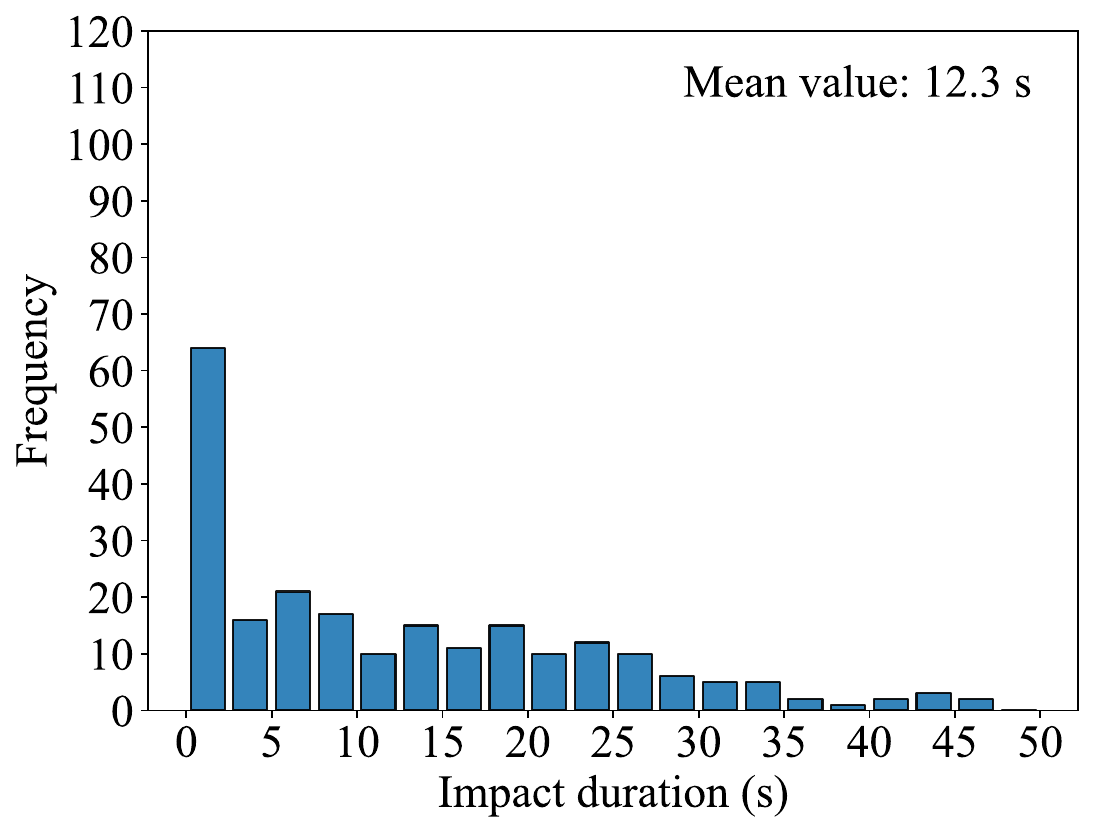}
        \caption{}
        \label{Fig8b}
    \end{subfigure}
    
    \medskip % 增加垂直间
    
    \begin{subfigure}{0.35\linewidth}
        \centering
        \includegraphics[width=\linewidth]{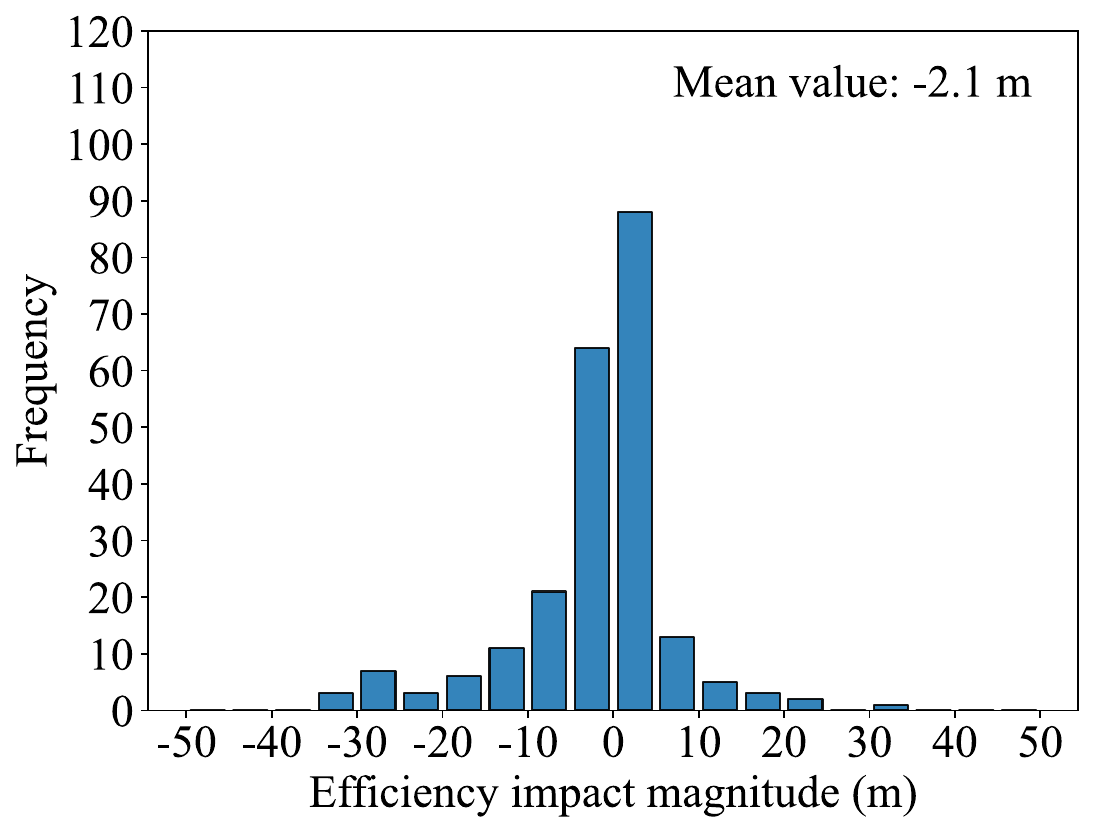}
        \caption{}
        \label{Fig8c}
    \end{subfigure}
    \hspace{0.01\linewidth} % 添加一些水平间距
    \begin{subfigure}{0.35\linewidth}
        \centering
        \includegraphics[width=\linewidth]{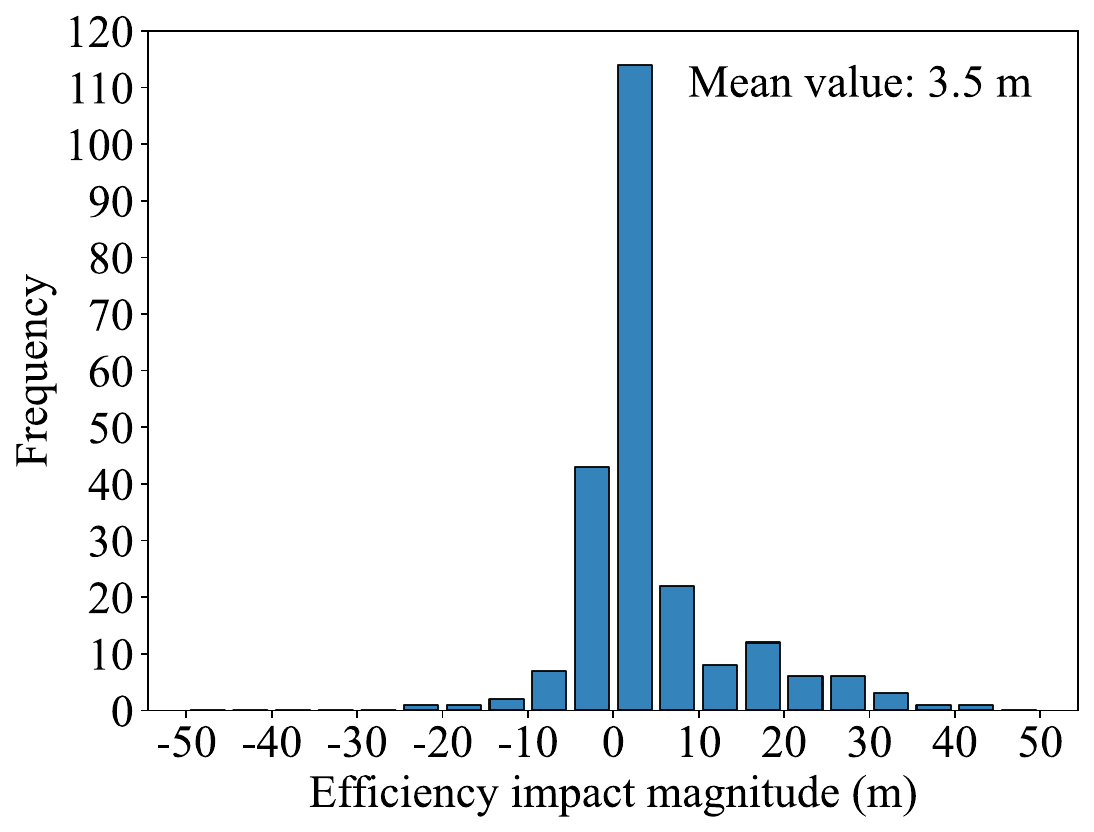}
        \caption{}
        \label{Fig8d}
    \end{subfigure}

    \medskip % 增加垂直间

    \begin{subfigure}{0.35\linewidth}
        \centering
        \includegraphics[width=\linewidth]{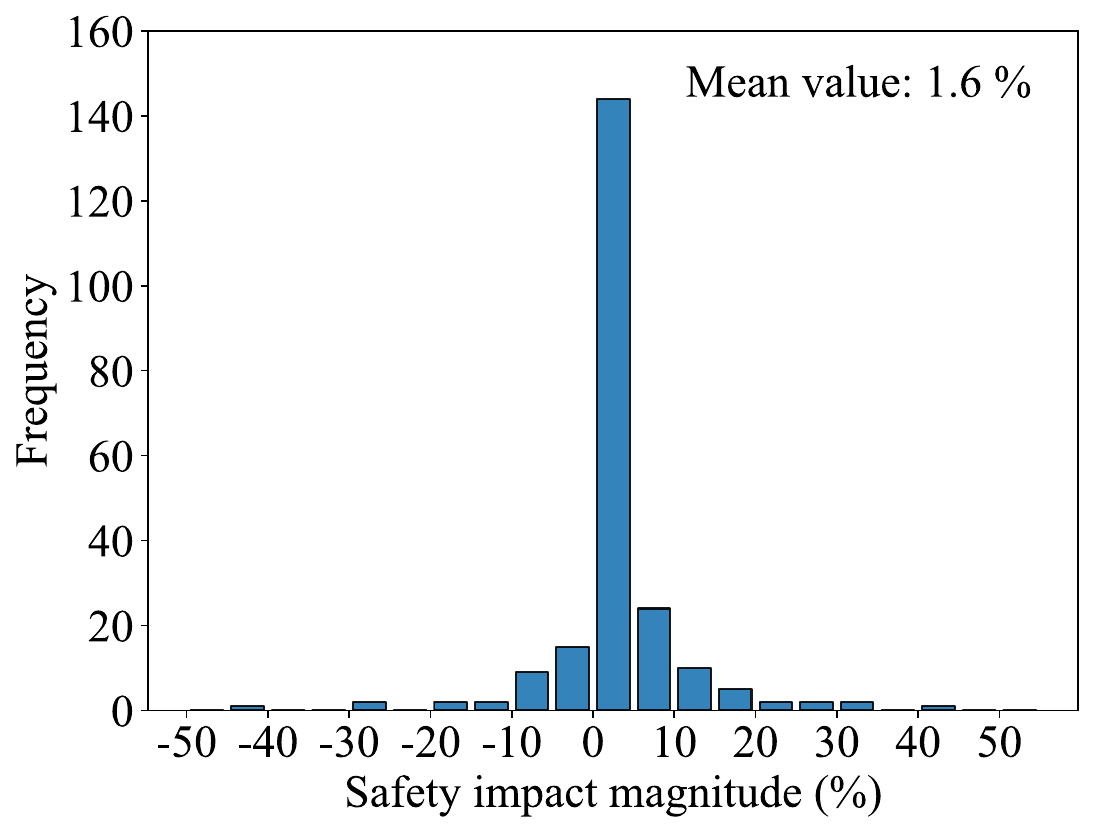}
        \caption{}
        \label{Fig8e}
    \end{subfigure}
    \hspace{0.01\linewidth} % 添加一些水平间距
    \begin{subfigure}{0.35\linewidth}
        \centering
        \includegraphics[width=\linewidth]{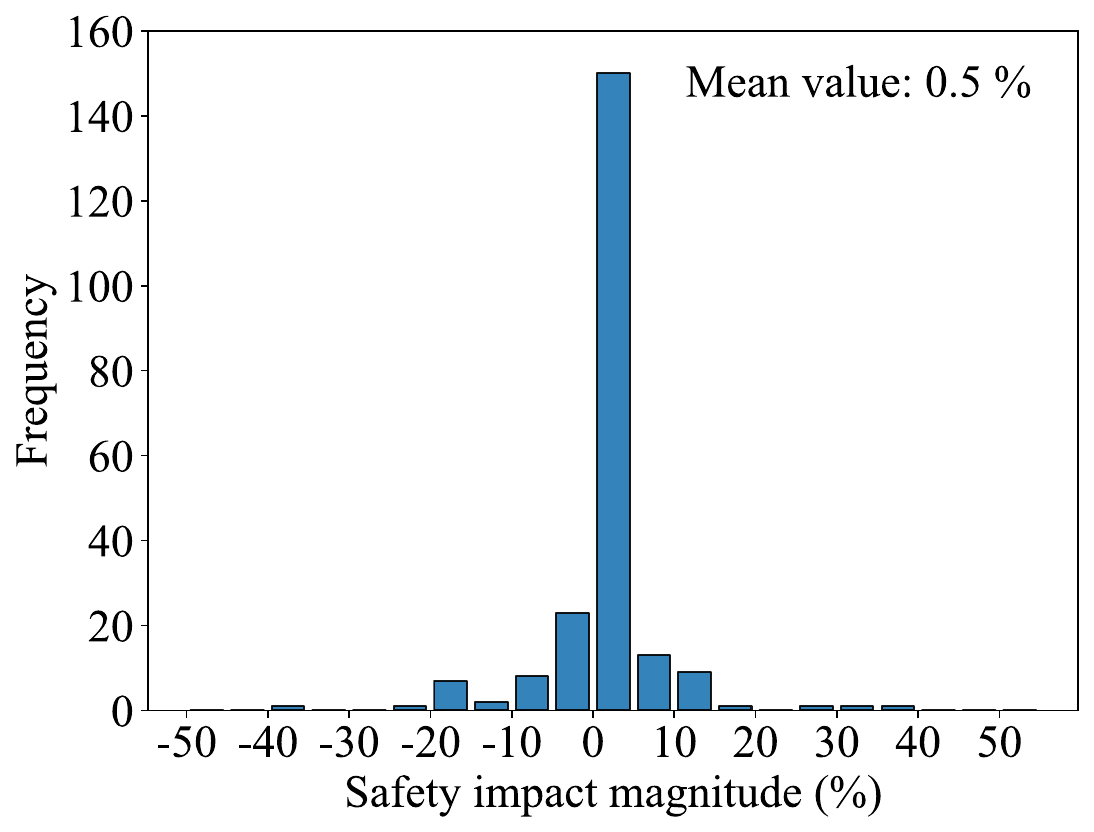}
        \caption{}
        \label{Fig8f}
    \end{subfigure}    
    \caption{Distribution of lane-change impact for the first following vehicles: (a) Impact durations for ${TFV}_1$; (b) Impact durations for ${FV}_1$; (c) EIM for ${TFV}_1$ and (d) EIM for ${FV}_1$; (e) SIM for ${TFV}_1$; (f) SIM for ${FV}_1$ }
    \label{Fig8}
\end{figure}

\hyperref[Fig8]{Fig.~\ref{Fig8}} illustrates the distribution of lane-change impact duration, efficiency impact magnitude (EIM) and safety impact magnitude (SIM) for the first following vehicles in both the target and original lanes. The average impact duration ($T_{\ell,i}^A$ as defined in \hyperref[Eq24]{Eq.~(\ref{Eq24})}), EIM ($w_{\ell,i}^A$ as defined in \hyperref[Eq35]{Eq.~(\ref{Eq35})}) and SIM ($r_{\ell,i}^A$ as defined in \hyperref[Eq40]{Eq.~(\ref{Eq40})}) for ${TFV}_1$ were approximately 14.2 seconds, -2.1 meters and 1.6\%, respectively. In contrast, for ${FV}_1$, these values were approximately 12.3 seconds, 3.5 meters, and 0.5\%, respectively.

A comparison of \hyperref[Fig8]{Figs.~\ref{Fig8c}} and \hyperref[Fig8d]{\ref{Fig8d}} shows that lane changes tend to reduce traffic efficiency for the immediate follower in the target lane but improve efficiency for the immediate follower in the original lane. This difference is consistent with the two lane-specific mechanisms: \textit{SV} insertion compresses the target-lane gap and often causes ${TFV}_1$ to travel a shorter distance relative to \textit{TLV}, whereas \textit{SV} departure releases space in the original lane and may allow ${FV}_1$ to accelerate or close the gap toward \textit{LV}. Similarly, \hyperref[Fig8]{Figs.~\ref{Fig8e}} and \hyperref[Fig8f]{\ref{Fig8f}} show that lane changes increase safety risks for immediate followers in both lanes, with a more pronounced rise in the target lane (1.6\% rise in unsafe proportion) than in the original lane (0.5\% rise). The original-lane risk increase is not contradictory to its positive EIM: acceleration and gap-closing behavior can increase relative speed to \textit{LV} and reduce TTC-based safety margins even when the travel-distance deviation indicates an efficiency gain. Overall, these findings reveal that the negative impacts of lane changes on both efficiency and safety are more pronounced for the immediate follower in the target lane than for the immediate follower in the original lane.

Importantly, this lane-specific asymmetry in efficiency and safety effects has not been captured by previous methods such as \citet{he2023impact}, which quantify only the affected duration or the number of affected vehicles without distinguishing positive and negative impacts. By explicitly incorporating spatiotemporal dimensions and differentiating between beneficial and adverse effects, our approach provides a richer and more nuanced characterization of lane-change impacts.

\subsubsection{Overall impact on upstream traffic}\label{Section5.2.2}

\begin{figure}[!ht]
    \centering
    \begin{subfigure}{0.35\linewidth}
        \centering
        \includegraphics[width=\linewidth]{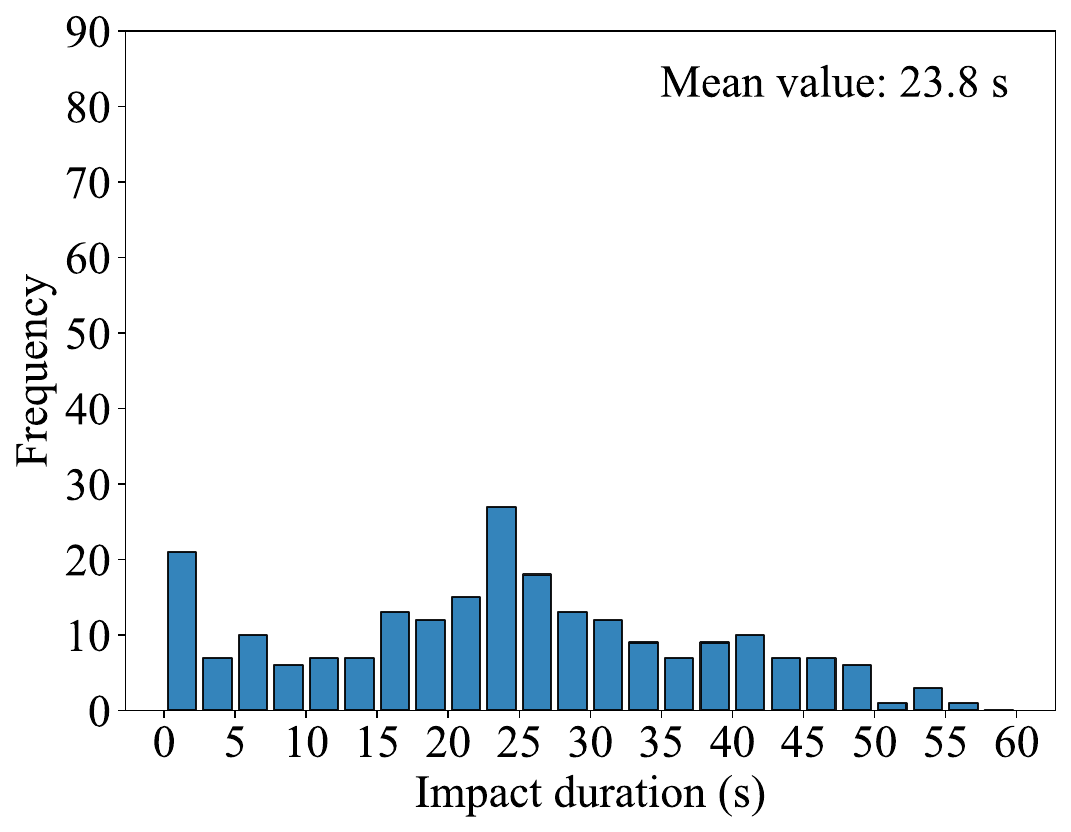}
        \caption{}
        \label{Fig9a}
    \end{subfigure}
    \hspace{0.01\linewidth} % 添加一些水平间距
    \begin{subfigure}{0.35\linewidth}
        \centering
        \includegraphics[width=\linewidth]{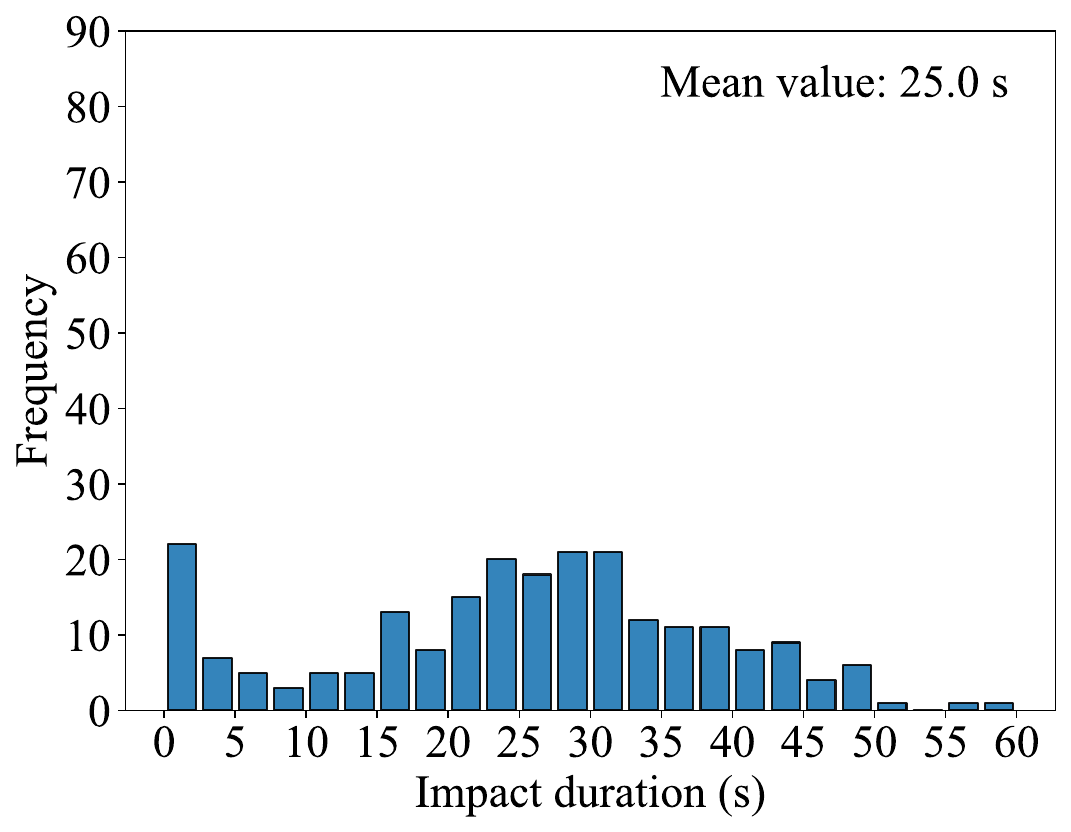}
        \caption{}
        \label{Fig9b}
    \end{subfigure}

    \medskip % 增加垂直间
    
    \begin{subfigure}{0.35\linewidth}
        \centering
        \includegraphics[width=\linewidth]{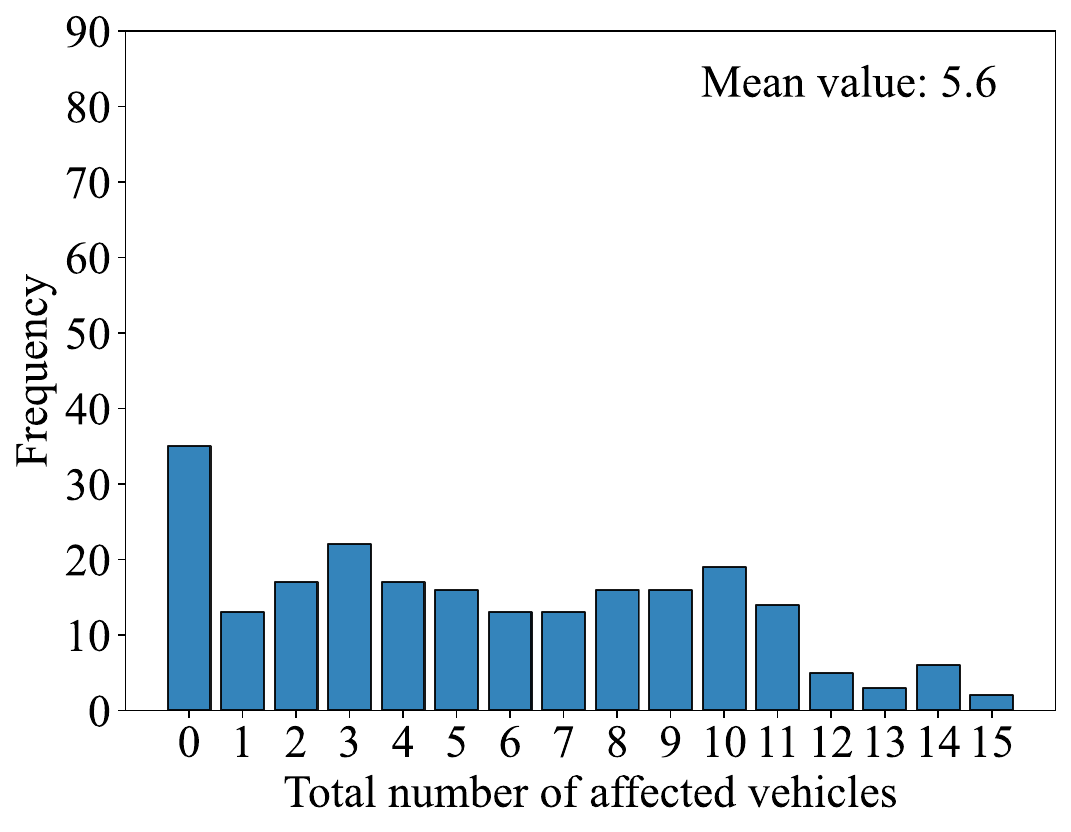}
        \caption{}
        \label{Fig9c}
    \end{subfigure}
    \hspace{0.01\linewidth} % 添加一些水平间距
    \begin{subfigure}{0.35\linewidth}
        \centering
        \includegraphics[width=\linewidth]{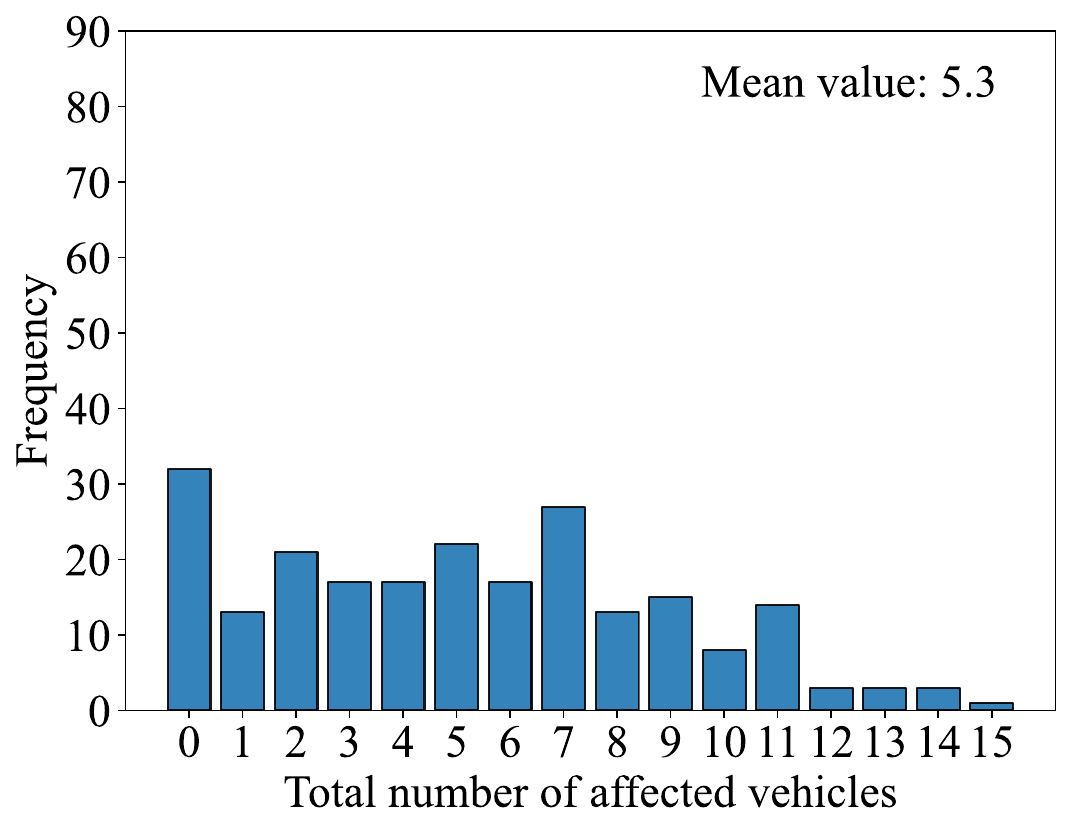}
        \caption{}
        \label{Fig9d}
    \end{subfigure}
    
    \medskip % 增加垂直间
    
    \begin{subfigure}{0.35\linewidth}
        \centering
        \includegraphics[width=\linewidth]{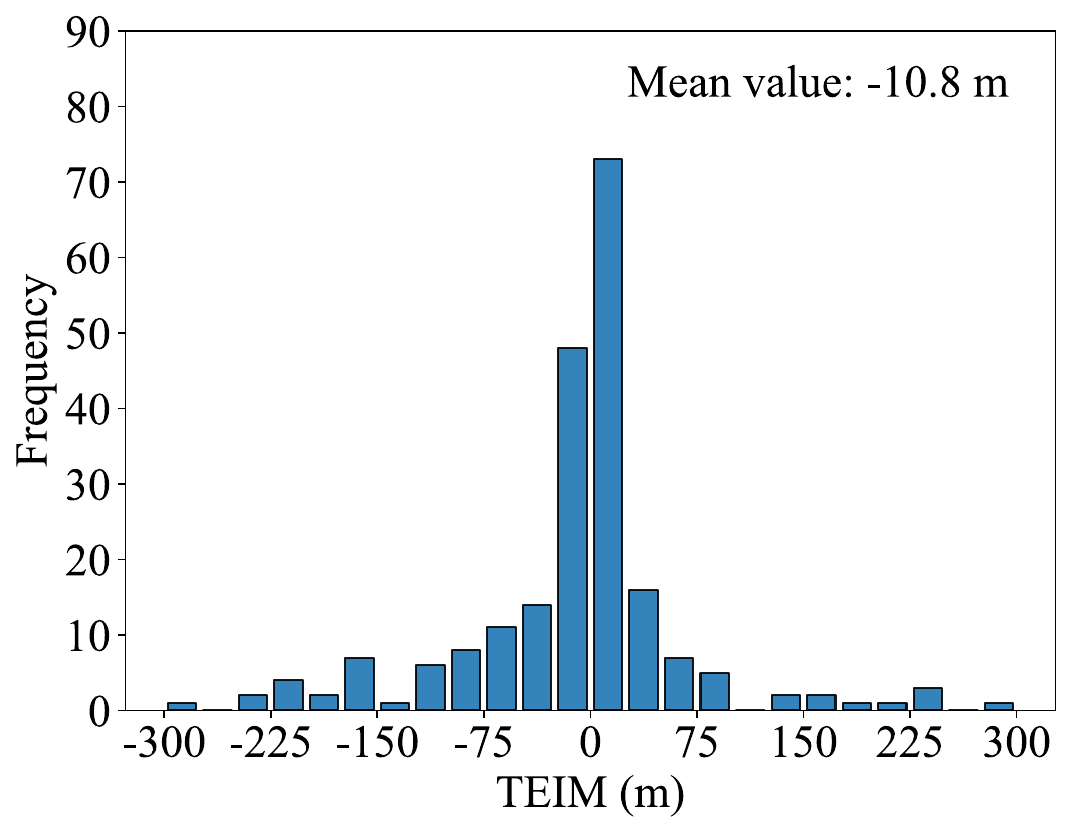}
        \caption{}
        \label{Fig9e}
    \end{subfigure}
    \hspace{0.01\linewidth} % 添加一些水平间距
    \begin{subfigure}{0.35\linewidth}
        \centering
        \includegraphics[width=\linewidth]{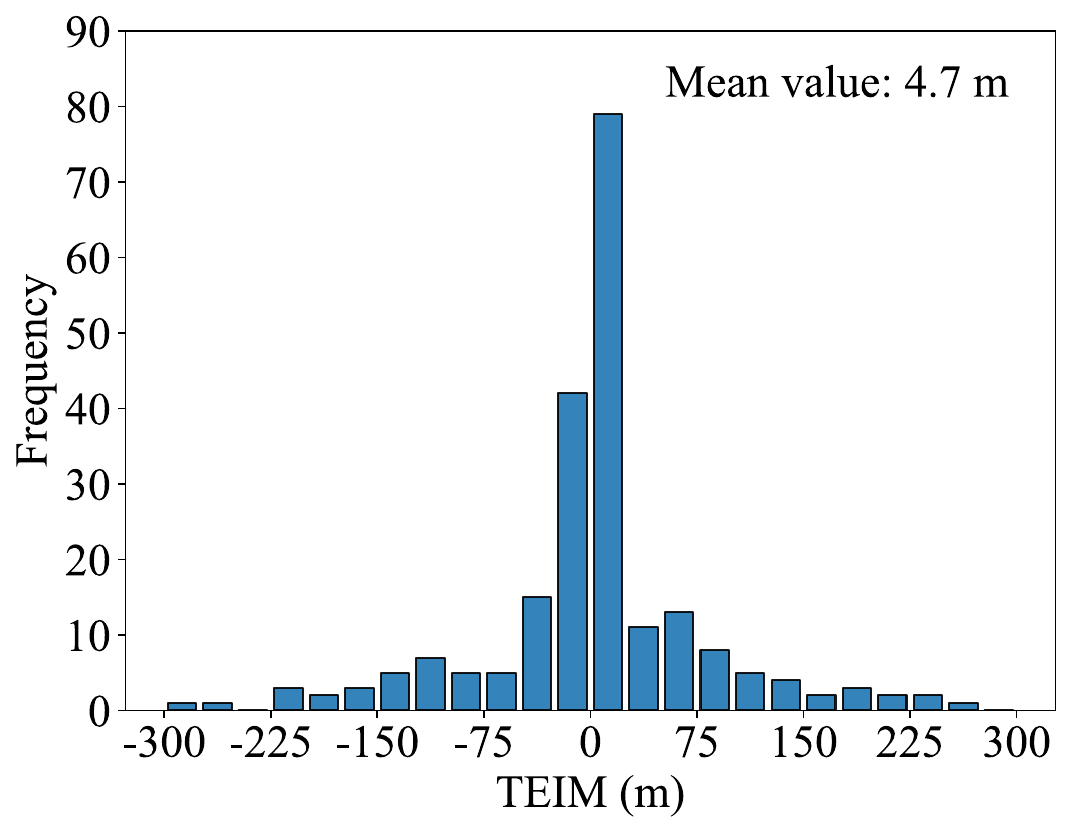}
        \caption{}
        \label{Fig9f}
    \end{subfigure}

    \medskip % 增加垂直间

    \begin{subfigure}{0.35\linewidth}
        \centering
        \includegraphics[width=\linewidth]{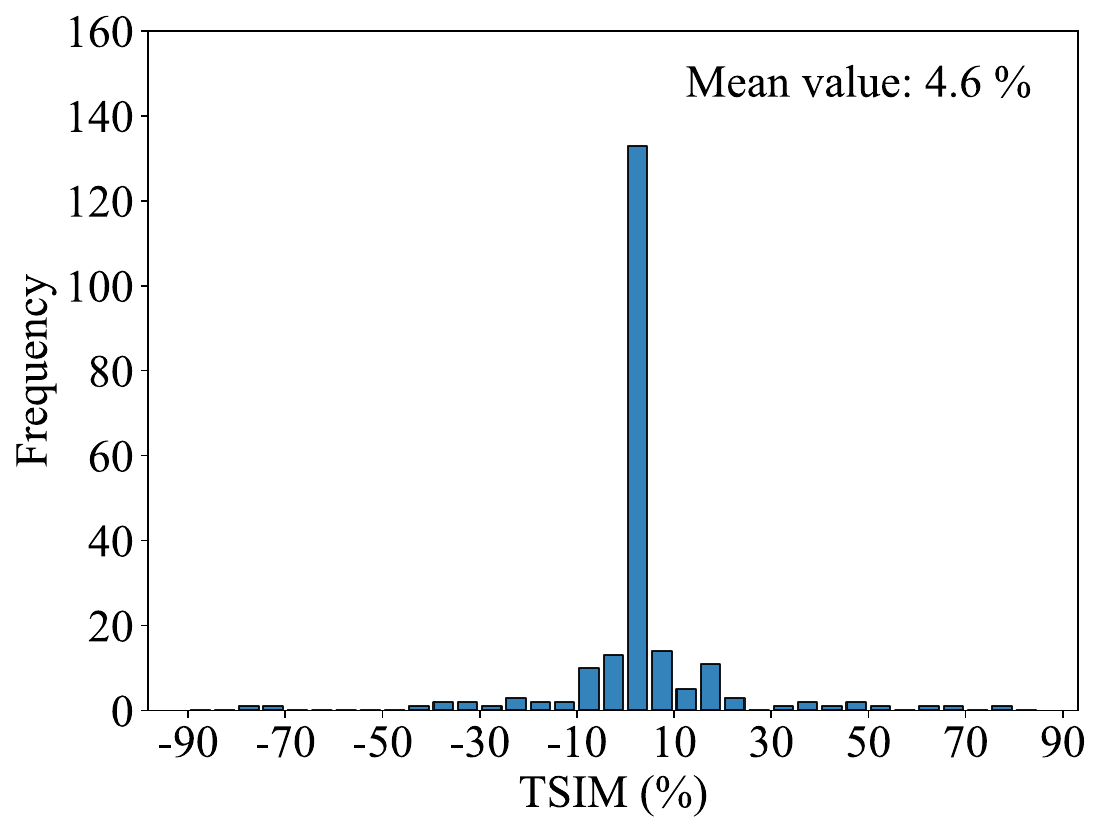}
        \caption{}
        \label{Fig9g}
    \end{subfigure}
    \hspace{0.01\linewidth} % 添加一些水平间距
    \begin{subfigure}{0.35\linewidth}
        \centering
        \includegraphics[width=\linewidth]{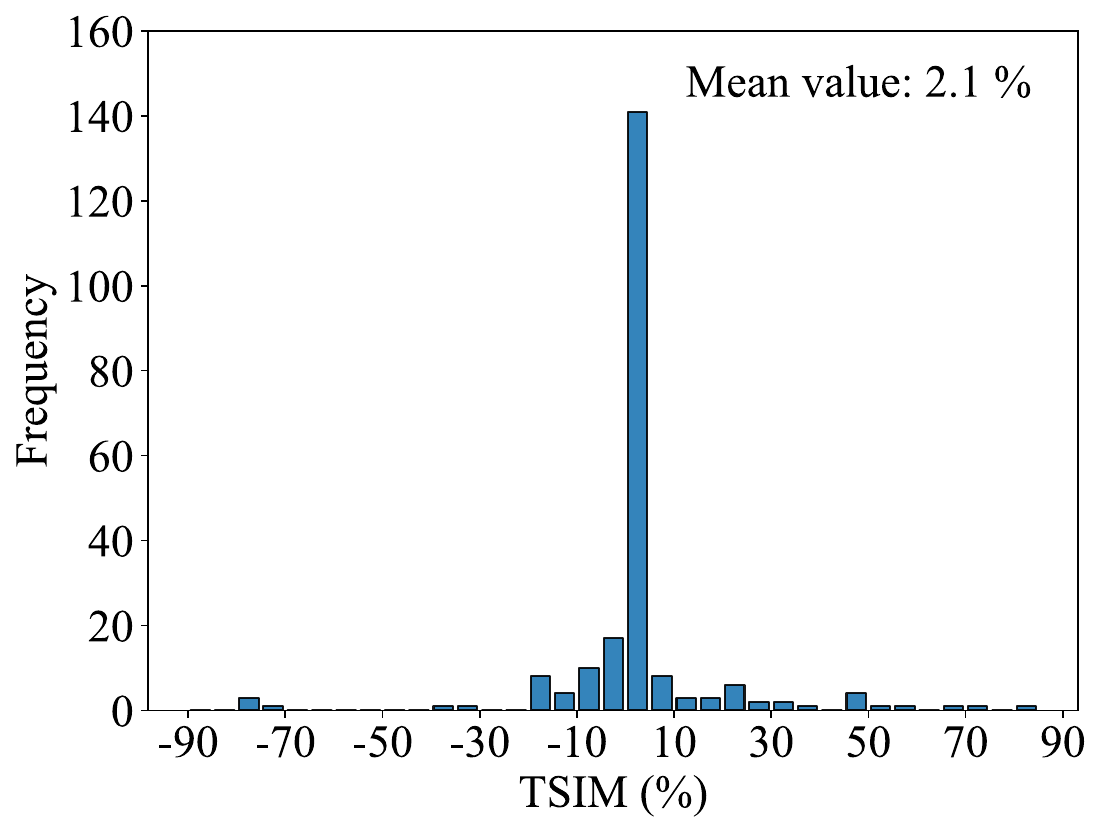}
        \caption{}
        \label{Fig9h}
    \end{subfigure}
    
    \caption{Distributions of lane-change impact duration, total number of affected vehicles, TEIM, and TSIM across lanes: (a)(c)(e)(g) Target lane; (b)(d)(f)(h) Original lane}
    \label{Fig9}
\end{figure}

We further analyze the average overall impact duration ($T_\ell^A$, \hyperref[Eq32]{Eq.~(\ref{Eq32})}), the total number of affected following vehicles ($N_\ell^A$, \hyperref[Eq27]{Eq.~(\ref{Eq27})}), the total efficiency impact magnitude (TEIM) ($W_{\ell,effi}^A$, \hyperref[Eq36]{Eq.~(\ref{Eq36})}), and the total safety impact magnitude (TSIM) ($R_{\ell,unsafe}^A$, \hyperref[Eq41]{Eq.~(\ref{Eq41})}) identified for a single lane-change event, as illustrated in \hyperref[Fig9]{Fig.~\ref{Fig9}}. The results show that the average impact durations are 23.8 seconds in the target lane and 25 seconds in the original lane. The average numbers of affected following vehicles are 5.6 in the target lane and 5.3 in the original lane. In terms of efficiency, the TEIM is -10.8 meters in the target lane and +4.7 meters in the original lane. For safety, the TSIM is 4.6\% in the target lane and 2.1\% in the original lane. 

These findings reveal several key insights: 1) Comparable durations between the target and original lanes indicate that the framework identifies a broad and sustained spatiotemporal disturbance extending across both lanes; 2) The negative TEIM in the target lane indicates efficiency loss (increased travel distance delay) for upstream vehicles and is consistent with gap compression and deceleration after \textit{SV} insertion, whereas the positive TEIM value in the original lane implies a slight efficiency gain (reduced travel distance delay) for upstream vehicles left behind after the \textit{SV}'s departure; 3) The consistently positive TSIM values in both lanes demonstrate that lane changes increase safety risks for upstream vehicles in both the target and original lanes regardless of efficiency effects, with a more pronounced effect in the target lane. In the original lane, this can occur because the released space encourages gap-closing behavior: the following vehicle may travel farther relative to \textit{LV}, but a higher closing speed or smaller remaining spacing can still increase the proportion of time with low TTC.

\begin{figure}[h]
  \centering
  \begin{subfigure}{0.35\linewidth}
        \centering
        \includegraphics[width=\linewidth]{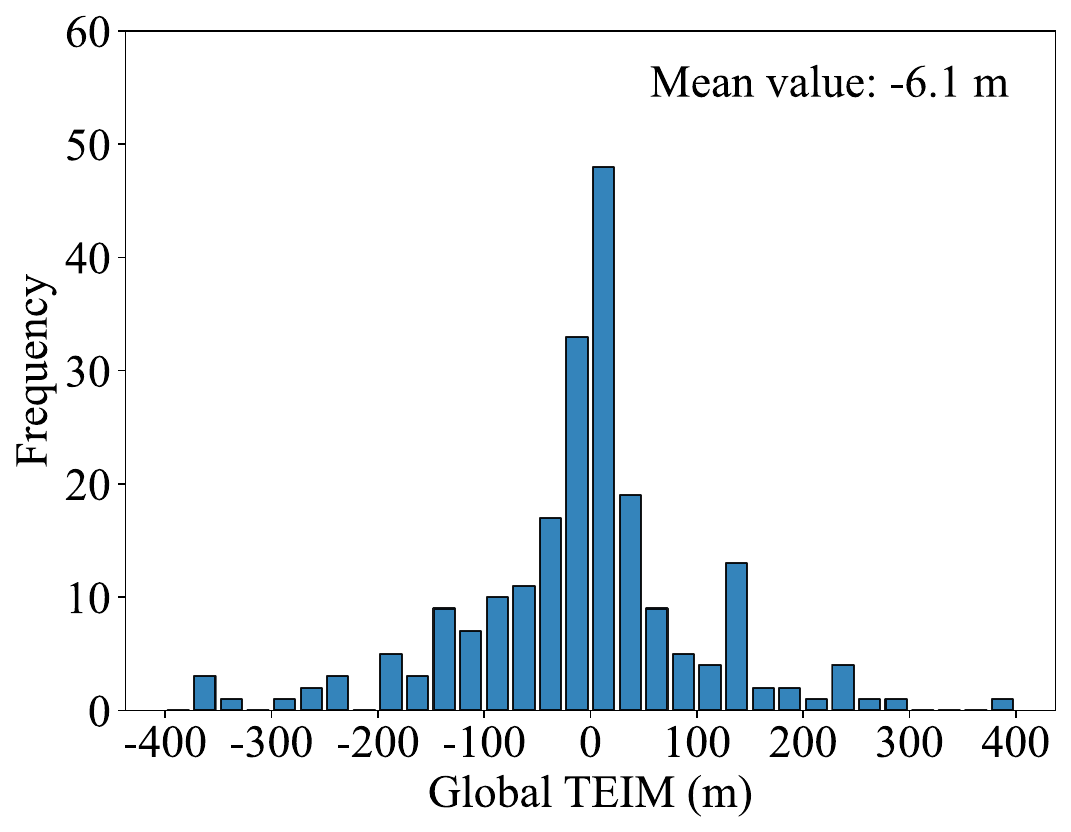}
        \caption{}
        \label{Fig10a}
    \end{subfigure}
    \hspace{0.01\linewidth} % 添加一些水平间距
    \begin{subfigure}{0.35\linewidth}
        \centering
        \includegraphics[width=\linewidth]{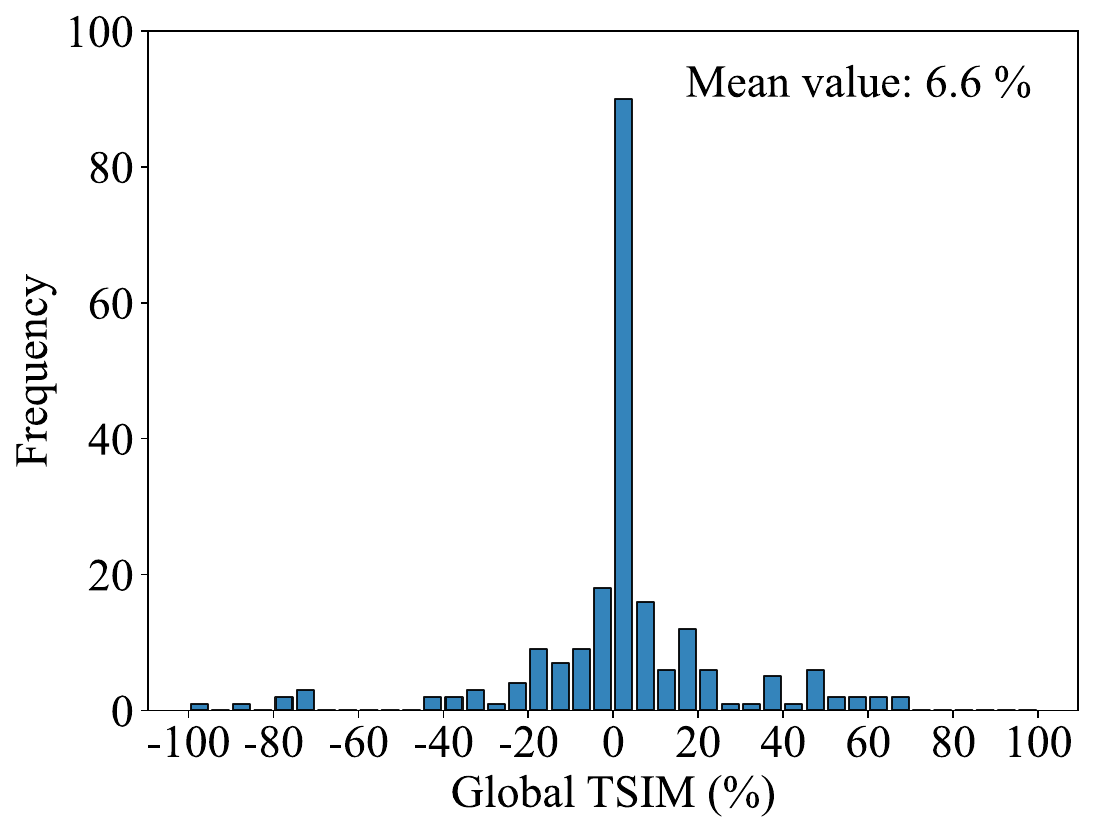}
        \caption{}
        \label{Fig10b}
    \end{subfigure}
    
  \caption{Global total impact magnitude for both the target and original lanes: (a) Global TEIM and (b) Global TSIM}
  \label{Fig10}
\end{figure}

To assess the joint impact of a single lane change on both the target and original lanes, we compute the global TEIM and global TSIM, defined as the sum of the impacts from the target and original lanes (\hyperref[Fig10]{Fig.~\ref{Fig10}}). On average, the global TEIM is -6.1 meters and the global TSIM is 6.6\%, indicating overall negative effects of lane changes on both traffic efficiency and safety. Nevertheless, the tails of both distributions contain positive TEIM or negative TSIM values, suggesting that under certain traffic conditions or driver behaviors, lane changes may benefit overall traffic flow, enhancing traffic efficiency or mitigating risk rather than deteriorating it. 

\rev{The event-level frequency distributions and mean values of the impact metrics are presented in Figs.~\ref{Fig9} and~\ref{Fig10}. These mean values summarize the average direction and magnitude of the identified lane-change impacts. Table~\ref{tableFig9Fig10Inference} complements these figures by reporting the mean, bootstrap 95\% confidence interval (CI) for the mean, and median with interquartile range (IQR) for each metric. The bootstrap CI was obtained from 10,000 event-level resamples with replacement, using the 2.5th and 97.5th percentiles of the bootstrap mean distribution as the lower and upper confidence limits. The median and IQR characterize the central tendency and variability across lane-change events. For the signed event-level aggregate magnitude metrics, namely TEIM and TSIM, a two-sided Wilcoxon signed-rank test was applied to the nonzero event-level values under the null hypothesis that their distribution is symmetric about zero. Zero-valued observations were retained in the mean, bootstrap CI, median, and IQR calculations but were excluded from the Wilcoxon calculation because they receive no positive or negative rank. Thus, the bootstrap interval estimates uncertainty in the mean across all events, whereas the Wilcoxon test assesses symmetry about zero using the signed ranks of the nonzero event-level values; consequently, the two statistics need not yield identical conclusions. Wilcoxon tests were not applied to impact duration or affected-vehicle count because these are nonnegative measures.}

\begin{table}[H]
  \centering
  \captionsetup{font=small}
  \caption{\rev{Summary statistics for lane-change impact metrics}}
  \label{tableFig9Fig10Inference}
  {
  \small
  \setlength{\tabcolsep}{3pt}
  \begin{tabular*}{\linewidth}{@{\extracolsep{\fill}}lllll@{}}
    \toprule[1.0pt]
    \textbf{Metric} & \textbf{Mean} & \textbf{95\% CI} & \rev{\textbf{Median [Q1, Q3]}} & \textbf{Wilcoxon $p$-value} \\
    \midrule[0.5pt]
    Target-lane impact duration (s) & 23.8 & [21.9, 25.6] & \rev{23.6 [14.5, 33.1]} & \rev{--} \\
    Original-lane impact duration (s) & 25.0 & [23.3, 26.7] & \rev{26.4 [16.5, 33.5]} & \rev{--} \\
    Target-lane affected vehicles & 5.6 & [5.1, 6.1] & \rev{5.0 [2.0, 9.0]} & \rev{--} \\
    Original-lane affected vehicles & 5.3 & [4.8, 5.8] & \rev{5.0 [2.0, 8.0]} & \rev{--} \\
    Target-lane TEIM (m) & $-10.8$ & [$-21.4$, $-0.2$] & \rev{0.0 [$-27.1$, 12.9]} & 0.086 \\
    Original-lane TEIM (m) & 4.7 & [$-6.4$, 15.9] & \rev{0.0 [$-12.4$, 21.2]} & 0.237 \\
    Target-lane TSIM (\%) & 4.6 & [1.1, 8.3] & \rev{0.0 [0.0, 0.7]} & 0.014 \\
    Original-lane TSIM (\%) & 2.1 & [$-0.6$, 4.9] & \rev{0.0 [0.0, 0.0]} & 0.263 \\
    Global TEIM (m) & $-6.1$ & [$-20.9$, 9.8] & \rev{0.7 [$-50.8$, 35.5]} & 0.618 \\
    Global TSIM (\%) & 6.6 & [1.7, 11.9] & \rev{0.0 [$-2.4$, 8.7]} & 0.025 \\
    \bottomrule[1.0pt]
  \end{tabular*}
  \par\vspace{2pt}
  \begin{minipage}{\linewidth}
  \footnotesize\rev{\textit{Note:} ``--'' indicates that the Wilcoxon test was not applied to the corresponding metric.}
  \end{minipage}
  }
\end{table}
\vspace{-1.0\baselineskip}

\rev{The statistics in Table~\ref{tableFig9Fig10Inference} describe both uncertainty in the estimated means and heterogeneity across lane-change events. For impact duration and affected-vehicle count, the means and medians are broadly similar in both lanes, while the bootstrap confidence intervals around the means are comparatively narrow. Their interquartile ranges nevertheless indicate appreciable event-to-event variation in the identified impact ranges. In contrast, for TEIM and TSIM, the means differ more noticeably from the medians, and the bootstrap confidence intervals are generally wider relative to the magnitudes of the means. Several confidence intervals for these signed aggregate magnitude metrics include zero, indicating greater uncertainty in the direction and magnitude of their estimated mean values.} These patterns are consistent with the construction of TEIM and TSIM, which accumulate signed efficiency- and safety-related changes over affected vehicles and time intervals and are therefore sensitive to event-level heterogeneity, cancellation between positive and negative values, and local driving responses. This is especially evident for global TEIM, where target-lane efficiency loss can be offset by original-lane efficiency gain.

\rev{A lane-specific pattern is also evident across both signed aggregate magnitude metrics. For TEIM, both lanes have a median of 0.0 m, but the original-lane interquartile range of [$-12.4$, 21.2] m is shifted upward relative to the target-lane range of [$-27.1$, 12.9] m. Meanwhile, the bootstrap confidence interval for mean target-lane TEIM [$-21.4$, $-0.2$] m is entirely negative, whereas that for original-lane TEIM [$-6.4$, 15.9] m includes zero. Taken together, these results show a clearer average efficiency-loss direction in the target lane. The original-lane distribution is shifted toward more positive, efficiency-improving values, but the direction of its estimated mean remains less certain. For TSIM, the target-lane median is 0.0\%, with an interquartile range of [0.0, 0.7]\%, and its bootstrap confidence interval [1.1, 8.3]\% excludes zero. By contrast, the original-lane median is 0.0\%, with an interquartile range of [0.0, 0.0]\%, and its confidence interval [$-0.6$, 4.9]\% includes zero. Together, these results indicate a clearer direction toward increased TTC-based unsafe exposure in the target lane, whereas the corresponding original-lane direction is weaker and less consistent across events.} This lane-specific contrast is consistent with target-lane following vehicles being more passively exposed to gap compression after \textit{SV} insertion. Original-lane following vehicles have greater behavioral discretion after \textit{SV} departure; some may accelerate or maintain speed, leading to different gap-closing patterns and TTC-based unsafe exposure across events. Such heterogeneous responses contribute to greater uncertainty in original-lane TEIM and TSIM. \rev{The Wilcoxon results provide supplementary information about the nonzero event-level values. The test was significant for target-lane TSIM ($p=0.014$) but not for original-lane TSIM ($p=0.263$). Together with the positive target-lane TSIM mean and its bootstrap confidence interval above zero, this contrast is consistent with a stronger directional tendency toward increased TTC-based unsafe exposure in the target lane.}

From an operational perspective, the reported TEIM and TSIM values should be interpreted as aggregate event-level magnitude indicators rather than direct per-vehicle operational outcomes. For example, a target-lane TEIM of $-10.8$ m represents the accumulated efficiency impact over the identified affected vehicles and impact duration, rather than a travel-distance loss for each individual vehicle. Similarly, a target-lane TSIM of 4.6\% indicates increased TTC-based unsafe exposure within the identified impact region, rather than a direct increase in crash probability. \rev{Therefore, impact duration and affected-vehicle count are descriptive measures of the temporal and spatial extent identified by the proposed framework, rather than independent evidence of ground-truth impact ranges. TEIM and TSIM provide additional magnitude information that should be interpreted together with their bootstrap confidence intervals, median and interquartile ranges, Wilcoxon signed-rank results, the event-level distributions shown in Figs.~\ref{Fig9} and~\ref{Fig10}, and lane-specific traffic mechanisms.}

%%%%%%%%%%%%%%%%%%%%%%%%%%%%%%%%%%%%%%%%%%%%%%%%%%%%%%%%%%%%%%%%%%%%%%%%%%%%%%%%%%%%%%%%%%%%%%%%%%%%%
\section{Validation and sensitivity analysis}\label{Section6}

\subsection{Construction of matched no-lane-change control segments}\label{Section6.1}

{
To evaluate the false positive rate of the proposed method under natural traffic fluctuations, matched no-lane-change control segments were constructed from the ZTD trajectory data. These controls were designed to approximate traffic evolution in the absence of an \textit{SV} lane-changing maneuver and to provide a matched reference for comparison with real lane-change events. The construction of the matched no-lane-change controls included the following two steps.

\textit{Step 1: Candidate no-lane-change pseudo-event construction.}
Candidate no-lane-change pseudo-events were first extracted from lane-keeping vehicles in the ZTD trajectories. For each candidate, one lane-keeping vehicle was designated as the pseudo-SV at a pseudo-event time $t_0$. The lane occupied by the pseudo-SV was defined as the original lane, and the adjacent mainline lane was defined as the target lane. A candidate was retained only if the pseudo-SV remained in the same mainline lane throughout the observation window $[t_0-50~\mathrm{s},t_0+50~\mathrm{s}]$ and no vehicle changed into or out of the paired original and target lanes within the local spatiotemporal region $[t_0-50~\mathrm{s},t_0+50~\mathrm{s}]$$\times[x(t_0)-500~\mathrm{m},x(t_0)+500~\mathrm{m}]$, where $x(t_0)$ is the longitudinal position of the pseudo-SV at $t_0$.

\textit{Step 2: Traffic-state matching.}
To make the no-lane-change controls comparable to real lane-change events, each real event and each candidate control segment were described using the same local traffic-state vector. The vector included the SV or pseudo-SV speed, local mean speed, local density, front and rear gaps in the original and target lanes, and relative speeds to the front vehicles in the original and target lanes. For real lane-change events, all matching variables were evaluated at the identified lane-change start time $T_{SV}^{s}$. For no-lane-change candidates, the corresponding variables were evaluated at the pseudo-event time $t_0$. Local mean speed and density were calculated within the 1 km mainline window centered at the SV or pseudo-SV position (500 m upstream and 500 m downstream), with density expressed in vehicles per kilometre.

Traffic-state regimes were assigned using empirical tertiles of the pooled local traffic states from real lane-change events and candidate no-lane-change controls. Let $\bar{v}_{33}$ and $\bar{v}_{67}$ denote the 33rd and 67th percentiles of local mean speed, and let $\rho_{33}$ and $\rho_{67}$ denote the corresponding percentiles of local density. An event was assigned to the high-speed/low-density regime if $\bar{v}\geq\bar{v}_{67}$ and $\rho\leq\rho_{33}$, the low-speed/high-density regime if $\bar{v}\leq\bar{v}_{33}$ and $\rho\geq\rho_{67}$, and the intermediate regime otherwise.

For a real lane-change event $e$ and a no-lane-change candidate $c$, the standardized matching distance was defined as
\begin{equation}\label{EqNoLCMatchDistance}
D(e,c)=\sum_{m=1}^{M}\frac{|x_{e,m}-x_{c,m}|}{s_m}
\end{equation}
where $x_{e,m}$ and $x_{c,m}$ are the $m$-th traffic-state variables of the real event and the candidate control, respectively, and $s_m$ is the corresponding scaling factor. Speed and relative-speed differences were scaled by 10 km/h, local-density differences by 20 veh/km, and gap differences by 50 m. For each real lane-change event, eligible candidates were first restricted to no-lane-change segments from the same dataset and with the same original-to-target lane-pair direction. Matching was then performed only within the same traffic-state regime. The matched control for event $e$ was selected as
\begin{equation}\label{EqNoLCMatchedControl}
c^{*}(e)=\arg\min_{c\in\mathcal{C}_e}D(e,c)
\end{equation}
where $\mathcal{C}_e$ denotes the eligible candidate set restricted to the same dataset and lane-change direction, the same traffic-state regime, and candidates satisfying the caliper constraint $D(e,c)\leq5$. The validation subset used one-to-one nearest-neighbor matching without replacement. To obtain high-quality matched controls, candidate controls with $D(e,c)>5$ were excluded. Because $D(e,c)$ is the sum of standardized absolute differences over the nine matching variables ($M=9$), this threshold corresponds to an average standardized discrepancy of approximately $5/9=0.56$ per variable. In the retained matched subset, the mean standardized matching distance was 2.48, well below the caliper, indicating that the caliper provided a practical balance between local traffic-state comparability and sample retention. Using this procedure, 150 matched real-LC/no-LC event pairs were obtained from the 228 valid lane-change instances.

These matched no-lane-change control segments were then analyzed using the proposed framework in the same manner as real lane-change events. If the proposed method identified an affected vehicle in a no-lane-change control segment, the detection was treated as a false positive. This matched subset was used as the common sample basis for the validation metrics, parameter sensitivity analysis, and method comparison in Sections~\ref{Section6.2}--\ref{Section6.4}. \rev{The resulting full-sample findings therefore constitute an in-sample assessment under the matched-control design. To further examine whether using the same matched event pairs materially affects the conclusions reported in Sections~\ref{Section6.2}--\ref{Section6.4}, an event-level 5-fold cross-validation analysis was conducted, with the results provided in Appendix A.}
}

\subsection{Validation metrics}\label{Section6.2}

{
To evaluate whether the framework can reduce the misattribution of background traffic fluctuations to lane-change-associated disturbances, two complementary metrics were defined for the matched no-lane-change controls and the paired real lane-change events.

Let $\mathcal{V}_c$ denote the set of candidate following vehicles in no-lane-change control segment $c$, and let $z_{c,i}=1$ if following vehicle $i\in\mathcal{V}_c$ is classified as affected by the framework and $z_{c,i}=0$ otherwise. Because no SV lane change occurs in a control segment, any candidate following vehicle classified as affected was counted as a false-positive identification. The no-lane-change false positive rate (No-LC FPR) is defined as:
\begin{equation}\label{EqNoLCFPR}
\mathrm{FPR}^{\mathrm{NoLC}}=
\frac{\sum_{c=1}^{N_C}\sum_{i\in\mathcal{V}_c}z_{c,i}}
{\sum_{c=1}^{N_C}|\mathcal{V}_c|}
\end{equation}
where $N_C$ is the number of matched no-lane-change control segments.

For real lane-change events, let $\mathcal{V}_e$ be the candidate following-vehicle set for event $e$, and let $z_{e,i}=1$ indicate that vehicle $i\in\mathcal{V}_e$ is classified as affected by the same framework. The real-LC affected-vehicle identification rate (Real-LC AIR) is defined as:
\begin{equation}\label{EqRealLCAIR}
\mathrm{AIR}^{LC}=
\frac{\sum_{e=1}^{N_L}\sum_{i\in\mathcal{V}_e}z_{e,i}}
{\sum_{e=1}^{N_L}|\mathcal{V}_e|}
\end{equation}
where $N_L$ is the number of real lane-change events in the corresponding validation set. Because objective ground-truth affected/unaffected labels are unavailable for real lane-change events, $\mathrm{AIR}^{LC}$ is not interpreted as recall, sensitivity, or accuracy. It quantifies the proportion of candidate following vehicles for which the framework identifies affected states under real lane-change conditions.

These two metrics were interpreted jointly in the subsequent sensitivity and method-comparison analyses. A lower No-LC FPR indicates fewer false detections under matched no-lane-change conditions, whereas a higher Real-LC AIR indicates stronger retention of affected-vehicle identification under real lane-change conditions. However, either metric alone can be misleading: a very low No-LC FPR may be obtained by suppressing impact identification under real lane-change conditions, whereas a very high Real-LC AIR may be obtained by over-identifying affected vehicles. Therefore, parameter settings and comparison methods were evaluated based on the trade-off between No-LC FPR and Real-LC AIR rather than by either metric alone.
}

\subsection{Sensitivity analysis}\label{Section6.3}

\subsubsection{Sensitivity analysis of key affected-vehicle identification parameters}\label{Section6.3.1}

\begin{figure}[H]
  \centering
  \begin{tabular}{@{}c@{\hspace{0.025\linewidth}}c@{\hspace{0.025\linewidth}}c@{}}
  \begin{subfigure}{0.31\linewidth}
        \centering
        \includegraphics[width=\linewidth]{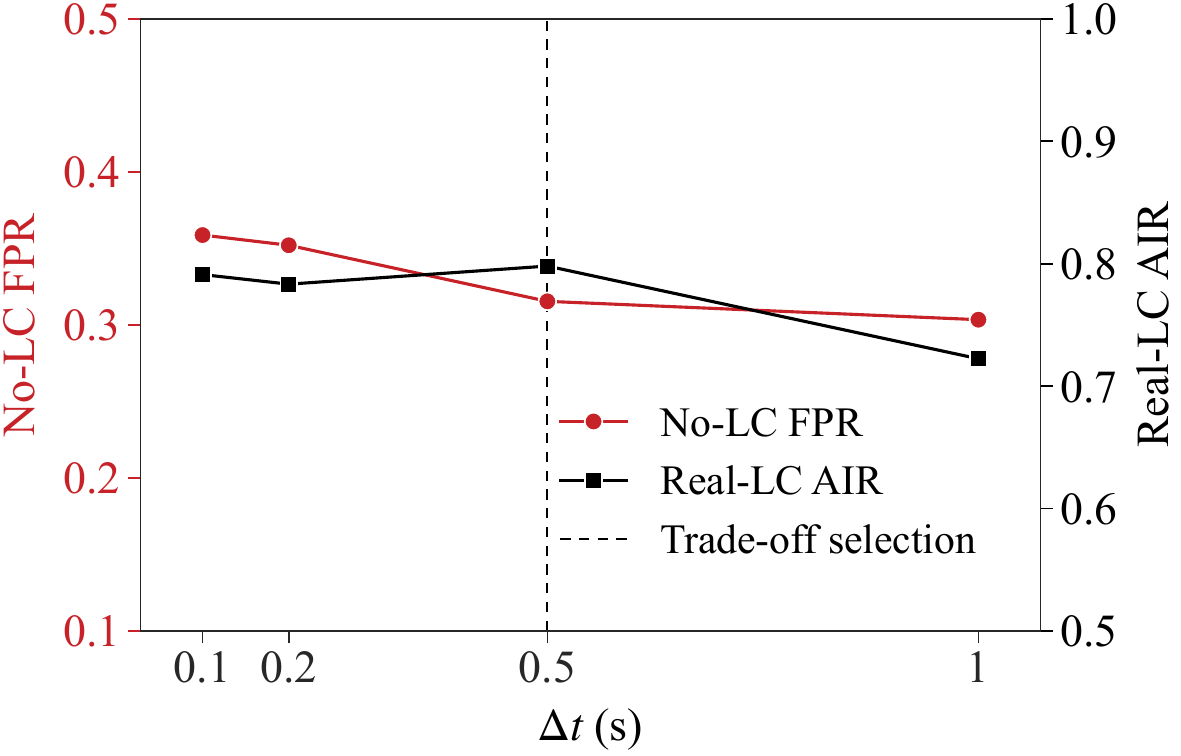}
        \caption{}
        \label{FigSensitivityA}
  \end{subfigure}
  &
  \begin{subfigure}{0.31\linewidth}
        \centering
        \includegraphics[width=\linewidth]{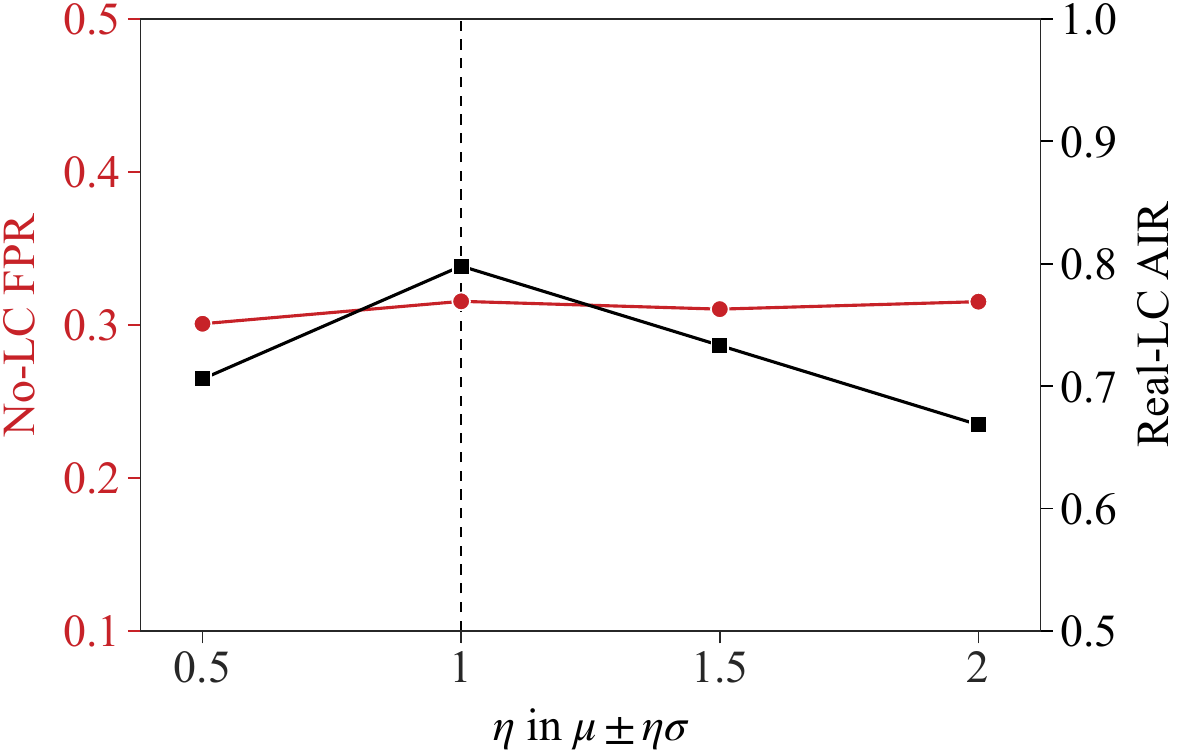}
        \caption{}
        \label{FigSensitivityB}
  \end{subfigure}
  &
  \begin{subfigure}{0.31\linewidth}
        \centering
        \includegraphics[width=\linewidth]{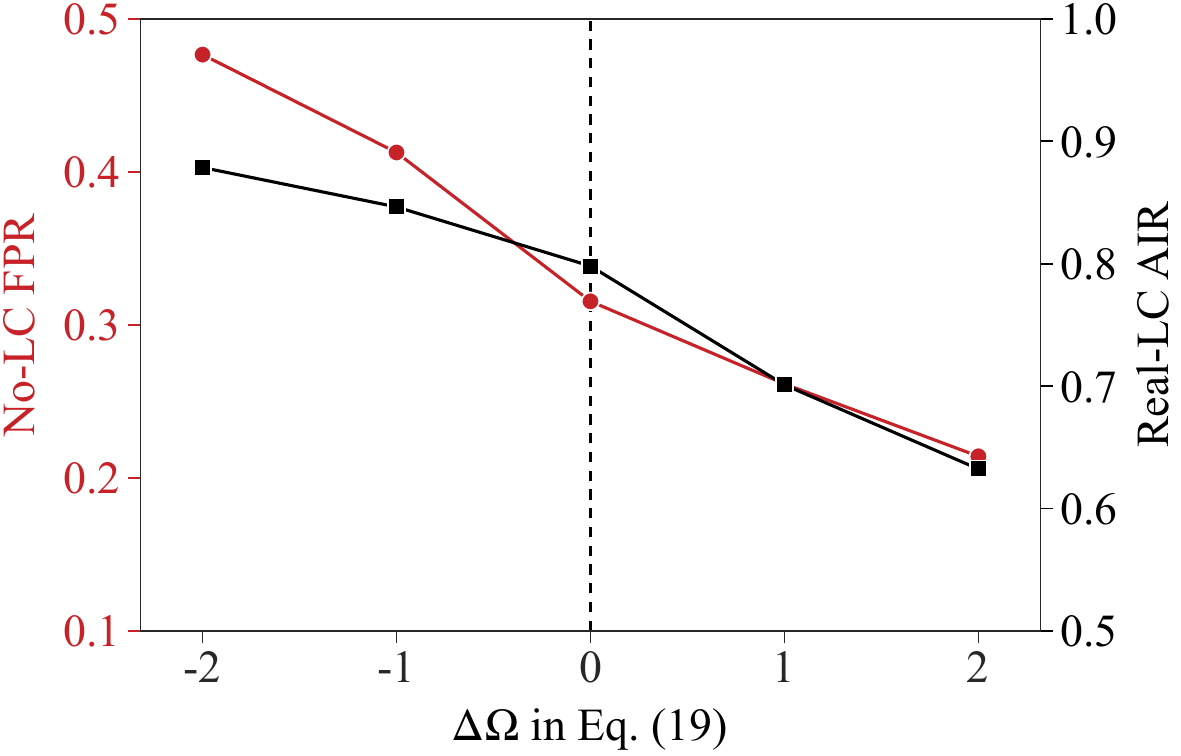}
        \caption{}
        \label{FigSensitivityC}
  \end{subfigure}
  \\[0.6em]
  \begin{subfigure}{0.31\linewidth}
        \centering
        \includegraphics[width=\linewidth]{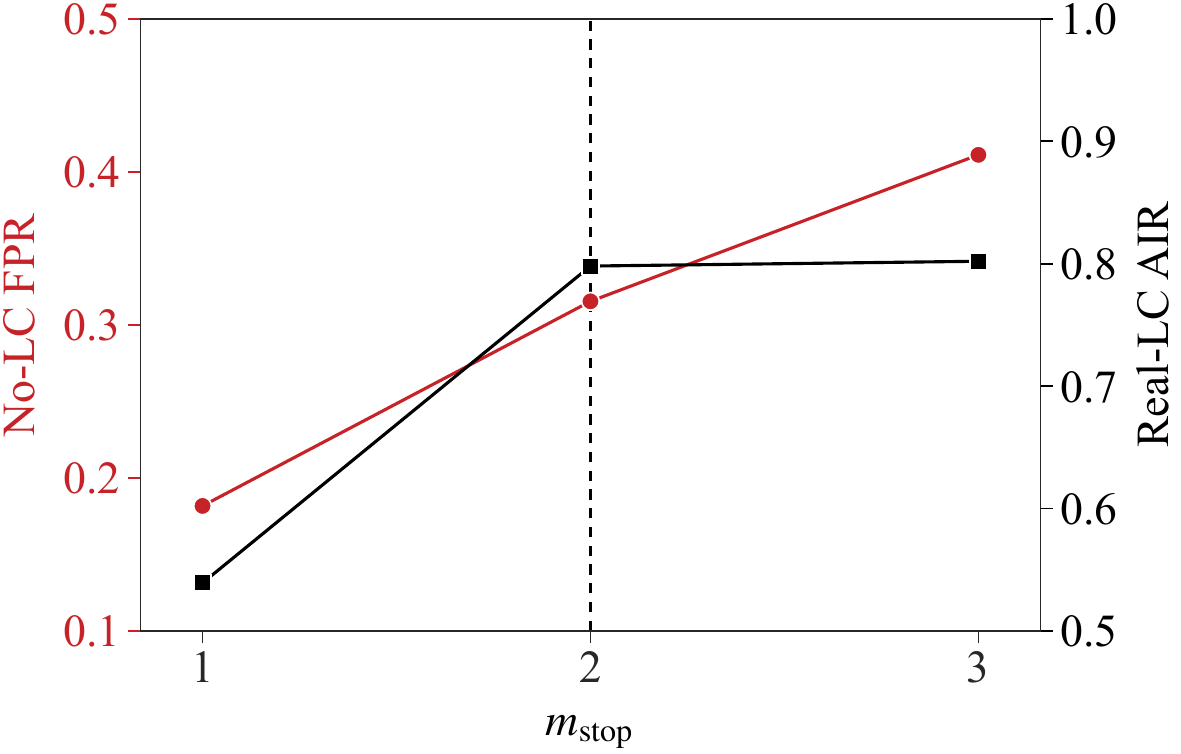}
        \caption{}
        \label{FigSensitivityD}
  \end{subfigure}
  &
  \begin{subfigure}{0.31\linewidth}
        \centering
        \includegraphics[width=\linewidth]{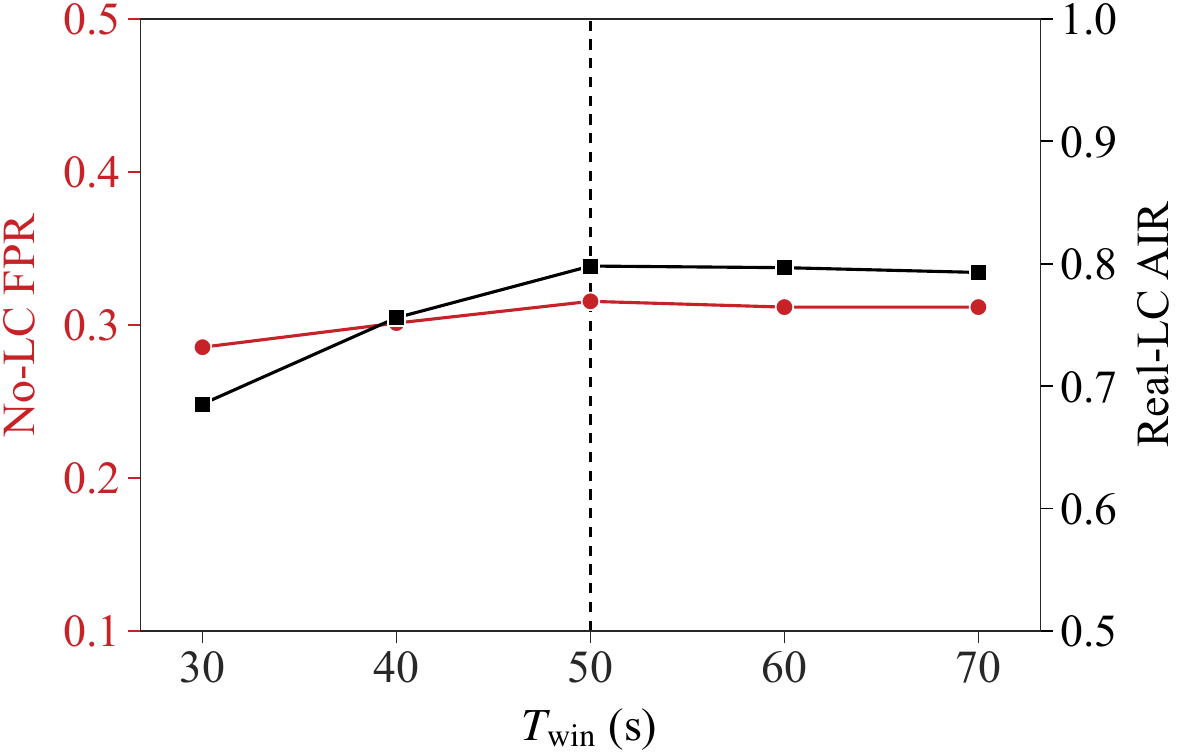}
        \caption{}
        \label{FigSensitivityE}
  \end{subfigure}
  &
  \begin{minipage}{0.31\linewidth}
        \vspace{0pt}
  \end{minipage}
  \end{tabular}
  \caption{Sensitivity of key affected-vehicle identification parameter settings: (a) $\Delta t$; (b) $\eta$ in $\mu \pm \eta\sigma$; (c) $\Delta\Omega$ in the run-length persistence filter; (d) $m_{\mathrm{stop}}$; and (e) $T_{\mathrm{win}}$}
  \label{FigSensitivity}
\end{figure}

This section examines the sensitivity of key affected-vehicle identification parameters, including $\Delta t$, $\eta$ in the $\mu\pm \eta\sigma$ threshold, $\Delta\Omega$ in the run-length persistence filter, $m_{\mathrm{stop}}$ in Eq.~(\ref{Eq28}), and $T_{\mathrm{win}}$ in Criterion 3 of Section~\ref{Section4.2}. The No-LC FPR and Real-LC AIR defined in Section~\ref{Section6.2} were used to evaluate how each parameter affects false-positive control and retention of impact identification under real lane-change conditions. For the run-length persistence filter, the threshold was perturbed as $\Omega_{\ell,i}^{f,\Delta}=\max(0,\Omega_{\ell,i}^{f\ast}+\Delta\Omega)$, where $\Delta\Omega\in\{-2,-1,0,+1,+2\}$ and $\Delta\Omega=0$ corresponds to the formal setting in Eq.~(\ref{Eq19}). The corresponding No-LC FPR--Real-LC AIR trade-off curves are summarized in Fig.~\ref{FigSensitivity}, where the black vertical dashed lines indicate the parameter settings adopted based on the No-LC FPR--Real-LC AIR trade-off.

For $\Delta t$, overly fine aggregation intervals increased No-LC FPR, whereas the 1.0 s interval reduced No-LC FPR but also reduced Real-LC AIR, indicating that excessive temporal smoothing suppresses affected-vehicle identification under real lane-change conditions. The adopted $\Delta t=0.5$ s produced the highest Real-LC AIR among the tested settings while keeping No-LC FPR moderate. For $\eta$ in the $\mu\pm \eta\sigma$ threshold, $\eta=1.0$ produced the highest Real-LC AIR (0.798), while its No-LC FPR (0.315) remained close to the other tested thresholds. Although $\eta=0.5$ slightly reduced No-LC FPR, it lowered Real-LC AIR to 0.706; larger $\eta$ values also reduced Real-LC AIR without meaningful No-LC FPR improvement. For the run-length persistence filter, increasing $\Delta\Omega$ imposed a stricter persistence requirement. As a result, No-LC FPR decreased monotonically from 0.477 to 0.214, but Real-LC AIR also decreased from 0.879 to 0.632. This pattern is consistent with the role of the filter: stricter thresholds remove more short abnormal sequences that may reflect background traffic fluctuations, but they can also suppress weaker or shorter impact identification under real lane-change conditions. Thus, $\Delta\Omega=0$, i.e., using the longest pre-event abnormal sequence as the threshold, provides an interpretable balance between filtering short background fluctuations and retaining impact identification under real lane-change conditions. For $T_{\mathrm{win}}$, the 30 s window reduced No-LC FPR but substantially lowered Real-LC AIR, whereas longer windows did not provide a simultaneous improvement. These results support $T_{\mathrm{win}}=50$ s as a balanced setting among the tested window lengths.

The stopping parameter $m_{\mathrm{stop}}$ produced a clear No-LC FPR--Real-LC AIR trade-off. Setting $m_{\mathrm{stop}}=1$ was too conservative, reducing No-LC FPR but also strongly suppressing Real-LC AIR. Setting $m_{\mathrm{stop}}=3$ produced only a marginal Real-LC AIR increase compared with $m_{\mathrm{stop}}=2$, while increasing No-LC FPR substantially. Thus, $m_{\mathrm{stop}}=2$ provides the most balanced setting among the tested alternatives. Taken together, these results support the adopted parameter settings of $\Delta t=0.5$ s, $\eta=1.0$, $\Delta\Omega=0$, $m_{\mathrm{stop}}=2$, and $T_{\mathrm{win}}=50$ s.

\subsubsection{Sensitivity analysis of the TTC threshold}\label{Section6.3.2}

The TTC threshold was examined separately because it affects only the unsafe-exposure calculation in SIM and TSIM, rather than the affected-vehicle identification procedure evaluated by No-LC FPR and Real-LC AIR. For this analysis, the TTC threshold was varied from 2.5 s to 6.5 s around the selected value of 4.5 s. Table~\ref{tableTTCThresholdRobustness} summarizes the resulting mean safety-impact metrics \rev{together with the proportions of events with positive SIM/TSIM values}.

\begin{table}[H]
  \centering
  \captionsetup{font=small}
  \caption{\rev{Sensitivity of safety-impact measures to the TTC threshold}}
  \label{tableTTCThresholdRobustness}
  {
  \small
  \setlength{\tabcolsep}{3.0pt}
  \resizebox{\textwidth}{!}{%
  \begin{tabular}{@{}lcccccccccc@{}}
    \toprule[1.0pt]
    \multirow{2}{*}{\textbf{Metric}} &
    \multicolumn{2}{c}{\textbf{2.5 s}} &
    \multicolumn{2}{c}{\textbf{3.5 s}} &
    \multicolumn{2}{c}{\textbf{4.5 s}} &
    \multicolumn{2}{c}{\textbf{5.5 s}} &
    \multicolumn{2}{c}{\textbf{6.5 s}} \\
    \cmidrule(lr){2-3}\cmidrule(lr){4-5}\cmidrule(lr){6-7}\cmidrule(lr){8-9}\cmidrule(lr){10-11}
    & \rev{\textbf{Mean}} & \rev{\textbf{Events $>0$ (\%)}} 
    & \rev{\textbf{Mean}} & \rev{\textbf{Events $>0$ (\%)}} 
    & \rev{\textbf{Mean}} & \rev{\textbf{Events $>0$ (\%)}} 
    & \rev{\textbf{Mean}} & \rev{\textbf{Events $>0$ (\%)}} 
    & \rev{\textbf{Mean}} & \rev{\textbf{Events $>0$ (\%)}} \\
    \midrule[0.5pt]
    $TFV_1$ SIM (\%) & 0.7 & \rev{14.0} & 1.2 & \rev{25.4} & 1.6 & \rev{28.1} & 2.0 & \rev{29.4} & 2.4 & \rev{35.5} \\
    $FV_1$ SIM (\%) & 0.1 & \rev{7.0} & 0.3 & \rev{10.1} & 0.5 & \rev{17.1} & 1.0 & \rev{23.2} & 1.6 & \rev{28.1} \\
    Target-lane TSIM (\%) & 5.5 & \rev{16.2} & 6.6 & \rev{20.6} & 4.6 & \rev{25.4} & 4.2 & \rev{28.9} & 6.2 & \rev{32.9} \\
    Original-lane TSIM (\%) & 0.4 & \rev{8.8} & 1.1 & \rev{14.0} & 2.1 & \rev{18.4} & 2.8 & \rev{23.7} & 2.8 & \rev{26.3} \\
    Global TSIM (\%) & 5.2 & \rev{22.4} & 7.0 & \rev{28.1} & 6.6 & \rev{34.6} & 7.2 & \rev{40.8} & 8.8 & \rev{47.8} \\
    \bottomrule[1.0pt]
  \end{tabular}
  }
  }
  \vspace{0.3em}

  \begin{minipage}{\textwidth}
    \footnotesize
    \rev{\textit{Note:} Events $>0$ (\%) denotes the percentage of lane-change events for which the corresponding event-level safety-impact magnitude (SIM or TSIM) is positive.}
  \end{minipage}
\end{table}

\rev{Table~\ref{tableTTCThresholdRobustness} shows that the proportion of events with positive SIM/TSIM values increased across all five metrics as the TTC threshold increased. This empirical pattern is consistent with larger TTC thresholds classifying a broader range of following situations as unsafe. However, the positive-event proportions remained below 50\% for all metrics and threshold settings. At the selected threshold of 4.5 s, they ranged from 17.1\% for $FV_1$ SIM to 34.6\% for global TSIM. Therefore, positive mean SIM/TSIM values should be interpreted as average tendencies toward increased TTC-based unsafe exposure rather than evidence of a uniform increase across all lane-change events. The two types of results provide complementary information: mean safety-impact magnitudes characterize the average direction and magnitude of the estimated effects, whereas the proportion of positive event-level values reveals substantial heterogeneity across lane-change events. Notably, this sensitivity analysis evaluates how the resulting safety-impact measures vary under alternative TTC thresholds, whereas the 4.5 s threshold was selected with reference to the TTC-threshold literature discussed in Section~\ref{Section3.2} below Eq.~(\ref{Eq38}) \citep{rahman2018longitudinal,rahman2019safety,lu2021performance}.}

\subsection{Method comparison}\label{Section6.4}

{
This section compares five affected-vehicle identification methods using the No-LC FPR and Real-LC AIR defined in Section~\ref{Section6.2}. The comparison focuses on whether each method can reduce false-positive detections in no-lane-change controls while retaining affected-vehicle identification under real lane-change conditions.

The five methods are defined as follows:
\begin{itemize}
  \item \textbf{Method 1: Proposed method.} This is the full framework described in Sections~\ref{Section2} and \ref{Section3}.
  \item \textbf{Method 2: Proposed method with a fixed 20 m/s free-speed reference.} This variant replaces the dynamic TLV/LV speed reference in Eq.~(\ref{Eq8}) with a fixed free-speed reference of 20 m/s.
  \item \textbf{Method 3: Proposed method without pre-event TDB envelope.} This variant removes the pre-event TDB fluctuation envelope in Eq.~(\ref{Eq17}), so any non-zero TDB after the demarcation time is treated as a potential disturbance before the run-length decision.
  \item \textbf{Method 4: Proposed method without run-length filter.} This variant removes the persistence filter in Eq.~(\ref{Eq19}), so short abnormal sequences are no longer suppressed.
  \item \textbf{Method 5: Method of \citet{he2023impact}.} This method implements the dynamic minimum-spacing method of \citet{he2023impact}, in which lane-change impact is identified from the temporal pattern of Newell-model minimum spacing along the upstream vehicle chain.
\end{itemize}
}

\begin{table}[H]
  \centering
  \captionsetup{font=small}
  \caption{Comparison of affected-vehicle identification methods}
  \label{tableMethodComparison}
  {
  \small
  \setlength{\tabcolsep}{4.0pt}
  \resizebox{\textwidth}{!}{%
  \begin{tabular}{@{}lllllll@{}}
    \toprule[1.0pt]
    \textbf{Method} & \shortstack[l]{\textbf{No-LC}\\\textbf{FPR}} & \shortstack[l]{\textbf{Real-LC}\\\textbf{AIR}} & \shortstack[l]{\textbf{Target-lane}\\\textbf{affected vehicles}} & \shortstack[l]{\textbf{Original-lane}\\\textbf{affected vehicles}} & \shortstack[l]{\textbf{Target-lane}\\\textbf{impact duration (s)}} & \shortstack[l]{\textbf{Original-lane}\\\textbf{impact duration (s)}} \\
    \midrule[0.5pt]
    Proposed method & 0.315 & 0.798 & 6.49 & 5.73 & 27.01 & 26.60 \\
    Proposed with fixed 20 m/s reference & 0.374 & 0.606 & 4.87 & 4.41 & 25.01 & 24.87 \\
    Proposed without pre-event TDB envelope & 0.880 & 0.997 & 8.08 & 7.20 & 36.10 & 36.97 \\
    Proposed without run-length filter & 0.844 & 0.996 & 8.06 & 7.19 & 34.81 & 35.74 \\
    He et al. (2023) method & 0.879 & 0.950 & 7.71 & 6.84 & 39.01 & 37.28 \\
    \bottomrule[1.0pt]
  \end{tabular}}}
\end{table}

Table~\ref{tableMethodComparison} shows that, among the tested methods, the proposed method achieved the lowest No-LC FPR of 0.315, with a Real-LC AIR of 0.798. This combination represents a comparatively favorable trade-off between reducing false-positive identification under matched no-lane-change conditions and retaining affected-vehicle identification under real lane-change conditions. \rev{Nevertheless, the No-LC FPR of 0.315 indicates that a non-negligible level of false-positive identification remains under the matched no-lane-change controls. Thus, the results indicate that the proposed method reduces the misattribution of background traffic fluctuations relative to the tested alternatives, but does not eliminate it.} The affected-vehicle counts and impact durations in Table~\ref{tableMethodComparison} are reported as mean values across the evaluated real lane-change events. Replacing the dynamic TLV/LV reference in Eq.~(\ref{Eq8}) with a fixed 20 m/s free-speed reference increased No-LC FPR to 0.374 and reduced Real-LC AIR to 0.606. This indicates that the dynamic reference contributes to both false-positive control and retention of impact identification under real lane-change conditions. In contrast, removing the pre-event TDB envelope (Eq.~(\ref{Eq17})) or the run-length filter (Eq.~(\ref{Eq19})) substantially increased No-LC FPR to 0.880 and 0.844, respectively. These two ablations also enlarged the estimated impact range in both lanes. Removing the pre-event TDB envelope increased the target/original affected-vehicle counts from 6.49/5.73 to 8.08/7.20 and the corresponding impact durations from 27.01/26.60 s to 36.10/36.97 s. Removing the run-length filter produced similar increases, with target/original affected-vehicle counts of 8.06/7.19 and impact durations of 34.81/35.74 s. These results indicate that the pre-event TDB envelope and run-length filter are the main components reducing false-positive affected-vehicle identification caused by background traffic oscillations.

The He et al. method also achieved a high Real-LC AIR (0.950), but its No-LC FPR was high (0.879), and it identified broader lane-specific impact ranges than the proposed method, with target/original affected-vehicle counts of 7.71/6.84 versus 6.49/5.73 and impact durations of 39.01/37.28 s versus 27.01/26.60 s. Thus, the minimum-spacing pattern can identify many candidate affected vehicles under real lane-change conditions, but without an explicit pre-event fluctuation envelope and run-length noise filter, it is more sensitive to background traffic oscillations and tends to produce broader upstream impact ranges. Overall, the comparison suggests that the dynamic reference, pre-event fluctuation envelope, and run-length filter play complementary roles, with the latter two being particularly important for suppressing false positives caused by normal traffic fluctuations.

\subsection{Robustness to perturbations in the identified lane-change start time}\label{Section6.5}

Because $T_{SV}^{s}$ is inferred from lateral trajectory patterns rather than directly observed, we examined whether the main spatial and temporal impact-range metrics are robust to plausible perturbations in $T_{SV}^{s}$. We perturbed $T_{SV}^{s}$ by $-3$ s, $-2$ s, $-1$ s, $+1$ s, and $+2$ s relative to the identified start time, while the 0 s column represents the identified start-time setting. The demarcation times $T_{\ell,i}^s$ were shifted consistently with $T_{SV}^{s}$. Positive perturbations were constrained to remain before $T_{SV}^{lane}$. The robustness analysis focused on four range-related metrics: target-lane impact duration, original-lane impact duration, target-lane number of affected vehicles, and original-lane number of affected vehicles.

\begin{table}[H]
  \centering
  \captionsetup{font=small}
  \caption{Start-time perturbation robustness results}
  \label{tableStartTimeRobustness}
  {
  \small
  \setlength{\tabcolsep}{4.5pt}
  \begin{tabular}{@{}lllllll@{}}
    \toprule[1.0pt]
    \textbf{Metric} & \textbf{-3 s} & \textbf{-2 s} & \textbf{-1 s} & \textbf{0 s} & \textbf{+1 s} & \textbf{+2 s} \\
    \midrule[0.5pt]
    Target-lane impact duration (s) & 27.1 & 25.7 & 24.7 & 23.8 & 22.4 & 20.5 \\
    Original-lane impact duration (s) & 28.6 & 27.1 & 26.1 & 25.0 & 23.0 & 21.8 \\
    Target-lane affected vehicles & 5.7 & 5.6 & 5.5 & 5.6 & 5.2 & 4.8 \\
    Original-lane affected vehicles & 5.5 & 5.4 & 5.4 & 5.3 & 4.9 & 4.7 \\
    \bottomrule[1.0pt]
  \end{tabular}
  }
\end{table}

Table~\ref{tableStartTimeRobustness} shows two main patterns. First, the estimated impact range remained within a comparable magnitude under all tested perturbations: the impact-duration metrics varied from 20.5 s to 28.6 s, and the affected-vehicle counts remained around five vehicles in both lanes. Second, the changes followed a systematic direction. Earlier start-time shifts generally enlarged the estimated impact range, whereas later shifts reduced it. This pattern is consistent with the demarcation-time structure of the proposed method. Because affected status is evaluated only after each vehicle-specific demarcation time, delaying $T_{SV}^{s}$ shifts the demarcation times of upstream following vehicles later and can exclude early lane-change disturbance responses from the affected-status evaluation. As a result, some vehicles may be classified as unaffected, and the retained affected intervals may become shorter. The largest reductions occurred under the $+2$ s perturbation, where the target-lane and original-lane impact durations decreased to 20.5 s and 21.8 s, and the corresponding affected-vehicle counts decreased to 4.8 and 4.7. Even under the $+2$ s perturbation, which produced the largest reductions among the tested shifts, the estimated impact range remained close to that obtained using the identified start time. This indicates that the identified temporal and spatial impact ranges were not highly sensitive to the tested start-time perturbations.

%%%%%%%%%%%%%%%%%%%%%%%%%%%%%%%%%%%%%%%%%%%%%%%%%%%%%%%%%%%%%%%%%%%%%%%%%%%%%%%%%%%%%%%%%%%%%%%%%%%%%
\section{Concluding remarks}\label{Section7}
In this study, we developed a trajectory-based framework for quantifying the spatiotemporal impacts of individual discretionary lane changes. Using vehicle trajectories before and after the lane-changing maneuver, the framework identifies affected upstream following vehicles, determines their impact durations, and quantifies the corresponding efficiency and safety impact magnitudes. By integrating both spatial (i.e., the number of affected vehicles) and temporal (i.e., the duration of impact) dimensions, we further propose two aggregate indicators—TEIM and TSIM—to characterize efficiency- and safety-related lane-change impact magnitudes. \rev{The framework is intended to identify lane-change-associated disturbances and to reduce, rather than eliminate, the misattribution of background traffic fluctuations.} \rev{To this end, it combines} local dynamic references, vehicle-specific pre-event TDB envelopes, run-length persistence filtering, and baseline-excess efficiency and safety corrections. Moreover, the proposed metrics allow signed characterization of beneficial and adverse impacts on traffic efficiency and safety.

The empirical results show that the framework identifies lane-change-associated disturbances among upstream following vehicles in both lanes, but the identified patterns differ between the target and original lanes. In the target lane, SV insertion generally compresses the available gap and is associated with average efficiency loss and increased TTC-based unsafe exposure for following vehicles. In the original lane, SV departure may release space and correspond to efficiency gains for some following vehicles, while gap-closing responses can still increase unsafe exposure. \rev{The bootstrap confidence intervals further indicate that the estimated mean impact duration and affected-vehicle count have lower sampling uncertainty than the aggregate TEIM and TSIM magnitudes}, which are more sensitive to heterogeneity across lane-change events, positive-negative cancellation, and lane-specific driving responses.

The matched no-lane-change controls and method comparison provide comparative evidence that the proposed framework reduces false-positive affected-vehicle identification relative to the tested alternatives. The ablation results further indicate that the pre-event TDB envelope and run-length filter play important roles in this reduction. \rev{Nevertheless, the proposed method yields a No-LC FPR of 0.315, indicating that a non-negligible level of false-positive identification remains.} The parameter sensitivity analysis, TTC-threshold analysis, and start-time perturbation robustness check further characterize how the results vary under the tested settings.

Future research can extend the proposed framework in two directions. First, the empirical sample was constructed using an isolated-event criterion that excluded lane changes occurring within $\pm 50$ s of another LC event. This design reduces confounding from overlapping lane-change interactions and makes the upstream impact of a single discretionary lane change more identifiable. However, it may also underrepresent dense, highly interactive traffic conditions. Therefore, future work can examine and extend the framework in more complex traffic scenarios, such as ramp merging and diverging areas, weaving sections, and dense multi-lane-change interactions, where the impacts of multiple lane changes may overlap and accumulate. \rev{Second, ground-truth affected-vehicle labels, impact ranges, and impact magnitudes are unavailable in the ZTD dataset. The framework also does not construct an event-specific no-lane-change counterfactual. Accordingly, the identified impact ranges and magnitudes should be interpreted as measures of lane-change-associated disturbances rather than exact causal effects.} Future work can further evaluate the framework using independent trajectory datasets with richer event annotations or controlled simulation scenarios in which affected-vehicle labels, impact ranges, and impact magnitudes are known by design.

%%%%%%%%%%%%%%%%%%%%%%%%%%%%%%%%%%%%%%%%%%%%%%%%%%%%%%%%%%%%%%%%%%%%%%%%%%%%%%%%%%%%%%%%%%%%%%%%%%%%%

%% BIBLIOGRAPHY
%%%%%%%%%%%%%%%%%%%%%%%%%%%%%%%%%%%%%%%%%%%%%%%%%%%%%%%%%%%%%

\section*{Declaration of Competing Interest}
The authors declare that they have no known competing financial interests or personal relationships that could have appeared to influence the work reported in this paper.

\section*{Acknowledgments}
This work was supported by Sichuan Province Science and Technology Support Program under Grant No. 2025YFHZ0193, National Natural Science Foundation of China (NSFC) under Grants 52572362 and 72671151, and Doctoral Innovation Fund Program of Southwest Jiaotong University under Grant No. CX2025YB13. 

\appendix
\setcounter{table}{0}
\renewcommand{\thetable}{A.\arabic{table}}
\section{\rev{Event-level 5-fold cross-validation analysis}}\label{AppendixA}

\rev{To examine whether the reuse of the same matched sample affects the parameter-selection and method-comparison conclusions, an event-level 5-fold cross-validation analysis was conducted using the 150 matched real-LC/no-LC event pairs. The matched pairs were divided into five folds, with each real-LC/no-LC pair retained within the same fold. In each iteration, the parameter configuration was determined using four training folds according to the same No-LC FPR--Real-LC AIR trade-off used in Section~\ref{Section6.3.1} and was then evaluated on the held-out fold. Across all five iterations, the same parameter configuration was consistently selected: $\Delta t=0.5$ s, $\eta=1.0$, $\Delta\Omega=0$, $m_{\mathrm{stop}}=2$, and $T_{\mathrm{win}}=50$ s. As summarized in Table~\ref{tableAppendixCVParameterSelection}, the held-out mean No-LC FPR and Real-LC AIR were 0.315 (SD = 0.051) and 0.797 (SD = 0.030), respectively, which closely matched the corresponding full-sample values of 0.315 and 0.798.}

\rev{The same fold assignment was used for the held-out comparison of the five affected-vehicle identification methods examined in Section~\ref{Section6.4}. As shown in Table~\ref{tableAppendixCVMethodComparison}, the rankings of the five methods based on the mean No-LC FPR and Real-LC AIR across the five held-out folds were consistent with the corresponding full-sample rankings. The proposed method achieved the lowest held-out mean No-LC FPR among the tested methods, with a mean Real-LC AIR of 0.797. These results indicate that the selected parameter configuration and the main comparative pattern were stable across the 5-fold partitions.}

\begin{table}[H]
  \centering
  \captionsetup{font=small}
  \caption{\rev{Full-sample and 5-fold held-out results for parameter selection}}
  \label{tableAppendixCVParameterSelection}
  {
  \small
  \setlength{\tabcolsep}{6pt}
  \begin{tabularx}{0.70\textwidth}{@{}>{\raggedright\arraybackslash}X*{3}{>{\centering\arraybackslash}X}@{}}
    \toprule[1.0pt]
    \rev{\multirow{2}{*}{\textbf{Metric}}} & \rev{\multirow[c]{2}{*}{\raisebox{-0.5ex}{\shortstack{\textbf{Full sample}\\\textbf{(in-sample)}}}}} & \multicolumn{2}{c}{\rev{\textbf{Held-out folds}}} \\
    \cmidrule(lr){3-4}
     & & \rev{\textbf{Mean}} & \rev{\textbf{SD}} \\
    \midrule[0.5pt]
    \rev{No-LC FPR} & \rev{0.315} & \rev{0.315} & \rev{0.051} \\
    \rev{Real-LC AIR} & \rev{0.798} & \rev{0.797} & \rev{0.030} \\
    \bottomrule[1.0pt]
  \end{tabularx}
  \par\smallskip
  \begin{minipage}{0.70\textwidth}
    \raggedright\footnotesize
    \rev{\textit{Note:} The held-out mean and standard deviation (SD) are calculated across the five held-out folds. The full-sample values correspond to the in-sample results reported in Section~\ref{Section6.3.1}.}
  \end{minipage}
  }
\end{table}

\begin{table}[H]
  \centering
  \captionsetup{font=small}
  \caption{\rev{Full-sample and 5-fold held-out comparison of affected-vehicle identification methods}}
  \label{tableAppendixCVMethodComparison}
  {
  \small
  \setlength{\tabcolsep}{2pt}
  \begin{tabularx}{\textwidth}{@{}>{\raggedright\arraybackslash}p{0.35\textwidth}*{6}{>{\centering\arraybackslash}X}@{}}
    \toprule[1.0pt]
    \rev{\multirow{3}{*}{\textbf{Method}}} & \multicolumn{3}{c}{\rev{\textbf{No-LC FPR}}} & \multicolumn{3}{c}{\rev{\textbf{Real-LC AIR}}} \\
    \cmidrule(lr){2-4}\cmidrule(lr){5-7}
     & \rev{\multirow[c]{2}{*}{\raisebox{-0.5ex}{\shortstack{\textbf{Full sample}\\\textbf{(in-sample)}}}}} & \multicolumn{2}{c}{\rev{\textbf{Held-out folds}}} & \rev{\multirow[c]{2}{*}{\raisebox{-0.5ex}{\shortstack{\textbf{Full sample}\\\textbf{(in-sample)}}}}} & \multicolumn{2}{c}{\rev{\textbf{Held-out folds}}} \\
    \cmidrule(lr){3-4}\cmidrule(lr){6-7}
     & & \rev{\textbf{Mean}} & \rev{\textbf{SD}} & & \rev{\textbf{Mean}} & \rev{\textbf{SD}} \\
    \midrule[0.5pt]
    \rev{Proposed method} & \rev{0.315} & \rev{0.315} & \rev{0.051} & \rev{0.798} & \rev{0.797} & \rev{0.030} \\
    \rev{Proposed with fixed 20 m/s reference} & \rev{0.374} & \rev{0.373} & \rev{0.049} & \rev{0.606} & \rev{0.605} & \rev{0.074} \\
    \rev{Proposed without pre-event TDB envelope} & \rev{0.880} & \rev{0.881} & \rev{0.018} & \rev{0.997} & \rev{0.997} & \rev{0.002} \\
    \rev{Proposed without run-length filter} & \rev{0.844} & \rev{0.845} & \rev{0.021} & \rev{0.996} & \rev{0.996} & \rev{0.002} \\
    \rev{He et al. (2023) method} & \rev{0.879} & \rev{0.879} & \rev{0.030} & \rev{0.950} & \rev{0.950} & \rev{0.016} \\
    \bottomrule[1.0pt]
  \end{tabularx}
  \par\smallskip
  \begin{minipage}{\textwidth}
    \raggedright\footnotesize
    \rev{\textit{Note:} For each method, the held-out mean and SD are calculated across the five held-out folds. The full-sample values correspond to the in-sample results reported in Table~\ref{tableMethodComparison}.}
  \end{minipage}
  }
\end{table}

% The authors declare no potential conflict of interest.

% \pagebreak

% \end{appendices}

%%%%%%%%%%%%%%%%%%%%%%%%%%%%%%%%%%%%%%%%%%%%%%%%%%%%%%%%%%%%%%%%%%%%%%%%%%%%%%%%%%%%%%%%%%%%%%%%%%%%%%%%%%%%%%%%%%%%%%%%%%%%%%%%%%%%%%%%%%%%%%%%%%%%%%%%%%%%%%%%%%%%%%%%%%%%%%%%%%%%%%%%%%%%%%%%%%%%%%%%%%%%%%%%%%%%%%%%%%%%%%%%%%%%%%%%%%%%%%%%%%%%%%%%%%%%%%%%%%%%%%%%%%%%%%%%%%%%%%%%%%%%%%%%%%%%%%%%%%%%%%%%%%%%%%%%%%%%
% \pagebreak
%%%%%%%%%%%%%%%%%%%%%%%%%%%%%%%%%%%%%%%%%%%%%%%%%%%%%%%%%%%%%%%%%%%%%%%%%%%%%%%%%%%%%%%%%%%%%%%%%%%%%%%%%%%%%%%%%%%%%%%%%%%%%%%%%%%%%%%%%%%%%%%%%%%%%%%%%%%%%%%%%%%%%%%%%%%%%%%%%%%%%%%%%%%%%%%%%%%%%%%%%%%%
%%%%%%%%%%%%%%%%%%%%%%%%%%%%%%%%%%%%%%%%%%%%%%%%%%%%%%%%%%%%%

\bibliography{references}

\end{sloppypar}

\end{document}